\documentclass[final,onefignum,onetabnum]{siamonline190516}

\usepackage{lipsum}
\usepackage{amsfonts}
\usepackage{epstopdf}
\ifpdf
  \DeclareGraphicsExtensions{.eps,.pdf,.png,.jpg}
\else
  \DeclareGraphicsExtensions{.eps}
\fi

\usepackage{datetime}
\newdateformat{monthyeardate}{%
  \monthname[\THEMONTH], \THEYEAR}

\usepackage{academicons}

\usepackage{xcolor}

\usepackage{amsmath}
\allowdisplaybreaks
\usepackage{amssymb}
\usepackage{commath}
\usepackage{mathtools}
\usepackage{bbm}
\usepackage{bm}

\usepackage{color}
\usepackage{graphicx}
\usepackage[small]{caption}
\usepackage{subcaption}

\usepackage{relsize}
\usepackage{adjustbox}
\usepackage{algorithm}
\usepackage[noend]{algpseudocode}
\usepackage{booktabs}
\usepackage{tikz}
\usepackage{pifont}

\usepackage{enumitem}
\setlist[enumerate]{leftmargin=.5in}
\setlist[itemize]{leftmargin=.5in}

\newsiamremark{remark}{Remark}
\newsiamremark{example}{Example}
\newsiamremark{hypothesis}{Hypothesis}
\crefname{hypothesis}{Hypothesis}{Hypotheses}
\newsiamthm{claim}{Claim}

\headers{Stabilizing Global RBF Methods}{Jan Glaubitz and Anne Gelb}

\title{Stabilizing Radial Basis Function Methods for Conservation Laws Using Weakly Enforced Boundary Conditions\thanks{\monthyeardate\today \corresponding{Jan Glaubitz ( \email{Jan.Glaubitz@Dartmouth.edu}
)} 
\funding{This work is partially supported by the German Research Foundation
(DFG, Deutsche Forschungsgemeinschaft) \#GL 927/1-1 (Glaubitz),
AFOSR \#F9550-18-1-0316 (Glaubitz and Gelb),  NSF-DMS \#1502640, NSF-DMS \#1912685, and ONR \#N00014-20-1-2595  (Gelb).}}}

\author{Jan Glaubitz\thanks{Department of Mathematics, Dartmouth College, Hanover, NH 03755, USA}
\and 
Anne Gelb\footnotemark[2] 
}

\usepackage{amsopn}

\usepackage[normalem]{ulem}

\renewcommand{\d}{\mathrm{d}} 
\newcommand{\intd}{\, \mathrm{d}}
\newcommand{\N}{\mathbb{N}}

\newcommand{\R}{\mathbb{R}} 
\newcommand{\fnum}{f^{\mathrm{num}}}
\ifpdf
\hypersetup{
  pdftitle={Glaubitz2021Stabilizing},
  pdfauthor={Jan Glaubitz}
}
\fi

\begin{document}

\maketitle

\begin{abstract}
	It is well understood that boundary conditions (BCs) may cause global radial basis function (RBF) methods to become unstable for hyperbolic conservation laws (CLs). 
Here we investigate this phenomenon and identify the strong enforcement of BCs as the mechanism triggering such stability issues. 
Based on this observation we propose a technique to weakly enforce  BCs in RBF methods. 
In the case of hyperbolic CLs, this is achieved by carefully building RBF methods from the weak form of the CL, rather than the typically enforced strong form. 
Furthermore, we demonstrate that global RBF methods may violate conservation, yielding physically unreasonable solutions when the approximation does not take into account these considerations. 
Numerical experiments validate our theoretical results.

\end{abstract}



\section{Introduction}
\label{sec:introduction}

RBFs have become powerful tools in multivariate interpolation and 
approximation theory, since they are easy to implement, allow arbitrary scattered data, 
and can be spectrally accurate. 
They are also often used to solve numerical partial differential equations (PDEs) 
\cite{kansa1990multiquadrics,fasshauer1996solving,iske1996structure,hon1998efficient,kansa2000circumventing,larsson2003numerical,platte2004computing,hesthaven2019entropy,hesthaven2020two}. 
In this regard, although RBFs are considered to be a viable alternative to  traditional methods such as finite 
difference (FD), finite element (FE) and spectral methods, 
investigations into their stability  are still underdeveloped and/or unsatisfactory.  
For instance, $L^2$ (energy) stability has not been thoroughly studied. 
Moreover, for time-dependent PDEs, differentiation matrices for RBF methods often have  
eigenvalues with positive real parts, \cite{platte2006eigenvalue,sarra2011numerical}.  
Hence due to rounding errors RBFs can become increasingly unstable in time unless a dissipative time integration method \cite{platte2006eigenvalue,sarra2011numerical,martel2016stability}, artificial dissipation 
\cite{flyer2016enhancing,glaubitz2016artificial,ranocha2018stability,glaubitz2019smooth,glaubitz2019analysis}, 
or some other stabilizing technique 
\cite{scarnati2018using,gelb2019numerical,glaubitz2019high,hesthaven2008filtering,glaubitz2018application,glaubitz2019shock,glaubitz2020shock,don2016hybrid}, 
is employed.
Such stabilizing techniques often result in reduced accuracy, however, \cite{kansa2000circumventing,platte2004computing,schaback1995error}.  

This investigation seeks to increase the understanding of the stability requirements for RBF methods, especially as they relate to  hyperbolic conservation laws (CLs).  In one dimension, we therefore consider
\begin{equation}\label{eq:cl}
  u_t + f(u)_x = 0, \quad x \in \Omega = [a,b] \subset \R, \ t > 0, 
\end{equation} 
equipped with an appropriate initial condition (IC) {$u(0,x) = u_0(x)$} and BCs ${u(t,a) = g_L(t)}$, ${u(t,b) = g_R(t)}$. 
In \cite{martel2016stability}, eigenvalue analysis was used to show that in order to guarantee stability for the usual RBF methods, that is those using conditionally positive definite kernels, no BCs could be imposed on the problem.  
We note that the analysis was restricted  to  scalar linear advection, i.\,e.\ $f(u) = u$ in \eqref{eq:cl}. 
Starting from these results, this investigation pinpoints the root of stability issues not to be the existence of BCs, but rather how they are implemented within the RBF framework.
In particular we demonstrate that the BCs should be {\em weakly} enforced.  This is consistent with stable boundary treatment in FD methods  \cite{kreiss1989initial,gustafsson1995time,gustafsson2007high,svard2014review,fernandez2014review}, as well as  FE  \cite{huynh2007flux,vincent2011new,jameson2012non,ranocha2016summation,abgrall2019analysis,abgrall2019analysis_2} and spectral \cite{hesthaven2000spectral} methods. 

Our analysis involves using the weak form to solve \eqref{eq:cl} given by (see e.\,g.\ \cite{randall1992numerical})
\begin{equation}\label{eq:cl-weak}
  \int_\Omega u_t v \intd x - \int_\Omega f(u) v_x \intd x + f(u)v\big|_{\partial \Omega} = 0, 
  \quad t > 0,
\end{equation}
with test function $v \in C^1(\Omega)$. 
Recall that \eqref{eq:cl-weak} is constructed from \eqref{eq:cl} by multiplying each term by  $v$, integrating over $\Omega$, and applying integration by parts. 
Observe that for \eqref{eq:cl-weak} less regularity is required for the solution $u$. 
This is important since even for smooth initial conditions solutions of \eqref{eq:cl} can develop jump discontinuities \cite{lax1973hyperbolic,dafermos2005hyperbolic}. 
Thus by using \eqref{eq:cl-weak} we permit the more general class of weak solutions, where \eqref{eq:cl} is satisfied in the sense of distribution theory, see \cite{lax1973hyperbolic,dafermos2005hyperbolic}. 
To distinguish the physically reasonable weak solution from all of the other possible weak 
solutions, \eqref{eq:cl} is augmented with an additional entropy condition 
\begin{equation}\label{eq:entropy-cond}
  U(u)_t + F(u)_x \leq 0.
\end{equation}
Here $U$ is an entropy function and $F$ is a corresponding entropy flux satisfying ${U'f'=F'}$. 
A strict inequality in \eqref{eq:entropy-cond} reflects the presence of a physically 
reasonable shock wave. 
For scalar conservation laws in one dimension, the square entropy ${U(u) = \frac{1}{2} u^2}$ 
is often a valid entropy function. 
In this case, from the entropy inequality \eqref{eq:entropy-cond}, we immediately get 
\begin{equation}\label{eq:L2-stab}
  \frac{\d}{\d t} \norm{u}^2_{L^2}
    = 2 \int_\Omega U(u)_t \intd x 
    \leq - 2 F(u)\big|_{\partial \Omega}
\end{equation}
for entropy solutions of \eqref{eq:cl}. 
In particular, the entropy should not increase over time for an isolated physical system, and a physically 
reasonable weak solution of \eqref{eq:cl} should therefore satisfy 
\begin{equation}\label{eq:L2-stab-periodic}
  \frac{\d}{\d t} \norm{u}^2_{L^2} \leq 0 
\end{equation}
for periodic BCs.
We refer to \eqref{eq:L2-stab} as \emph{$L^2$} or \emph{energy stability}. 
Together with the property of \emph{conservation}, given by
\begin{equation}\label{eq:conservation}
  \frac{\d}{\d t} \int_\Omega u \intd x = - f(u)\big|_{\partial \Omega},
\end{equation} 
energy stability often is considered an important design criteria for a numerical method to produce 
physically reasonable solutions. 

In what follows we show that it can be beneficial to build RBF methods from the weak form \eqref{eq:cl-weak} instead of the strong form \eqref{eq:cl}, which is the usual 
approach. 
We prove that RBF methods based on the weak form,  which we will refer to as \emph{weak RBF methods}, are conservative as long as constants are included in the RBF approximation, which will be explained in \S \ref{sec:preliminaries}. 
They are also energy stable when appropriate numerical fluxes are used for the (weak) treatment of BCs. 
In contrast, we also demonstrate that usual RBF methods based on the strong form, which we will refer to as \emph{strong RBF methods}, violate conservation as well as energy stability and might produce physically unreasonable solutions.
Our approach is closely related to the idea behind discontinuous Galerkin (DG) methods 
\cite{cockburn1991runge,cockburn1989tvb,cockburn1989tvb2,cockburn1990runge,cockburn1998runge,hesthaven2007nodal}. 
For these, a resembling but different energy stability analysis was performed in \cite{jiang1994cell}. 
Details on energy stability for DG methods and related schemes can be found in,  e.\,g., \cite{gassner2013skew,svard2014review,chen2017entropy,ranocha2018stability,offner2018stability,glaubitz2019smooth,glaubitz2020stableDG} and references therein. 
To the best of our knowledge, none of these investigations prove energy stability properties for RBF methods for hyperbolic CLs.

The rest of this work is organized as follows.   
In Section \ref{sec:preliminaries}, we collect all necessary preliminaries on RBF approximations. 
The heart of this investigation is Section \ref{sec:stable_RBFs},  where we prove that the weak RBF method for CLs is conservative and energy stable. 
We further describe two different realizations of the resulting weak RBF methods, the {\em weak RBF analytical method} and the more efficient {\em weak RBF collocation method}. 
{In Section \ref{sec:connections} we provide a comparison of the weak RBF method with some commonly used techniques.}
Section \ref{sec:numerical} compares numerical results for our new method with the traditional RBF method, and some concluding remarks are offered in  Section \ref{sec:summary}.

The MATLAB code corresponding to this manuscript can be found at \cite{glaubitz2020weakRBFcode}.
\section{Preliminaries}
\label{sec:preliminaries}

This section collects all necessary concepts and results regarding RBF approximations. 
More details may be found in the survey articles 
\cite{schaback1995error,schaback1995multivariate,schaback2005multivariate}.

\subsection{Method of Lines} 

In this investigation we consider only spatial discretization of the hyperbolic CL (\ref{eq:cl}), so that the problem remains continuous in time.  
The resulting system of ordinary differential equations (ODEs), often referred to as the semi-discrete formulation, is given by
\begin{equation}\label{eq:ODE}
  \frac{\d}{\d t} u = L(u),  
\end{equation} 
where $L(u)$ is a discretization of the spatial operator. 
This approach, i.\,e.\ where time dependent PDEs are reduced to a system of ODEs, is often called the 
\emph{method of lines}, see \cite[Chapter 10.4]{leveque2002finite}. 
Time integration techniques used for solving (\ref{eq:ODE})  will be further discussed in Section \ref{sub:time-int}.

\subsection{{RBF} Approximations} 

We now consider approximations of a function $u:\Omega \to \R$ with $\Omega \subset \R^d$ by \emph{{RBF} interpolants} 
\begin{equation}\label{eq:u-RBF}
  u_N(\boldsymbol{x}) = \sum_{n=1}^N \alpha_n \phi\left( \varepsilon \norm{\boldsymbol{x}-\mathbf{x}_n} \right), 
\end{equation} 
where $\phi:\R \to \R$ is a \emph{basis function (kernel)} and the coefficients $\alpha_k$ are 
calculated such that the interpolation condition 
\begin{equation}\label{eq:interpol-cond}
  u_N(\mathbf{x}_n) = u(\mathbf{x}_n), \quad n=1,\dots,N,   
\end{equation} 
holds. 
The interpolation points $\mathbf{x}_n \in \Omega$ are called \emph{centers} and $\varepsilon > 0$ is the 
\emph{shape parameter}. 
The interpolation condition \eqref{eq:interpol-cond} yields a system of linear equations,  
\begin{equation}\label{eq:LS-RBF-interpol}
  \underbrace{\begin{pmatrix} 
    \phi\left( \varepsilon \norm{\mathbf{x}_1-\mathbf{x}_1} \right) & \dots & \phi\left( \varepsilon \norm{\mathbf{x}_1-\mathbf{x}_N} 
\right) \\ 
    \vdots & & \vdots \\ 
    \phi\left( \varepsilon \norm{\mathbf{x}_N-\mathbf{x}_1} \right) & \dots & \phi\left( \varepsilon \norm{\mathbf{x}_N-\mathbf{x}_N} 
\right)
  \end{pmatrix}}_{=: \Phi}
  \underbrace{
  \begin{pmatrix} 
    \alpha_1 \\ \vdots \\ \alpha_N
  \end{pmatrix}
  }_{=: \boldsymbol{\alpha}}
  = 
  \underbrace{
  \begin{pmatrix} 
    u(\mathbf{x}_1) \\ \vdots \\ u(\mathbf{x}_N)
  \end{pmatrix}
   }_{=: \mathbf{u}},
\end{equation}
which can be solved for the vector of coefficients $\boldsymbol{\alpha} \in \R^N$ if the matrix $\Phi$ is invertible. 
Popular examples for basis functions (kernels) $\phi$ are
\begin{alignat}{2}
  	\phi(r) & = e^{-r^2} \qquad && \text{(Gaussian)}, \label{eq:kernel_G} \\ 
  	\phi(r) & = \sqrt{1+r^2} \qquad && \text{(multiquadric)}, \label{eq:kernel_MQ} \\ 
  	\phi(r) & = \frac{1}{(1+r^2)} \qquad && \text{(inverse quadratic)}, \label{eq:kernel_IQ} \\ 
  	\phi(r) & = 
		\left\{
		\begin{array}{c c l} 
			r^{k} & ; & k \in 2\N+1, \\ 
			r^k \log r & ; & k \in 2\N,
		\end{array} 
		\right. 
		\qquad && \text{(polyharmonic splines)}, \label{eq:kernel_PHS}
\end{alignat}
More details may be found in 
\cite{schaback1995creating,buhmann2003radial,wendland2004scattered,schaback2005multivariate,platte2005polynomials,fornberg2007runge} 
and references therein.

\subsection{Stability of RBF Methods for Time-Dependent Problems} 
\label{sub:stability}

Experience suggests that RBF approximations will produce discretizations that are 
unstable in time unless highly dissipative time stepping is used. 
It was shown in \cite{platte2006eigenvalue} that under a variety of conditions,
differentiation matrices obtained with RBF collocation methods have eigenvalues 
with positive real parts. 
In particular, this was demonstrated for a simple one-dimensional linear advection equation, 
suggesting its unsuitability for nonlinear hyperbolic CLs. 
A related observation was made in \cite{fornberg2002observations}, where 
it was proposed that one source of instability might be inaccuracy of RBF approximations 
near boundaries. 
On the flip side it was also proved in \cite{platte2006eigenvalue} that RBF 
collocation methods are time-stable (in the sense of eigenvalues for linear problems) for all 
conditionally positive definite RBFs and node distributions when \emph{no BCs are needed}. 
Hence while RBFs perform well in periodic domains, such as circles or unit spheres, they are evidently not suitable in applications where periodicity of the computational domain cannot be assumed.
In Section \ref{sec:numerical} we will also demonstrate that conservation and energy stability are both violated by usual RBF methods when applied to hyperbolic CLs, possibly leading to physically irrelevant solutions.

\subsection{RBF Approximations With Polynomials}

RBF interpolants \eqref{eq:u-RBF} are often modified to include polynomials along with matching constraints on 
the expansion coefficients, 
\cite{schaback1995error,buhmann2000radial,flyer2016role,flyer2016enhancing}. 
For example, for ${\Omega \subset \R^d}$, let us define $\{p_k\}_{k=1}^K$ as a basis for the space of polynomials of degree at most $P-1$ in $d$ variables, denoted by $\mathbb{P}_{P-1}(\R^d)$, where
$K = \binom{P-1+d}{d}$. 
The resulting RBF interpolants for  polynomials of degree up to $P-1$ are then
\begin{equation}\label{eq:u-RBF-polynomial}
  u_N(\boldsymbol{x}) 
    = \sum_{n=1}^N \alpha_n \phi\left( \varepsilon \norm{\boldsymbol{x}-\mathbf{x}_k} \right) 
    + \sum_{k=1}^K \beta_k p_k(\boldsymbol{x}) 
\end{equation} 
with constraints 
\begin{equation}\label{eq:constraints}
  \sum_{n=1}^N \alpha_n p_k(\mathbf{x}_n) = 0, \quad k=1,\dots,K.
\end{equation} 
Let us also define 
\begin{equation}\label{eq:P} 
  P = 
  \begin{pmatrix} 
    p_1\left( \mathbf{x}_1 \right) & \dots & p_1\left( \mathbf{x}_N \right) \\ 
    \vdots & & \vdots \\ 
    p_K\left( \mathbf{x}_1 \right) & \dots & p_K\left( \mathbf{x}_N \right)
  \end{pmatrix}, \quad 
  \boldsymbol{\beta} = 
  \begin{pmatrix} 
    \beta_1 \\ 
    \vdots \\ 
    \beta_K
  \end{pmatrix}.
\end{equation}
Then, given the interpolation condition \eqref{eq:interpol-cond}, the counterpart to \eqref{eq:LS-RBF-interpol} is 
\begin{equation}\label{eq:LS-RBF-polynomial} 
  \underbrace{\begin{pmatrix} 
    \Phi & P^T \\ 
    P & 0
  \end{pmatrix}}_{=: V}
  \begin{pmatrix} 
    \boldsymbol{\alpha} \\ \boldsymbol{\beta}
  \end{pmatrix}
  = 
  \begin{pmatrix} 
    \mathbf{u} \\ \mathbf{0}
  \end{pmatrix}. 
\end{equation}
 
There are various reasons for including polynomials in RBF interpolants
\cite{schaback1995error,buhmann2000radial,flyer2016role,flyer2016enhancing}:

\begin{enumerate}
  \item 
  Polynomial terms can ensure that \eqref{eq:LS-RBF-polynomial} is uniquely solvable when working with conditionally positive definite basis functions (kernels), assuming the set of centers $\{\mathbf{x}_k\}_{k=1}^N$ is $\mathbb{P}_{P-1}(\R^q)$-unisolvent. 
  See for instance \cite[Chapter 7]{fasshauer2007meshfree}.
  
  \item 
  Numerical tests demonstrate that including a constant improves the accuracy of derivative approximations. 
  In particular, adding a constant avoids oscillatory representations of constant functions. 
  
  \item 
  Including polynomial terms of low order can also improve the accuracy of RBF interpolants near domain boundaries 
due to regularizing the far-field growth of RBF interpolants \cite{fornberg2002observations}. 
\end{enumerate}

\begin{figure}[!htb]
  \centering
  \begin{subfigure}[b]{0.45\textwidth}
    \includegraphics[width=\textwidth]{%
      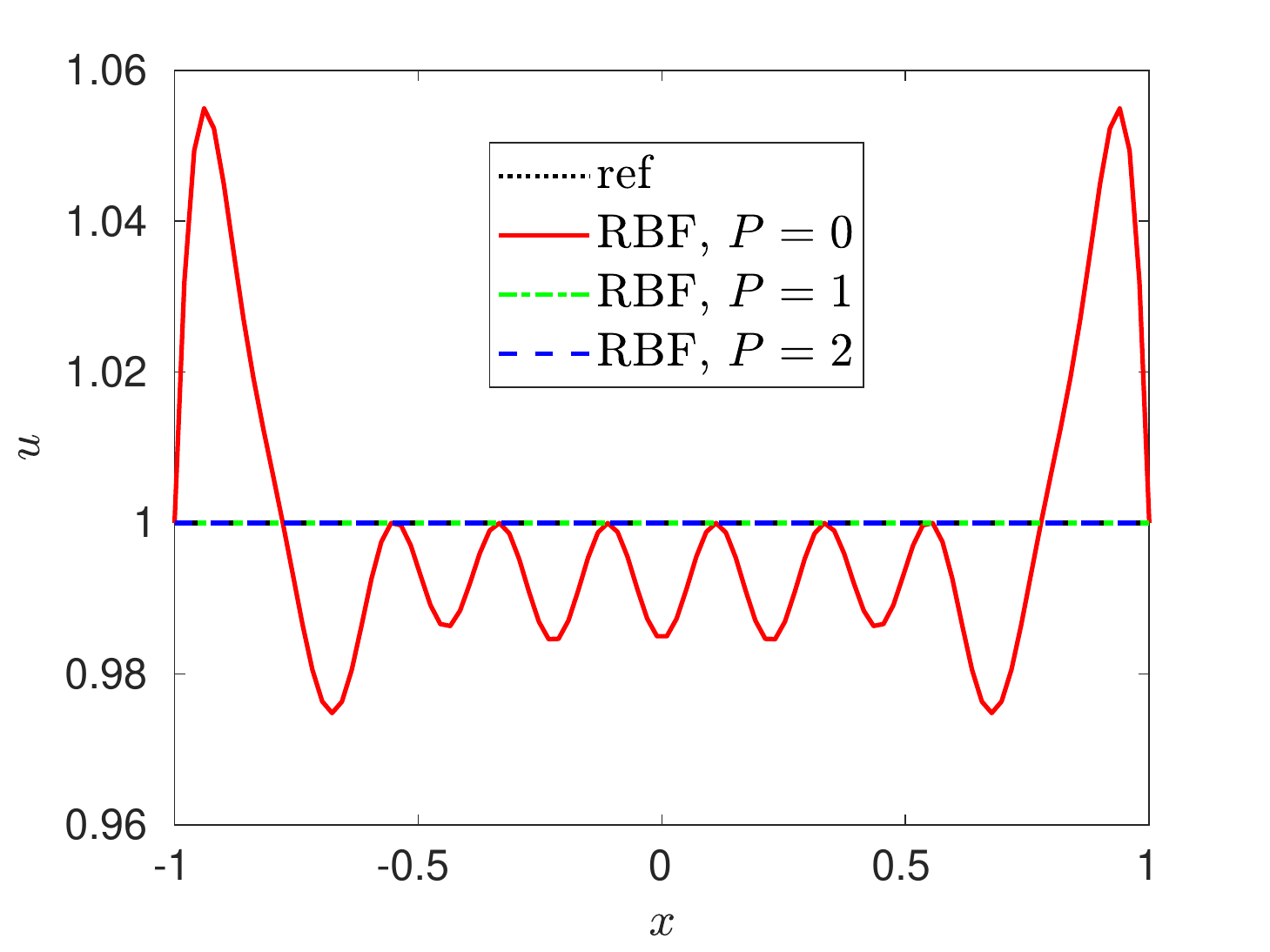}
    \caption{RBF approx. ($\varepsilon=6$) of $u(x)=1$}
    \label{fig:RBF_approx_const}
  \end{subfigure}%
  ~ 
  \begin{subfigure}[b]{0.45\textwidth}
    \includegraphics[width=\textwidth]{%
      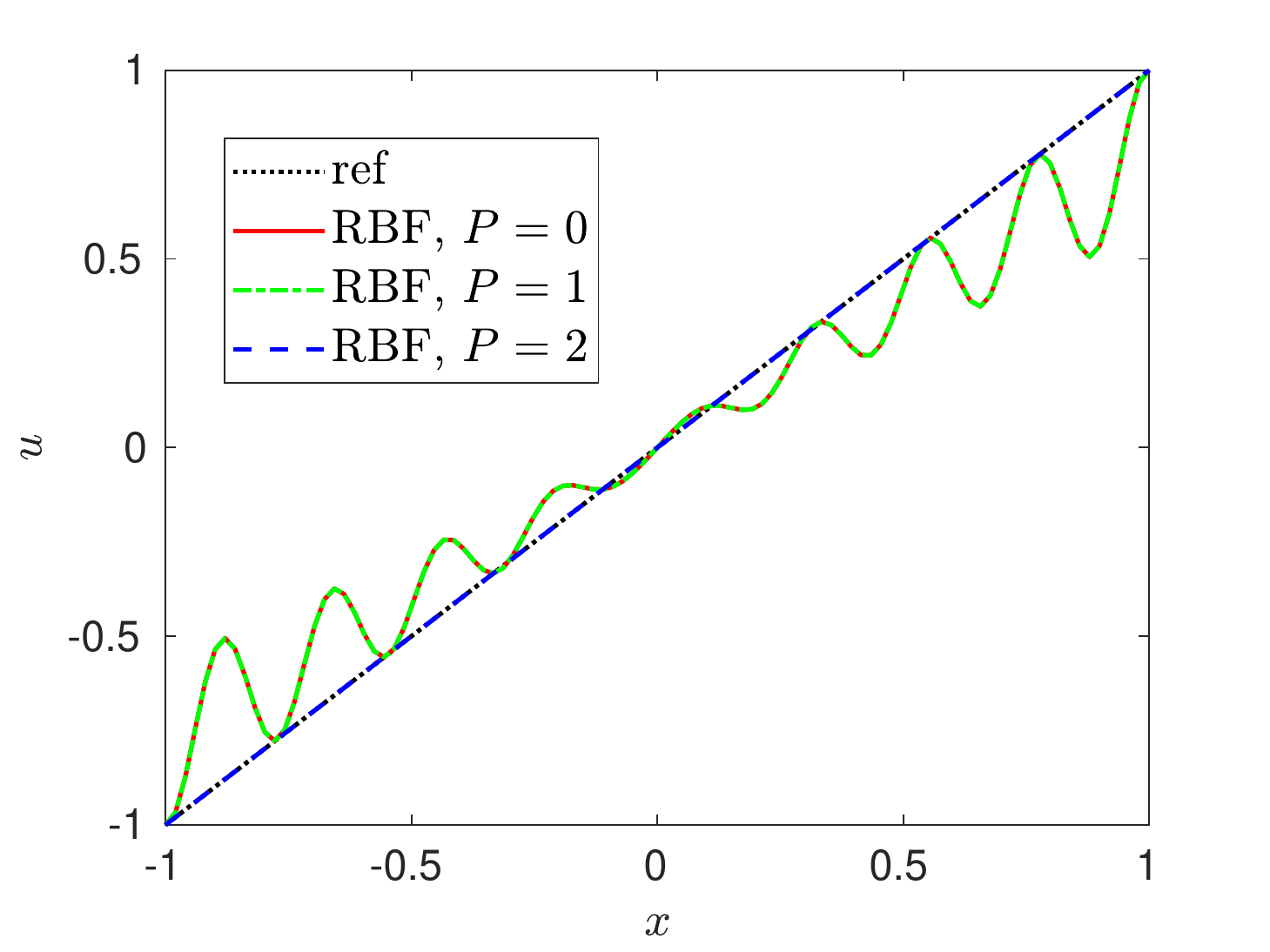}
    \caption{RBF approx. ($\varepsilon=10$) of $u(x) = x$}
    \label{fig:RBF_approx_linear}
  \end{subfigure}
  \caption{RBF approximations including polynomials up to different degrees. 
  In both cases Gaussian kernels were used. 
    }
  \label{fig:RBF_approx}
\end{figure}

For our purposes, the main advantage in including polynomials in the RBF interpolants is that the constraints in
\eqref{eq:constraints} enforce the RBF interpolants \eqref{eq:u-RBF-polynomial} to reproduce polynomials up to degree 
$P-1$: 
\begin{equation*}
  u_N = u \quad \forall u \in \mathbb{P}_{P-1}(\R^d)
\end{equation*} 
For example, Figure \ref{fig:RBF_approx} demonstrates in one dimension ($d=1$) that constant functions can be reconstructed exactly by RBF interpolants for $P \geq 1$.
This property will be crucial to prove conservation for the stable RBF methods proposed in Section \ref{sec:stable_RBFs}. 

\begin{remark} 
	We stress that  the above discussion is specific to global RBFs.   Polynomials play a different role in local RBF (RBF-FD) methods, \cite{flyer2016role}
\end{remark}

\section{Energy Stable RBF Methods}
\label{sec:stable_RBFs}

RBF methods typically use collocation to discretize \eqref{eq:cl}. 
That is, $u$ and $f$ are both approximated by RBF interpolants with respect to the same set of centers $\mathbf{x}_n$, $n=1,\dots,N$. 
As discussed in Section \ref{sub:stability}, this yields unstable methods in the presence of BCs. 
Here, however, we prove that stability as well as conservation can be ensured if RBF methods are built from the weak form. 
For ease of presentation, we perform our analysis in one dimension ($d = 1$). 
As will be demonstrated in Section \ref{sub:2d}, the method can be implemented  in higher dimensions.  No attempt has been made to prove stability for $d > 1$, however.

In one dimension, the weak form \eqref{eq:cl-weak} is equivalent to
\begin{equation}\label{eq:weak-form}
  \int_{\Omega} u_t v \intd x - \int_{\Omega} f(u) v_x \intd x 
    + f(u(t,b))v(b) - f(u(t,a))v(a) = 0
\end{equation} 
with $v \in C^1(\Omega)$ and $t > 0$.
%
In what follows we describe two different RBF methods built from \eqref{eq:weak-form}. 
In both cases the solution $u$ is approximated by an RBF interpolant \eqref{eq:u-RBF-polynomial}, 
which as we noted earlier can include polynomials.

The method described in Section \ref{sub:ana-RBF-method} uses the {\em analytical}
flux function $f$ applied to the RBF interpolant $u_N$.  As a consequence,
the resulting approximation $f(u_n) \approx f(u)$ still satisfies the interpolation 
condition but is no longer an RBF approximation.
By contrast, the technique described in Section \ref{sub:coll-RBF-method} utilizes the idea of collocation, where $u$ 
and the flux $f(u)$ are both approximated by RBF interpolants.

\subsection{Weak RBF Analytical Methods}
\label{sub:ana-RBF-method}

Let $u$ and $v$ in the weak form \eqref{eq:weak-form} be replaced by RBF approximations 
${u_N, v_N \in V_{N,P}}$ with 
\begin{equation}\label{eq:RBFspace} 
\resizebox{.9\textwidth}{!}{$\displaystyle 
  V_{N,P} := \left\{ \sum_{n=1}^N \alpha_n \phi\left( \varepsilon \norm{x-x_n} \right) 
    + \sum_{k=1}^K \beta_k p_k(x) \mid 
    \boldsymbol{\alpha} \in \R^N, \boldsymbol{\beta} \in \R^{K}, \text{ and } \eqref{eq:constraints} \text{ holds} 
\right\}, 
$}
\end{equation} 
where $K = \binom{P-1+1}{1} = P$. 
Note that for $P= 0$ no polynomials are included in the RBF interpolant and the approximation space reduces to 
\begin{equation}
  V_{N,0} = \operatorname{span}\left\{ \phi( \varepsilon \norm{x-x_n} ) \mid n=1,\dots,N \right\}.
\end{equation}
Next observe that while one or both BCs may be given as part of \eqref{eq:cl}, 
i.\,e.\ ${u(t,a) = g_L(t)}$ and ${u(t,b) = g_R(t)}$, 
it is also possible to assign these values with the RBF approximations evaluated there as
\begin{equation}
	u(t,a) = u_R := u_N(b), \quad u(t,b) = u_L := u_N(a).
\end{equation}
Hence to ensure well-defined boundary terms, we compute a single valued numerical flux at the boundaries as
\begin{equation}  
  \fnum_{L} = \fnum\left( g_L(t), u_L \right), \quad 
  \fnum_{R} = \fnum\left( u_R, g_R(t) \right),
\end{equation} 
and therefore enforce the BCs in a weak sense.
The numerical flux is chosen to be (i) consistent, that is we require $\fnum(u,u) = f(u)$; (ii) Lipschitz continuous; and (iii) monotone, meaning that $\fnum$ is nondecreasing in the first argument and nonincreasing in the second argument. 
Examples of commonly used numerical fluxes can be found in \cite{cockburn1989tvb,toro2013riemann}. 
We are now ready to define the \emph{weak RBF analytical method} as 

\begin{definition}[Weak RBF analytical method]\label{def:weakRBFanalytic}
 Determine $u_N \in V_{N,P}$ such that all $v_N \in V_{N,P}$ satisfies
  \begin{equation}\label{eq:anal-RBF} 
    \int_{\Omega} (u_N)_t v_N \intd x 
      - \int_{\Omega} f(u_N) (v_N)_x \intd x 
      + \left( \fnum_R v_R - \fnum_L v_L \right) = 0,
  \end{equation} 
where $v_L$ and $v_R$ respectively denote $v_N(a)$ and $v_N(b)$.
\end{definition}

Note that in \eqref{eq:anal-RBF} all integrals as well as the flux $f(u_N)$ are assumed to be 
evaluated exactly. 
Next we consider the properties of the weak RBF analytical method (\ref{eq:anal-RBF}) for the 
one-dimensional CL \eqref{eq:cl}.

\subsubsection{Conservation} 

The rate of change of the total amount of the conserved variable $u$ is given by \eqref{eq:conservation}, which establishes that the total amount of change in $u$  is due to the flux across the domain boundaries. 
For periodic BCs  conservation implies that
\begin{equation}\label{eq:conserve_u-periodic}
  \frac{\d}{\d t} \int_{\Omega} u \intd x = 0.
\end{equation} 
The highly celebrated Lax Wendroff theorem states that if a conservative numerical scheme converges, then it will converge toward a weak solution, \cite{randall1992numerical}. 
To prove conservation for \eqref{eq:anal-RBF}, we choose $P\geq1$ in order to include polynomials of degree 
$P-1$ in the approximation space $V_{N,P}$ defined by \eqref{eq:RBFspace}. 
Thus $1 \in V_{N,P}$, and since \eqref{eq:anal-RBF} holds for $v_N = 1$, we have
\begin{equation}
  \frac{\d}{\d t} \int_{\Omega} u_N \intd x  = \int_{\Omega} (u_N)_t \intd x
    = - \left( \fnum_R - \fnum_L \right), 
\end{equation}
which is the discrete counterpart to \eqref{eq:conservation}. 
Note that for periodic BCs, the numerical fluxes are given by ${\fnum_L = \fnum(u_R,u_L)}$ and 
${\fnum_R = \fnum(u_R,u_L)}$, yielding 
\begin{equation}
  \frac{\d}{\d t} \int_{\Omega} u_N \intd x = 0. 
\end{equation}
Observe that for periodic BCs, conservation of the continuous equation \eqref{eq:cl} is exact.

\subsubsection{Energy Stability}

Recall that \eqref{eq:L2-stab} implies that the rate of change of the squared $L^2$ norm is given by 
\begin{equation}  
	\frac{\d}{\d t} \norm{u}_{L^2}^2 = 2 \int_{\Omega} u_t u \intd x.
\end{equation}
Hence by choosing $v_N = u_N$ in \eqref{eq:anal-RBF} we obtain 
\begin{equation}
\begin{aligned}
  	 \frac{1}{2} \frac{\d}{\d t} \norm{u_N}_{L^2}^2 
    		= & \int_{\Omega} f(u_N) (u_N)_x \intd x 
      	- \left( \fnum_R u_R - \fnum_L u_L \right) \\
    		= & - \int_{\Omega} f(u_N)_x u_N \intd x 
      	+ \left( f(u_R) u_R - f(u_L) u_L \right) 
       	- \left( \fnum_R u_R - \fnum_L u_L \right),  
\end{aligned}
\end{equation}
with the second equality resulting from applying integration by parts. 
Observe that for the square entropy  $U(u)=\frac{u^2}{2}$ with corresponding entropy flux $F(u)$ 
satisfying $U'f'=F'$ we have 
\begin{equation}
  F(u)_x = F'(u) u_x  = u f'(u) u_x = u f(u)_x,
\label{eq:fluxentropy}
\end{equation} 
yielding 
\begin{equation}\label{eq:eq1} 
\begin{aligned}
  \frac{1}{2} \frac{\d}{\d t} \norm{u_N}_{L^2}^2 
    = & - \left( F(u_R) - F(u_L) \right) 
      + \left( f(u_R) u_R - f(u_L) u_L \right) \\ 
       & - \left( \fnum_R u_R - \fnum_L u_L \right). 
\end{aligned}
\end{equation} 
Further, by defining
\begin{equation}  
  {\gamma}(u) := \int^u f(v) \intd v, 
\end{equation}
the entropy flux $F(u)$ can be written as (see \cite{jiang1994cell})
\begin{equation}\label{eq:gamma}  
  F(u) 
    = \int^u f'(v) v \intd v 
    = f(u)u - \int^u f(v) \intd v
    = f(u)u - {\gamma}(u) 
\end{equation}
so that \eqref{eq:eq1} becomes 
\begin{equation}
\begin{aligned}
  \frac{1}{2} \frac{\d}{\d t} \norm{u_N}_{L^2}^2 
    = & \left( \gamma(u_R) - \gamma(u_L) \right) 
      - \left( \fnum_R u_R - \fnum_L u_L \right) \\ 
    = & \left( \gamma(u_R) - \gamma(g_R) \right) 
      - \left( \gamma(u_L) - \gamma(g_L) \right) \\ 
    & + \left( \gamma(g_R) - \gamma(g_L) \right) 
      - \left( \fnum_R u_R - \fnum_L u_L \right),
\end{aligned} 
\end{equation} 
where $g_L$ and $g_R$ are the BCs given as part of \eqref{eq:cl}. 
By the mean value theorem, there exists a $u_L^*$ between $u_L$ and $g_L$ as well as a $u_R^*$ between $u_R$ and $g_R$ 
such that 
\begin{equation}
\begin{aligned}
  \gamma(u_L) - \gamma(g_L)  &= \left( u_L - g_L \right) f(u_L^*), \\
  \gamma(u_R) - \gamma(g_R)  &= \left( u_R - g_R \right) f(u_R^*).
\end{aligned} 
\end{equation}
In this case we have 
\begin{equation}
\begin{aligned}
  \frac{1}{2} \frac{\d}{\d t} \norm{u_N}_{L^2}^2 
   = & \left( u_R - g_R \right) f(u_R^*) 
      - \left( u_L - g_L \right) f(u_L^*) 
       + \left( \gamma(g_R) - \gamma(g_L) \right) \\ 
   & - \left( \fnum_R u_R - \fnum_L u_L \right) \\ 
   = & \left( g_R - u_R \right) \left( \fnum_R - f(u_R^*) \right) 
      + \left( u_L - g_L \right) \left( \fnum_L - f(u_L^*) \right) \\
    & + \left( \gamma(g_R) - \gamma(g_L) \right) 
      - \left( g_R \fnum_R - g_L \fnum_L \right),
\end{aligned} 
\end{equation}
where  the numerical fluxes are given respectively by 
\begin{equation}
	\fnum_L = \fnum(g_L,u_L), \quad
	\fnum_R = \fnum(u_R,g_R).
\end{equation}
Thus, by employing an \emph{E-Flux} (see \cite{osher1984riemann}) so that
\begin{equation} 
  (b - a) \left( \fnum(a,b) - f(u) \right) \leq 0 
\end{equation}
for all $u$ between $a$ and $b$, we have 
\begin{equation}
  \frac{1}{2} \frac{\d}{\d t} \norm{u_N}_{L^2}^2 
    \leq \left( \gamma(g_R) - \gamma(g_L) \right) 
      - \left( g_R \fnum_R - g_L \fnum_L \right)
\end{equation}
Finally, utilizing \eqref{eq:gamma} results in 
\begin{equation}
  \frac{\d}{\d t} \norm{u_N}_{L^2}^2 
    \leq -2F(u_N)\big|_{\partial \Omega} 
      + 2 g_R \left( f(g_R) - \fnum_R \right) 
      - 2 g_L \left( f(g_L) - \fnum_L \right),
\end{equation} 
which is consistent with \eqref{eq:L2-stab} since the numerical flux $\fnum$ is consistent with the flux $f$. 
In particular, the above inequality implies \eqref{eq:L2-stab-periodic} for periodic BCs. 
This yields a conservative and energy stable RBF method for general one dimensional scalar CLs.

\subsection{Weak RBF Collocation Methods} 
\label{sub:coll-RBF-method}

Depending on the nonlinearity of $f$, the exact evaluation of $f(u_N)$ and resulting 
integrals may be impractical or even impossible. 
We therefore  extend our analysis from Section \ref{sub:ana-RBF-method} to a collocation based alternative to the weak RBF analytic method given in Definition \ref{def:weakRBFanalytic}.
As before, we replace $u$ and $v$ with their RBF approximations $u_N,v_N \in V_{N,P}$ for $P \geq 1$.
In the collocation case, $f(u)$ is approximated using an RBF interpolant $f_N \in V_{N,P}$ such that
\begin{equation}  
  f_N(x_n) = f(u_N(x_n)), \quad n=1,\dots,N. 
\end{equation}
We can now proceed as in the weak RBF analytical method and define
\begin{definition}[Weak RBF collocation method]\label{def:weakRBFcollocation}
Find $u_N \in V_{N,P}$ such that
\begin{equation}\label{eq:collo-RBF} 
  \int_{\Omega} (u_N)_t v_N \intd x 
    - \int_{\Omega} f_N (v_N)_x \intd x 
    + \left( \fnum_R v_R - \fnum_L v_L \right) = 0
\end{equation} 
for all $v_N \in V_{N,P}$.
\end{definition}

\subsubsection{Conservation}

As in the weak RBF analytical case, conservation follows by including constants in the RBF interpolants, i.\,e.\ by 
choosing $P \geq 1$.

\subsubsection{Energy Stability}

For the weak RBF collocation method, we can only prove energy stability for the linear advection equation, given by  
\begin{equation}\label{eq:lin-adv}
  u_t + {\lambda} u_x = 0.
\end{equation} 
From (\ref{eq:fluxentropy}) we obtain the entropy flux  $F(u) = (\lambda/2) u^2$. 
Here we pick constant velocity $\lambda > 0$ and note that the case for $\lambda < 0$ can be treated analogously.
By choosing $v_N = u_N$ in \eqref{eq:collo-RBF}, we obtain
\begin{equation}
\begin{aligned}
  	\frac{1}{2} \frac{\d}{\d t} \norm{u_N}_{L^2}^2  
  		& = \int_{\Omega} (u_N)_t u_N \intd x\\
    		& = \lambda \int_{\Omega} u_N (u_N)_x \intd x 
      	- \left( \fnum_R u_R - \fnum_L u_L \right) \\ 
   		& = - \lambda  \int_{\Omega} (u_N)_x u_N \intd x
      	+ \lambda \left( u_R^2 - u_L^2 \right)
      	- \left( \fnum_R u_R - \fnum_L u_L \right),
\end{aligned}
\end{equation}
where we have used integration by parts. 
Summing up the second and third equations above yields
\begin{equation}\label{eq:weak-RBF-coll-L2}
  \frac{\d}{\d t} \norm{u_N}_{L^2}^2 
     = \lambda \left( u_R^2 - u_L^2 \right) - 2 \left( \fnum_R u_R - \fnum_L u_L \right),
\end{equation} 
which can be rewritten as 
\begin{equation}
  \frac{\d}{\d t} \norm{u_N}_{L^2}^2 
     = -2F(u_N)\big|_{\partial \Omega} 
      + 2\lambda u_R \left( u_R - \fnum_R \right) 
      - 2\lambda u_L \left( u_L - \fnum_L \right).
\end{equation}
By now employing a simple \emph{upwind flux}, ${\fnum(a,b) = \lambda a}$, we have 
\begin{equation}
\begin{aligned}
  \fnum_L & = \fnum(g_L,u_L) = \lambda g_L, \\ 
  \fnum_R & = \fnum(u_R,g_R) = \lambda u_R,
\end{aligned}
\end{equation}
and therefore 
\begin{equation}
  \frac{\d}{\d t} \norm{u_N}_{L^2}^2 
     = - 2F(u_N)\big|_{\partial \Omega} 
      - 2 \lambda u_L \left( u_L - g_L \right). 
\end{equation}
The above equation is consistent with \eqref{eq:L2-stab}. 
Note that for the linear advection equation \eqref{eq:lin-adv} no shock waves arise and the inequalities 
\eqref{eq:entropy-cond} and \eqref{eq:L2-stab} become equalities. 
Moreover, for periodic BCs, \eqref{eq:weak-RBF-coll-L2} reduces to 
\begin{equation}
\begin{aligned}
  \frac{\d}{\d t} \norm{u_N}_{L^2}^2 
     & = \lambda \left( u_R^2 - u_L^2 \right) - 2 \left( \fnum(u_R,u_L) u_R - \fnum(u_R,u_L) u_L \right) \\
     & = - \lambda u_R^2 + 2 \lambda u_L u_R - \lambda u_L^2 \\ 
     & = - \lambda (u_R-u_L)^2 \\ 
     & \leq 0.
\end{aligned}
\end{equation}

\begin{remark}
	Recall that for general CLs $u_t + f(u)_x = 0$, $L^2$ stability for the weak RBF analytical method in Definition \ref{def:weakRBFanalytic} was shown by utilizing the relation  
	\begin{equation}\label{eq:chain}
  		F(u_N)_x = (u_N)_x F'(u_N) = (u_N)_x U'(u_N) f'(u_N) = u_N f(u_N)_x,
	\end{equation} 
	for the square entropy $U(u) = \frac{u^2}{2}$.  For the weak RBF collocation method in Definition \ref{def:weakRBFcollocation}, $f(u_N)$ in (\ref{eq:chain}) is replaced by $f_N$ and the final equality does not hold.  
	Thus we are unable to prove energy stability for general nonlinear CLs.
\end{remark}

\subsection{Numerical Fluxes}
\label{sub:E-fluxes}

There are several options for choosing numerical fluxes that result in energy stable weak RBF methods for one-dimensional scalar CLs.  
Some examples include
\begin{enumerate}
  \item  {\bf Upwind flux}: For linear advection, $u_t + \lambda u_x = 0$, with constant velocity $\lambda \neq 0$, the general upwind 
flux, given by
\begin{equation}   
 	\fnum(a,b) = 
    \begin{cases}
      \lambda a &; \ \lambda > 0, \\ 
      \lambda b &; \ \lambda < 0,
    \end{cases},
\end{equation}
  yields energy stability for both the analytical and collocation forms.
  \item  {\bf E-Flux}: For the nonlinear case we can use an E-Flux as defined in \cite{osher1984riemann} 
(see also \cite{cockburn1989tvb} and references therein). 
  For example, the Godunov flux is given by
\begin{equation}   
    \fnum(a,b) = 
    \begin{cases}
      \min_{a \leq u \leq b} f(u) &; \ a \leq b, \\ 
      \max_{a {\geq} u {\geq} b} f(u) &; \ a > b.
    \end{cases}
\end{equation}
\end{enumerate}

\subsection{Time Integration}
\label{sub:time-int}

Once we obtain the spatial discretization for the hyperbolic CL using one of the above methods, we then solve the semi-discrete formulation in \eqref{eq:ODE}.
Popular choices of time integration methods include explicit total variation
diminishing (TVD) Runge--Kutta (RK) methods \cite{shu1988total,gottlieb1998total}, also known as 
 strong stability preserving (SSP) RK methods 
\cite{gottlieb2001strong,ketcheson2008highly}. 
For our numerical experiments we will use the explicit TVD/SSP-RK method of third order using three stages (SSPRK(3,3)), \cite{gottlieb1998total}. 
We note that energy stability for SSP-RK methods is guaranteed for all time if it holds for the standard first order explicit Euler method, \cite{gottlieb1998total}. 
In \cite{levy1998semidiscrete} it was shown in the case of 
linear CLs that the energy stability is preserved in time for some choices of SSP-RK methods, 
including SSPRK(3,3).\footnote{This is unfortunately generally not true in the nonlinear case, as the energy might increase 
after one iteration of the explicit Euler method if no dissipation is added to the numerical solution.} Thus we see that at least in the case of linear advection, both the weak RBF analytical method as well as the weak RBF collocation method can be used with SSPRK(3,3)  and have guaranteed energy stability. 
For the time step $\Delta t$ we use  
\begin{equation}  
  \Delta t = C \cdot \frac{|\Omega|}{N \max{|f'(u)|}}
\end{equation}
with $C = 0.1$ in the later numerical tests. 
Here, $\max{|f'(u)|}$ is calculated for all $u$ between $\min_{x \in \Omega} u_0(x)$ and $\max_{x \in \Omega} u_0(x)$. 
Note that for the linear advection equation we simply have $\max{|f'(u)|} = |\lambda|$.

\subsection{Implementation}
\label{sub:implementation} 

Since the implementation mainly consists of standard techniques, we omit any detailed discussion.  Additional information may be found in  \cite[Chapter 7.2.7]{glaubitz2020shock}. 
 
\section{Relationship to Other Methods}
\label{sec:connections}

For additional context, we now provide some comparisons to some techniques commonly used for solving hyperbolic conservation laws.

\subsection{DG Methods}
\label{sub:DG}
DG methods, see \cite{hesthaven2007nodal} and references therein, are perhaps the most obviously comparable. 
DG methods use a partition of the domain $\Omega$ into smaller 
elements $\Omega_i$ with $\bigcup_i \Omega_i = \Omega$. 
In each element the problem is discretized in a weak form similar to 
\eqref{eq:RBFspace}, where the numerical solution $u$ and the test 
functions $v$ are typically replaced by polynomials in every element 
$\Omega_i$. 
These polynomials are allowed to be discontinuous at the element 
interfaces and numerical fluxes are utilized to couple neighboring 
elements and to weakly enforce BCs.  
In this context, the proposed weak RBF method might be interpreted as 
a DG method in which a single big element 
${\Omega_i = \Omega}$ is used and the polynomial approximations are 
replaced with RBF interpolants. 
In a nodal approach this allows the use of more sophisticated sets of 
interpolation points, especially in higher dimensions (although these 
are not considered in this work). 
Note that by the Mairhuber--Curtis theorem 
\cite[Theorem 2]{iske2011scattered} polynomial interpolation in general is not 
well-defined in more than one dimension. 

\subsection{Spectral Galerkin Tau Methods}
\label{sub:galerkintau}

Spectral Galerkin methods solve the PDE in form of an integral equation as well, only {\em without} including the BCs in the integral equation.
The BCs can, for instance, be enforced directly by choosing suitable trial functions to span the approximation space, 
e.\,g.\ by choosing ${V_N = \mathrm{span}\{ \, \sin(\pi n x) \mid n=1,\dots,N \, \} }$ in case of homogeneous Dirichlet BCs on ${\Omega = [0,1]}$. 
The so-called \emph{spectral Galerkin tau methods}, see \cite{canuto2006spectral} and 
references therein,
use trial functions that do not have to individually satisfy the BCs, but rather some additional equations are 
imposed to ensure the numerical solution satisfies BCs.
To maintain a well-posed discretization, i.\,e., the number of equations being equal to the number of degrees of 
freedom, some of the integral integrations corresponding to the highest order test functions are dropped in favor of 
the BC equations. 
In the weak RBF method, on the other hand, these BC equations include numerical flux functions and are incorporated 
into the integral equations corresponding to the test functions. 
As a consequence, we do not need to remove any test functions from the integral equations, yielding higher order of 
accuracy. 

\subsection{Penalty-Type Boundary Treatment in Pseudospectral Methods}
\label{sub:penalty}

As with strong RBF methods, classical pseudospectral methods typically are built from bases of Fourier, Chebyshev or 
Legendre polynomials, and require that the BCs are strongly (exactly) imposed, see \cite{gottlieb2001spectral} and 
references therein. 
Penalty methods, i.\,e.\ using a penalty term for treating BCs, have been used both for spectral 
methods in the weak \cite{canuto1982error} and strong \cite{funaro1988new,funaro1991convergence} forms.
The basic idea behind penalty methods is that it suffices to impose the BCs to the order of the given scheme, 
which can be done by introducing a penalty term into the discretized equation. 
In particular, the BCs have to be satisfied exactly by the numerical solution only in the limit of infinite order. 
Depending on the method and problem under consideration it may be challenging to construct suitable penalty terms. 

In the weak RBF method, such penalty terms are derived somewhat naturally by utilizing numerical flux functions. 
As a consequence, a large class of penalty terms may be available for practical use.  
Future work will address the development of stable RBF methods in strong form. 
As discussed above, a bottleneck for such an investigation will be the development of suitable penalty 
terms for the boundary treatment in a strong RBF method. 
This is consistent with the observation that classic strong RBF methods (in which BCs are imposed strongly), so far, 
have only been observed to be stable if no BCs were present \cite{platte2006eigenvalue}. 

\section{Possible Extensions for the Proposed Boundary Treatment}
\label{sec:extension} 

We now address some possible extension of the proposed boundary treatment in global RBF methods.

\subsection{Formulation in Multiple Dimensions} 
\label{sub:2d}

Let $\Omega \subset \R^m$ be a bounded region with piecewise smooth boundary $\partial \Omega$. 
The $m$ dimensional equivalent of the one dimensional CL \eqref{eq:cl} is given by 
\begin{equation}\label{eq:cl-2d}
	u_t + \nabla \cdot \boldsymbol{F}(u) = 0, 
	\quad \boldsymbol{x} \in \Omega, \ t > 0,
\end{equation}
where $F: \R \to \R^m$, $\nabla = (\partial_{x_1}, \dots, \partial_{x_m})$ is the formal \emph{nabla operator},  and $\cdot$ denotes their inner product. 
We also assume we are given suitable IC and BCs.  After applying the divergence theorem,
the weak form of \eqref{eq:cl-2d} reads 
\begin{equation*}
	\int_{\Omega} u_t v \intd V - \int_{\Omega} \boldsymbol{F}(u) \cdot \nabla v \intd V + \oint_{\partial \Omega} v \boldsymbol{F}(u) \cdot \mathbf{n} \intd S = 0
\end{equation*}
with test function $v \in C^1(\Omega)$. 
It should be stressed that the closed manifold $\partial \Omega$ is assumed to be oriented by outward pointing normals, and $\mathbf{n}$ denotes the outward pointing unit normal at each point on the boundary $\partial \Omega$. 

Following the ideas discussed in \S \ref{sec:stable_RBFs}, the corresponding ($m$ dimensional) weak RBF collocation method is defined as follows: 
Find $u_N \in V_{N,P}$ such that 
\begin{equation}\label{eq:weak-RBF-mD}
	\int_{\Omega} (u_N)_t v_N \intd V - \int_{\Omega} \boldsymbol{F}_N \cdot \nabla v_N \intd V + \oint_{\partial \Omega} v_N \mathbf{F}^{\text{num}} \cdot \mathbf{n} \intd S = 0
\end{equation}
for all $v_n \in V_{N,P}$. 
Note that in this case $u_N$ and $v_N$ still denote scalar-valued RBF approximations.  At the same time $\boldsymbol{F}_N$ denotes a vector-valued function for which every component has been replaced by an RBF approximation. 
Consequently, $\mathbf{F}^{\text{num}}$ also denotes an $m$-dimensional numerical flux function.

\subsection{Stability in Multiple Dimensions} 
\label{sub:2d-stability}

A similar analysis to the one in \S \ref{sub:coll-RBF-method} can be used in the linear case, that is for $\boldsymbol{F}(u) = \boldsymbol{\lambda} u$ with $\boldsymbol{\lambda} \in \R^m$. 
In particular, by choosing $v_N = u_N$ in \eqref{eq:weak-RBF-mD} and applying Gauss's divergence theorem, we obtain
\begin{equation}\label{eq:weak-RBF-mD-L2}
	\frac{\d}{\d t} \norm{u_N}_{L^2}^2 
		= \oint_{\partial \Omega} u_N \left[ u_N \boldsymbol{\lambda} - 2 \mathbf{F}^{\text{num}} \right] \cdot \mathbf{n} \intd S.
\end{equation} 
This equation can be considered as the $m$-dimensional analogue of \eqref{eq:weak-RBF-coll-L2}. 
It is unfortunately less clear in general how the boundary contributions sum up in the higher-dimensional setting. 
Indeed, the boundary integral in \eqref{eq:weak-RBF-mD-L2} strongly depends on the bounded region $\Omega$ as well as the sign of the different components of the constant velocity vector $\boldsymbol{\lambda} \in \R^m$. 
That said, Example \ref{ex:2dcube} suggests that in theory similar stability results as in \S \ref{sub:coll-RBF-method} are also obtainable in multiple dimensions.  They might be more cumbersome to formulate, however. 

\begin{example}
\label{ex:2dcube}
	Suppose we are given the two-dimensional cube $\Omega = [a,b]^2$ and a nonnegative velocity vector $\boldsymbol{\lambda} = (\lambda_1,\lambda_2)^T$ with $\lambda_1,\lambda_2 \geq 0$. 
	In this case the boundary $\partial \Omega$ can be partitioned into the four following lines: 
	\begin{equation*}
	\begin{aligned} 
		& \partial \Omega_W = \{ \, (a,y)^T \in \R^2 \mid a \leq y \leq b \, \},  
		&& \partial \Omega_E = \{ \, (b,y)^T \in \R^2 \mid a \leq y \leq b \, \}, \\ 
		& \partial \Omega_S = \{ \, (x,a)^T \in \R^2 \mid a \leq x \leq b \, \},  
		&& \partial \Omega_N = \{ \, (x,b)^T \in \R^2 \mid a \leq x \leq b \, \}.
	\end{aligned}
	\end{equation*} 
	Observe that  $\partial \Omega_E$ and $\partial \Omega_N$ are the outflow part of the boundary (no BC is given there), while $\partial \Omega_W$ and $\partial \Omega_S$ are the inflow part (BCs are given there). 
	Focusing on periodic BCs, for which know that the energy should not increase over time, we have 
	\begin{equation}
\label{eq:periodbc}
	\begin{aligned}
		& u(t,x,y) = u(t,x+b-a,y) && \text{for } (x,y) = (a,y) \in \partial \Omega_W, \\ 
		& u(t,x,y) = u(t,x,y+b-a) && \text{for } (x,y) = (x,a) \in \partial \Omega_S.
	\end{aligned} 
	\end{equation}
	For simplicity we choose the upwind flux $\mathbf{F}^{\text{num}} = \mathbf{F}^{\text{num}}(a,b)$, satisfying 
	\begin{equation}
\label{eq:upwind}
		\mathbf{F}^{\text{num}}(a,b) \cdot \mathbf{n} = 
		\begin{cases}
		(\boldsymbol{\lambda} \cdot \mathbf{n}) a &; \ \boldsymbol{\lambda} \cdot \mathbf{n} 	\geq 0, \\ 
		(\boldsymbol{\lambda} \cdot \mathbf{n}) b &; \ \boldsymbol{\lambda} \cdot \mathbf{n} < 0.
		\end{cases} 
	\end{equation} 
	Substituting  \eqref{eq:periodbc} and \eqref{eq:upwind}  into \eqref{eq:weak-RBF-mD-L2} we obtain
	\begin{equation*}
	\begin{aligned} 
		\frac{\d}{\d t} \norm{u_N}_{L^2}^2 
		= & \int_{\partial \Omega_W} u_N(t,a,y) \left[ u_N(t,a,y) - 2u_N(t,b,y) \right] \boldsymbol{\lambda} \cdot \mathbf{n} \intd S \\ 
		& + \int_{\partial \Omega_S} u_N(t,x,a) \left[ u_N(t,x,a) - 2u_N(t,x,b) \right] \boldsymbol{\lambda} \cdot \mathbf{n} \intd S \\ 
		& - \int_{\partial \Omega_E} u_N^2(t,b,y) \boldsymbol{\lambda} \cdot \mathbf{n} \intd S  
			- \int_{\partial \Omega_E} u_N^2(t,x,b) \boldsymbol{\lambda} \cdot \mathbf{n} \intd S \\ 
		= & - \lambda_1 \int_a^b u_N^2(t,a,y) - 2 u_N(t,a,y) u_N(t,b,y) \intd y \\  
			& - \lambda_2 \int_a^b u_N^2(t,x,a) - 2 u_N(t,x,a) u_N(t,x,b) \intd x \\
			& - \lambda_1 \int_a^b u_N^2(t,b,y) \intd y - \lambda_2 \int_a^b u_N^2(t,x,b) \intd x \\ 
		= & - \lambda_1 \int_a^b \left[ u_N(t,a,y) - u_N(t,b,y) \right]^2 \intd y \\
		 	& - \lambda_2 \int_a^b \left[ u_N(t,x,a) - u_N(t,x,b) \right]^2 \intd x \\ 
		\leq & 0.
	\end{aligned}
	\end{equation*}
\end{example}
Hence we observe from Example \ref{ex:2dcube} that linear stability for the weak RBF method might also hold in higher dimensions as well as  more general domains.\footnote{A more rigorous study is clearly needed and will be included in  future investigations.}

\subsection{Numerical Integration}
\label{sub:num-int}

Constructing the mass and stiffness matrices requires  computing integrals which may be costly depending on the number of degrees of freedom and the dimension.
Preliminary tests presented in \S \ref{sec:numerical} indicate that it is possible to increase efficiency without  reducing accuracy, either by using trapezoidal,  Gauss-Legendre, or Gauss-Lobatto rules (in one dimension), and their tensor products in higher domains when a rectangular domain is assumed,
see for example \cite{haber1970numerical,stroud1971approximate,engels1980numerical,davis2007methods,trefethen2017cubature} for general discussions on numerical quadrature.
Such techniques are not readily available for non-standard (non-rectangular) domains.  In this case an alternative
 might be to use  classical (quasi-)Monte Carlo methods, \cite{metropolis1949monte,niederreiter1992random,caflisch1998monte,dick2013high}, or more recently developed high-order least squares cubature rules, \cite{glaubitz2020stableCF,glaubitz2020constructing}, which are based on one-dimensional  approaches developed in \cite{wilson1970necessary,wilson1970discrete,huybrechs2009stable,glaubitz2020stable}.
Future work will address the advantages and potential difficulties in replacing these integrals by various numerical formulas.

\subsection{Local Radial Basis Function Methods} 

We have thus far only considered global RBF methods.  One obvious concern in using global RBFs is the associated computational cost. 
Specifically, determining a global RBF interpolant as well as  calculating the corresponding  differentiation matrix each cost $\mathcal{O}(N^3)$ operations for $N$ nodes. 
While for the discussed methods this can be done offline, that is once before time stepping commences  (assuming the nodes do not change over time),  there are additional $\mathcal{O}(N^2)$ operations to be performed each time a differentiation matrix is applied  during time stepping. 
Local RBF-FD are designed to remedy this problem.\footnote{The conference presentation \cite{tolstykh2000using} by Tolstykh in 2000 seems to be the earliest reference to RBF-FD methods.} 
Conceptually, these methods can be interpreted as an extreme case of overlapping domain decomposition, with a separate domain surrounding each node. 
The basic idea is to center a local RBF-FD stencil at each of the $N$ global nodes, and let it include the $n-1$ nearest neighbors, where $n \ll N$. 
For every node, and based on its surrounding stencil, a local FD formula that is exact for all RBF interpolants on that stencil---potentially including polynomials---is then derived from a system of linear equations similar to \eqref{eq:LS-RBF-interpol}. 
The main difference is that the right hand side of the linear system is replaced by the nodal values of a linear differentiation operator.
For more details, see \cite[Chapter 5]{fornberg2015primer} and references therein. 

We note that going from the strong to weak formulation of the underlying conservation law is also possible for local RBF-FD methods.
Although the conservation and energy stability proofs do not immediately follow, such results may be possible at least in the linear case when 
replacing exact integrals and differentiation operators by their discrete counterparts,  as long as certain summation-by-parts (SBP) properties are satisfied, \cite{svard2014review,fernandez2014review}. 
In this case, many stability properties which are based on integration by parts, i.e.~the continuous analogue of SBP, would still be satisfied in a discrete norm. 
This idea is also left for future investigations. 

\section{Numerical Results}
\label{sec:numerical}

We now demonstrate our theoretical findings for the weak RBF analytical and collocation methods. 
In most tests we focus on the cubic and quintic kernel, $\varphi(r) = r^3$ and $\varphi(r) = r^5$, which belong to the class of polyharmonic splines (PHSs). 
Although they  yield algebraic rather than spectral accuracy\footnote{For a discussion on the accuracy of infinitely smooth kernels as well as the role of the shape parameter $\varepsilon$ see for instance \cite{madych1992miscellaneous,schaback1995error,buhmann2003radial,fornberg2015primer} and references therein.}, there are several advantages associated with PHSs, see \cite{iske2003radial,fornberg2015primer,iske2020ten}. 
In particular, PHSs satisfy certain optimality results \cite{duchon1977splines,powell1992theory} that can be interpreted as multidimensional scattered node analogues of the one-dimensional result that the natural cubic spline, among all possible interpolants $s$, minimizes $\int [s''(x)]^2 \intd x$ over the interval spanned by the nodes. 
Essentially this means that PHSs interpolate scattered data with the fewest spurious oscillations. 
Finally, PHSs do not require a (sometimes cumbersome) selection of the shape parameter $\varepsilon$. 
The MATLAB code used to generate the subsequent numerical tests can be found at \cite{glaubitz2020weakRBFcode}.

\subsection{Linear Advection Equation} 
\label{sub:linear_adv}

Let us consider the linear initial value problem (IVP) 
\begin{equation}\label{eq:linear_1d}
	u_t + u_x = 0, \quad 
	u(0,x) = \exp \left( -20 x^2 \right)
\end{equation} 
with $x \in \Omega = [-1,1]$ and $t > 0$.
We will also consider periodic and inflow BCs given respectively by 
\begin{align}
	u(t,-1) & = u(t,1), && \text{(periodic BC)} \label{eq:linear_BC_periodic} \\
	u(t,-1) & = u(0,1-\text{mod}(t,2)) && \text{(inflow BC)} \label{eq:linear_BC_inflow}
\end{align}
Note that both BCs yield the same exact solution.

\subsubsection{Solution, Momentum and Energy Profiles}
\label{sub:profiles}
We start by comparing numerical solutions given by the weak RBF methods for $P=0$ (no polynomials included) and $P =1$ (constants included) to the standard RBF collocation method. 
The latter will be subsequently referred to as the \em{strong} or \em{standard} RBF (collocation) method. 
Note that in the linear advection case the weak RBF analytical and collocation methods are the same.

\begin{figure}[!htb]
  	\centering
 	\begin{subfigure}[b]{0.33\textwidth}
    		\includegraphics[width=\textwidth]{%
      		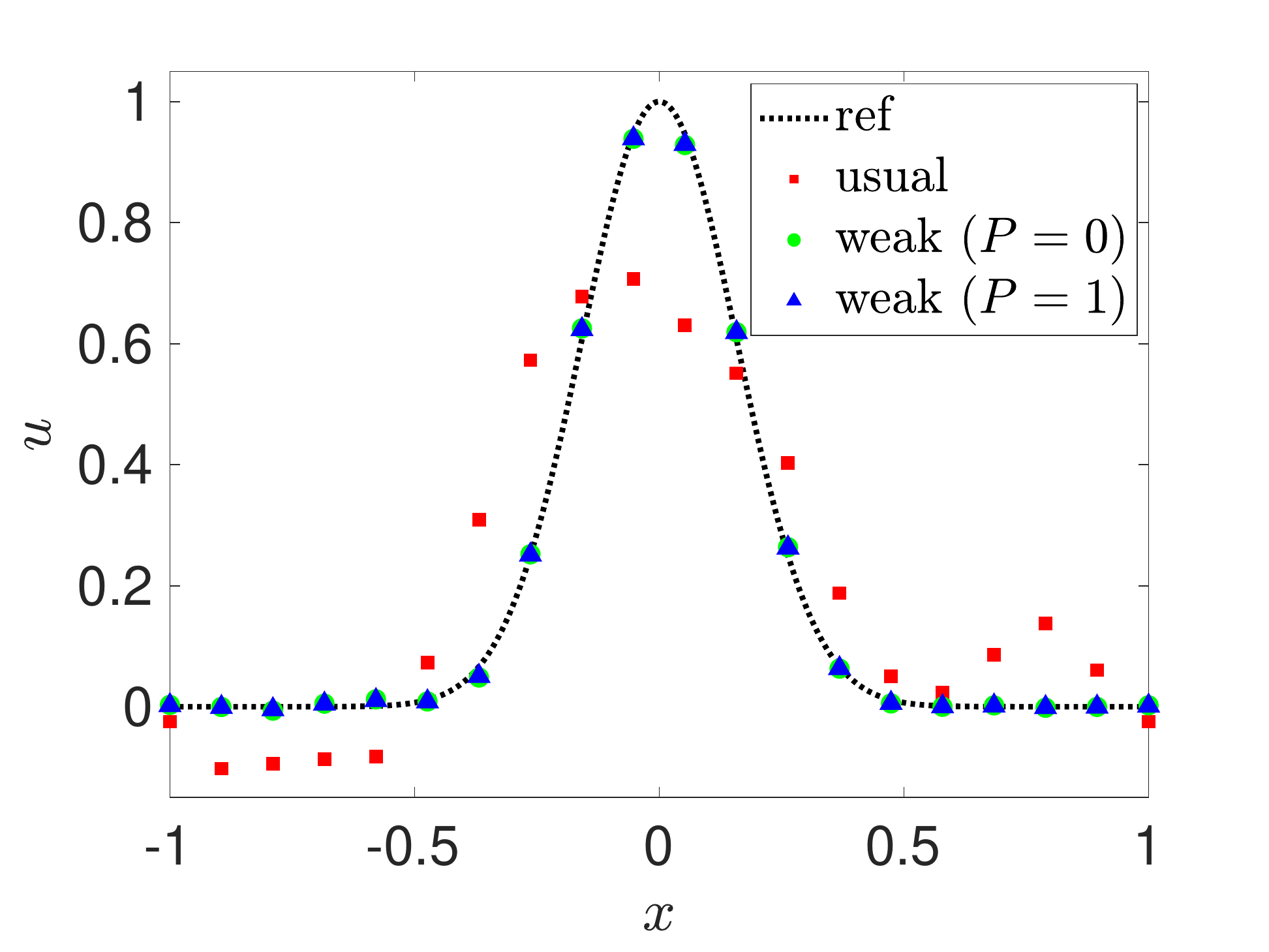}
    		\caption{Solution, cubic}
    		\label{fig:linear_1d_A_cubic_sol} 
  	\end{subfigure}%
  	~ 
  	\begin{subfigure}[b]{0.33\textwidth}
    		\includegraphics[width=\textwidth]{%
      		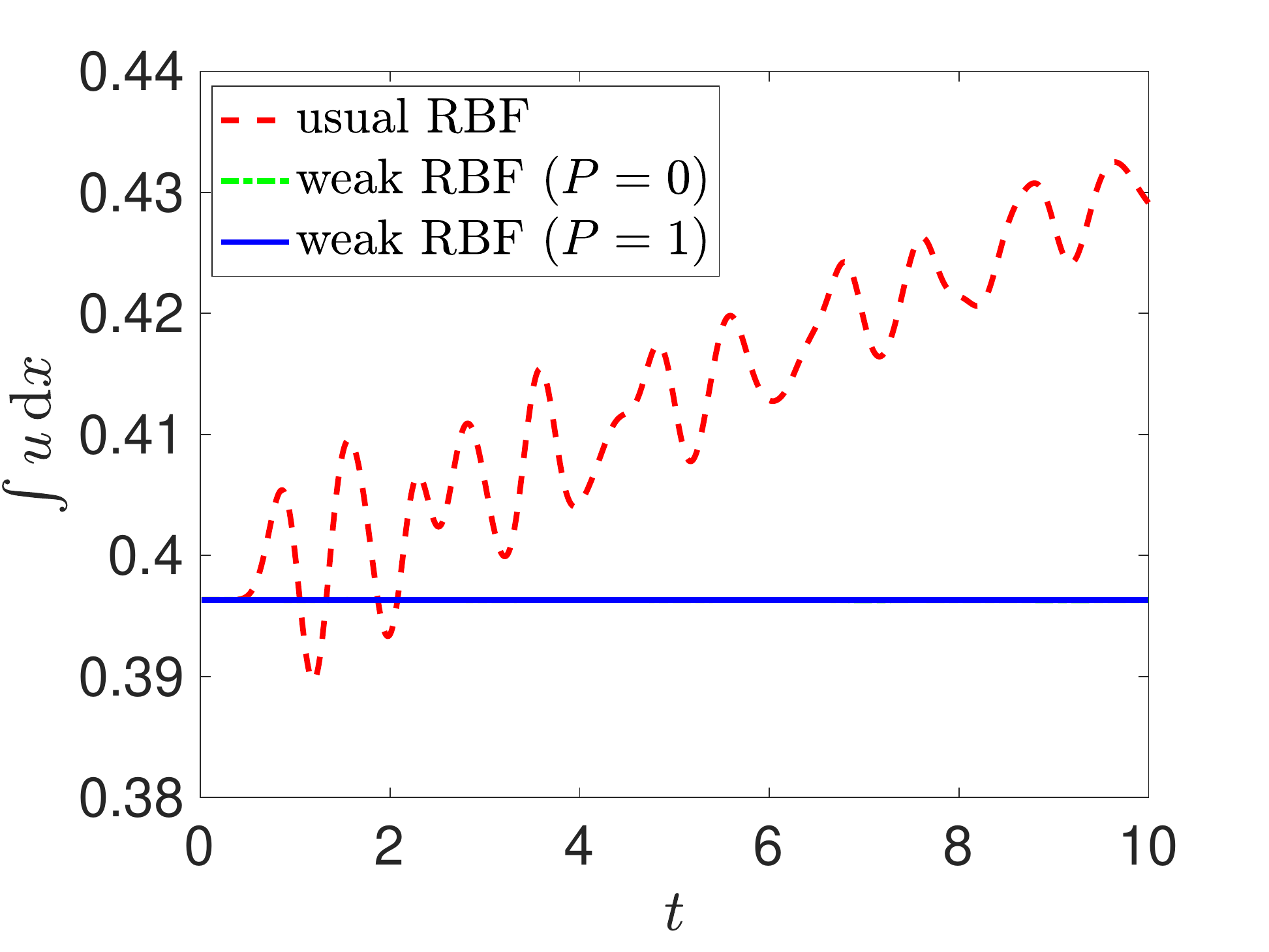}
    		\caption{Momentum, cubic}
    		\label{fig:linear_1d_A_cubic_mom} 
  	\end{subfigure}%
	~ 
  	\begin{subfigure}[b]{0.33\textwidth}
    		\includegraphics[width=\textwidth]{%
      		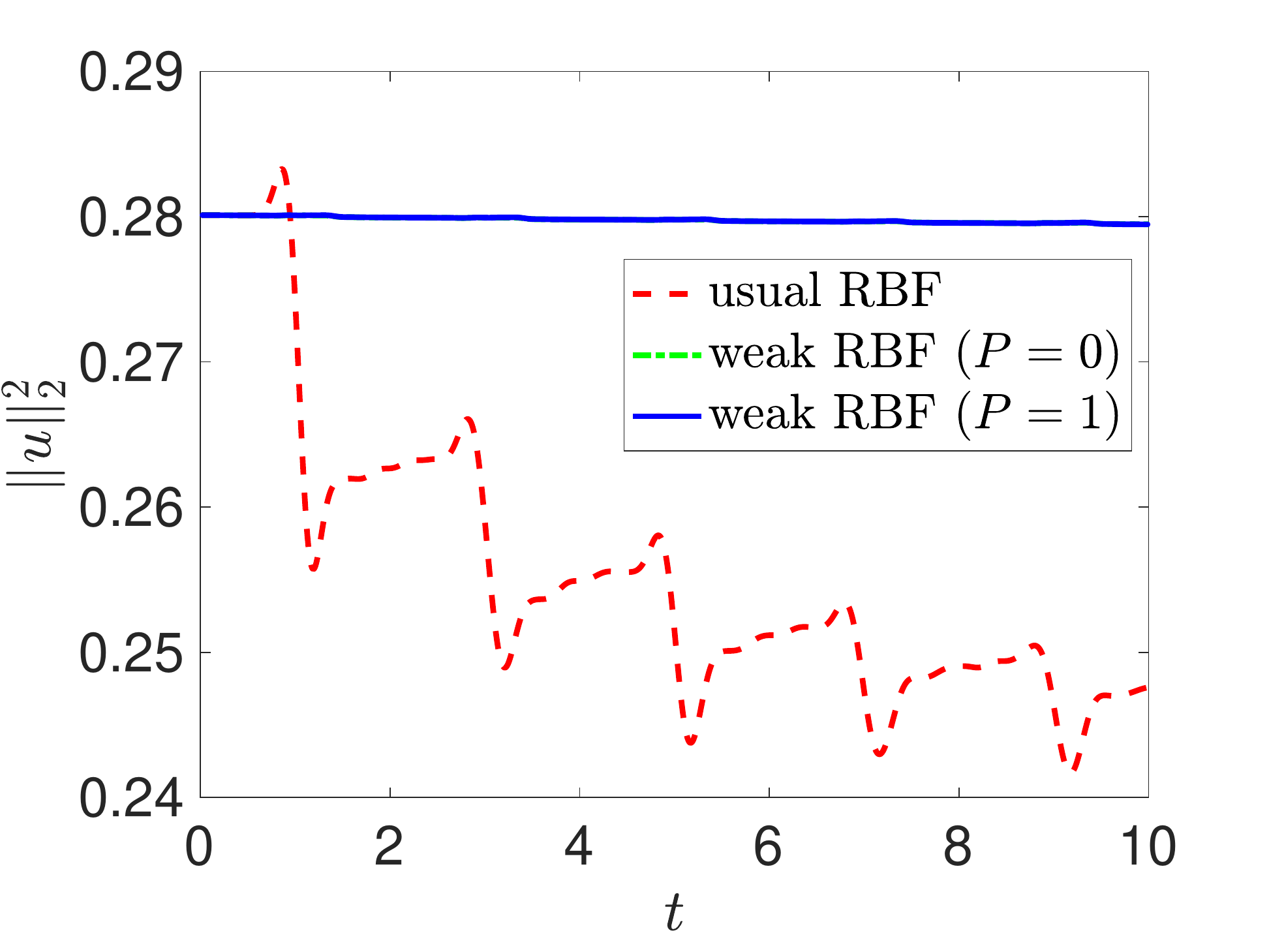}
    		\caption{Energy, cubic}
    		\label{fig:linear_1d_A_cubic_en} 
  	\end{subfigure}%
	\\
	\begin{subfigure}[b]{0.33\textwidth}
    		\includegraphics[width=\textwidth]{%
      		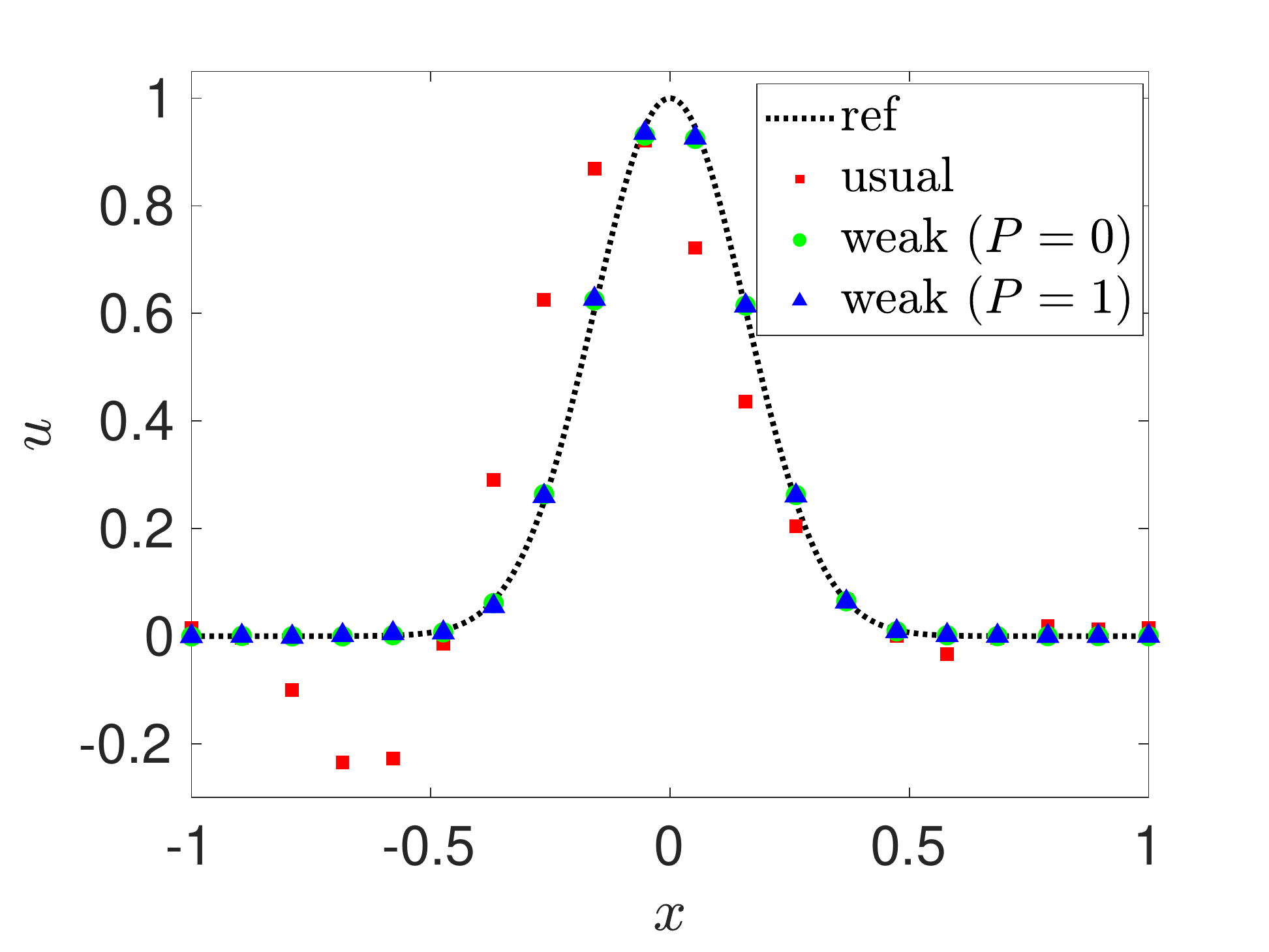}
    		\caption{Solution, quintic}
    		\label{fig:linear_1d_A_quintic_sol} 
  	\end{subfigure}%
  	~ 
  	\begin{subfigure}[b]{0.33\textwidth}
    		\includegraphics[width=\textwidth]{%
      		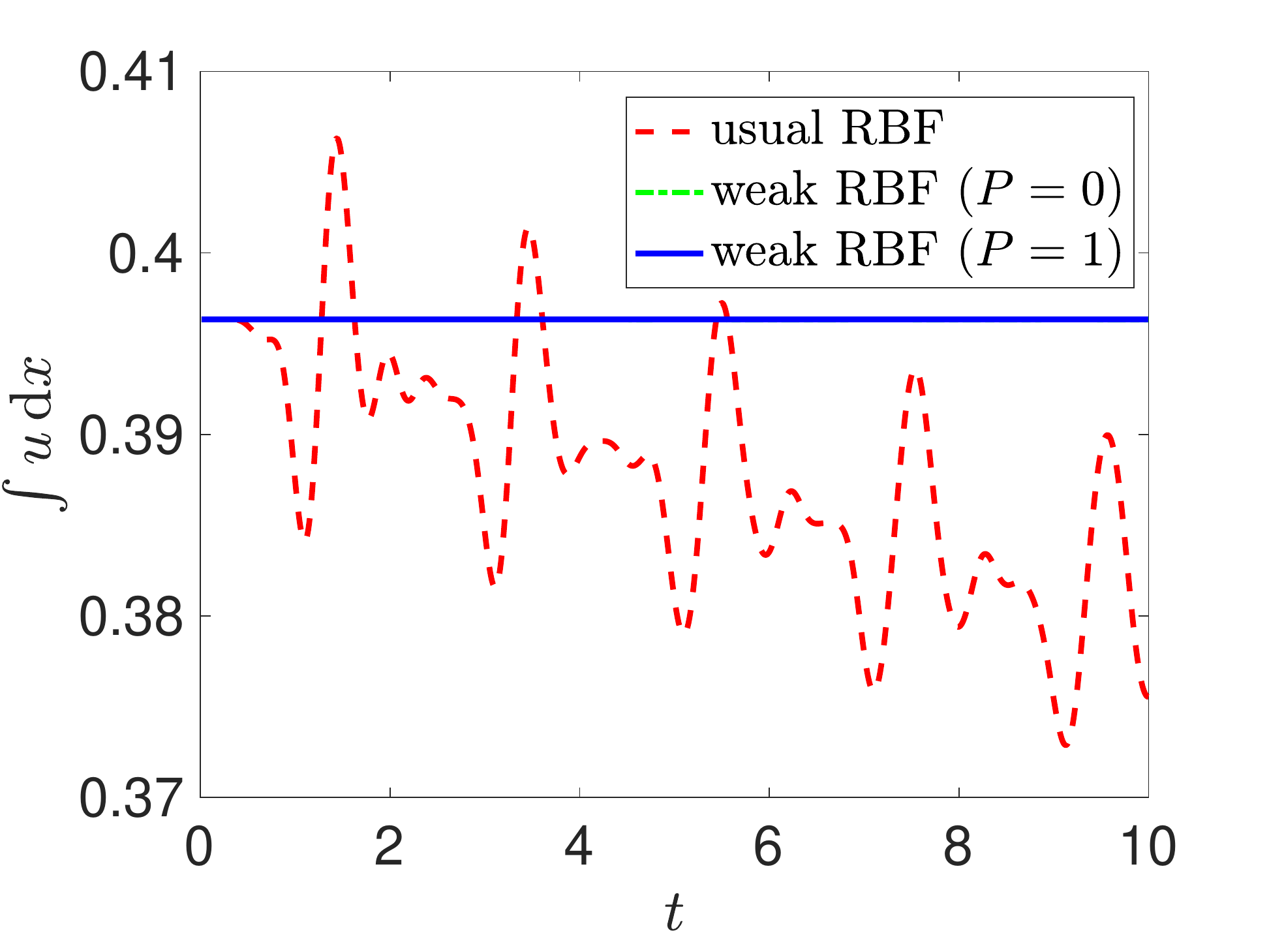}
    		\caption{Momentum, quintic}
    		\label{fig:linear_1d_A_quintic_mom} 
  	\end{subfigure}%
	~ 
  	\begin{subfigure}[b]{0.33\textwidth}
    		\includegraphics[width=\textwidth]{%
      		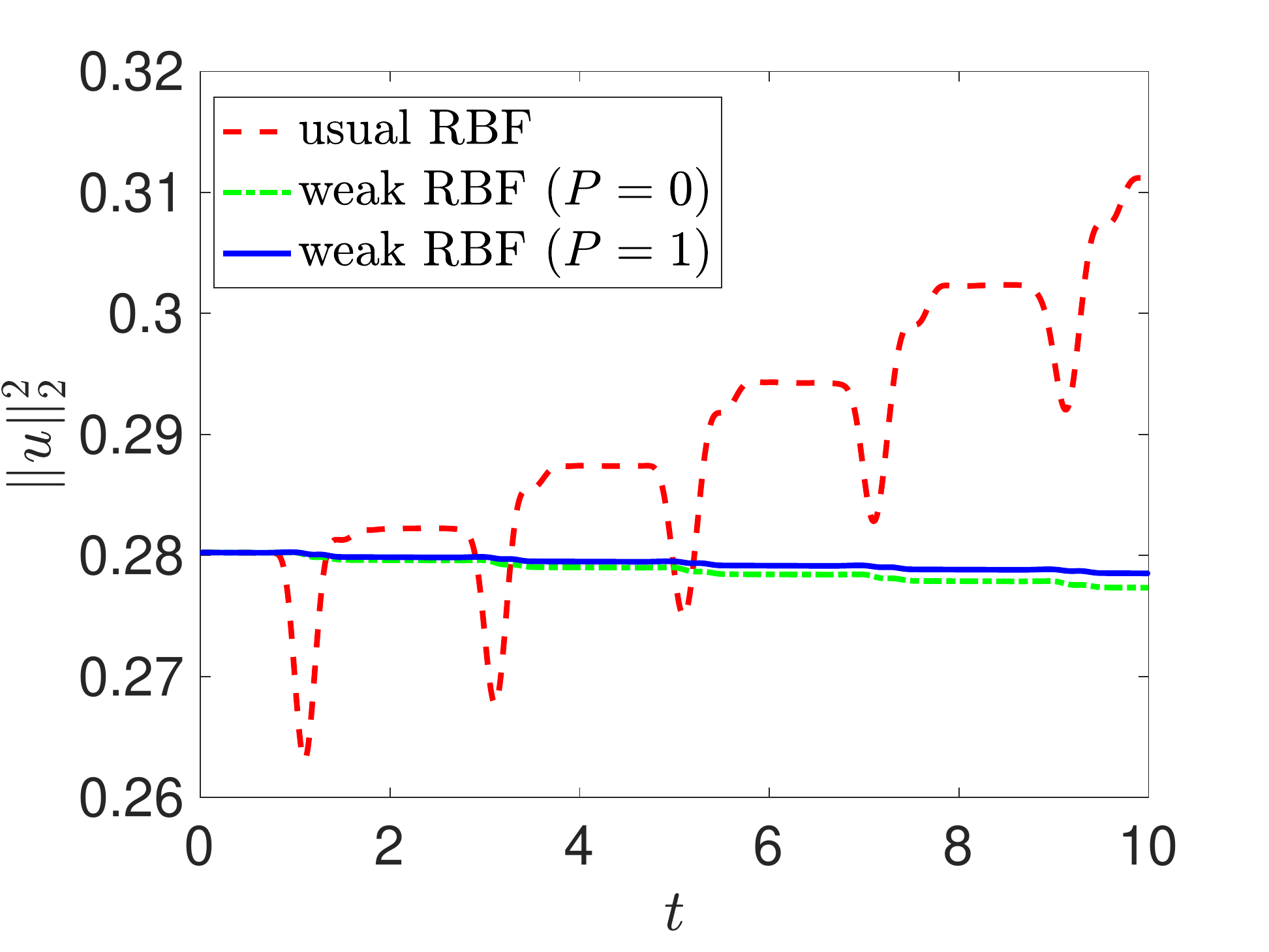}
    		\caption{Energy, quintic}
    		\label{fig:linear_1d_A_quintic_en} 
  	\end{subfigure}%
  	\caption{
  		Numerical solutions at $t=10$ (left); their momentum (middle); and energy (right) over time for $u_t + u_x = 0$ with periodic BC \eqref{eq:linear_BC_periodic}. 
  		In all cases, $N = 20$ equidistant nodes and shape parameter $\varepsilon=5$ were used. 
  	}
  	\label{fig:linear_1d_A}
\end{figure} 

\begin{figure}[!htb]
  	\centering
	\begin{subfigure}[b]{0.33\textwidth}
    		\includegraphics[width=\textwidth]{%
      		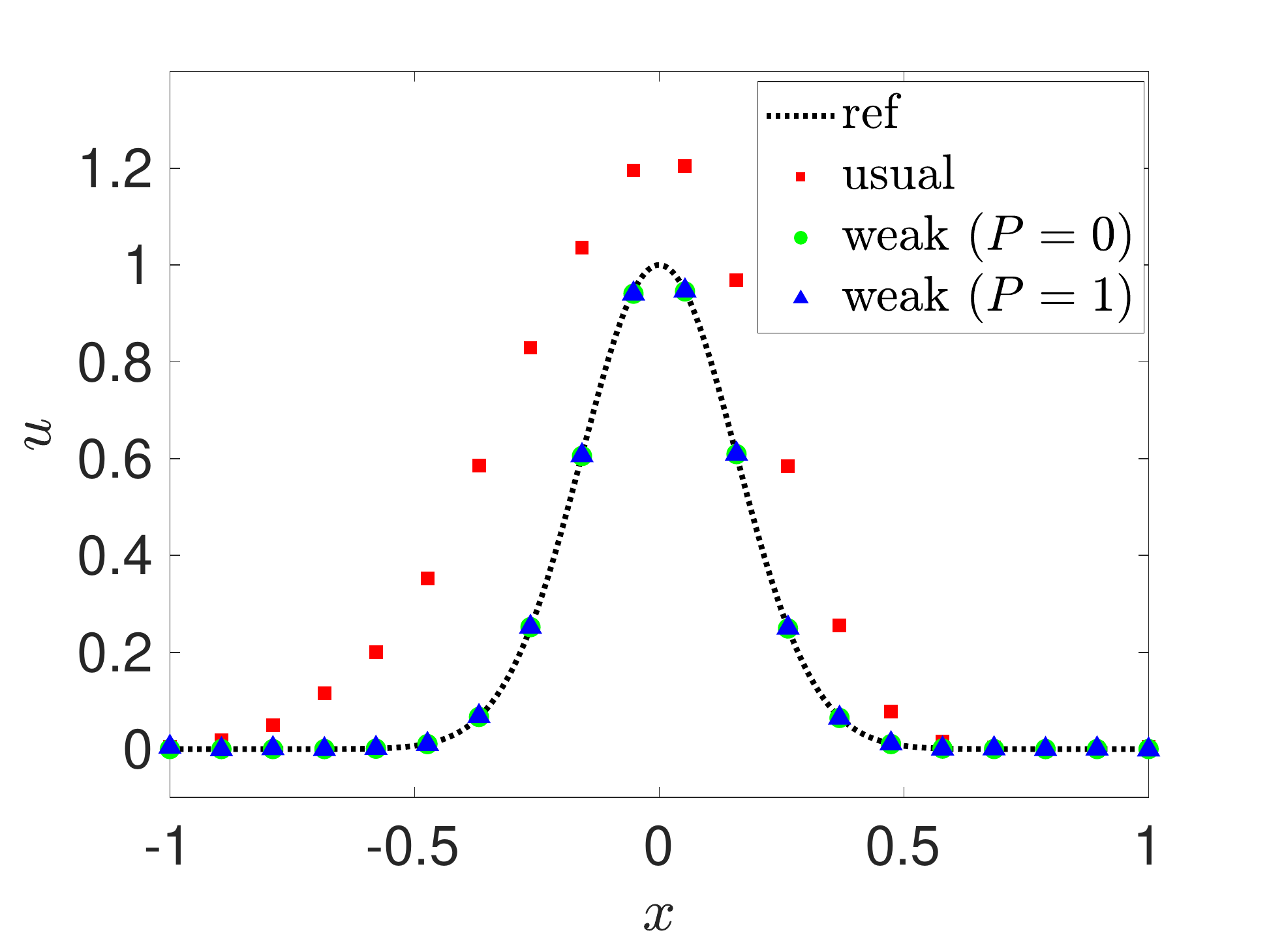}
    		\caption{Solution, Gauss}
    		\label{fig:linear_1d_A_G_sol} 
  	\end{subfigure}%
  	~ 
  	\begin{subfigure}[b]{0.33\textwidth}
    		\includegraphics[width=\textwidth]{%
      		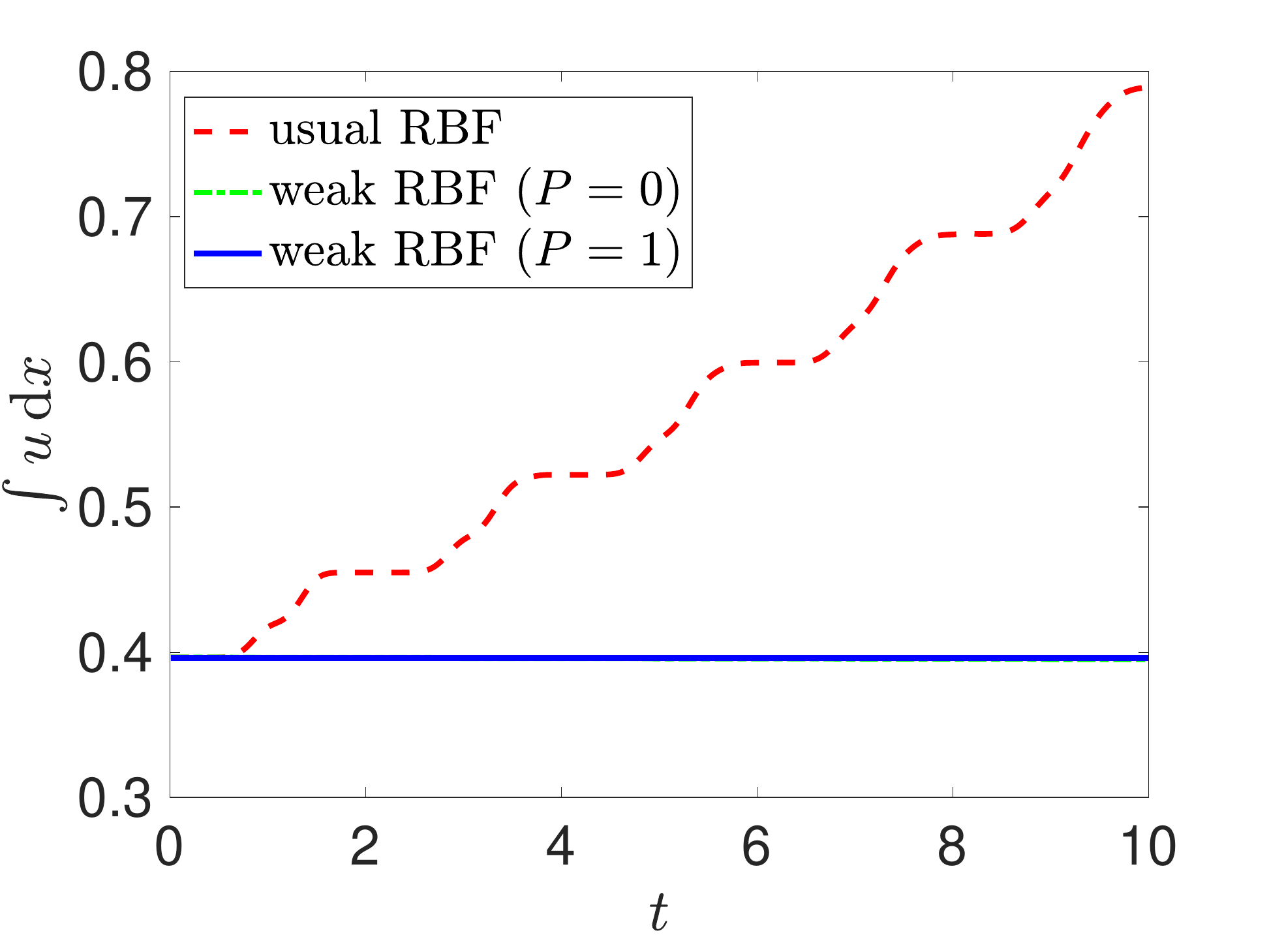}
    		\caption{Momentum, Gauss}
    		\label{fig:linear_1d_A_G_mom} 
  	\end{subfigure}%
	~ 
  	\begin{subfigure}[b]{0.33\textwidth}
    		\includegraphics[width=\textwidth]{%
      		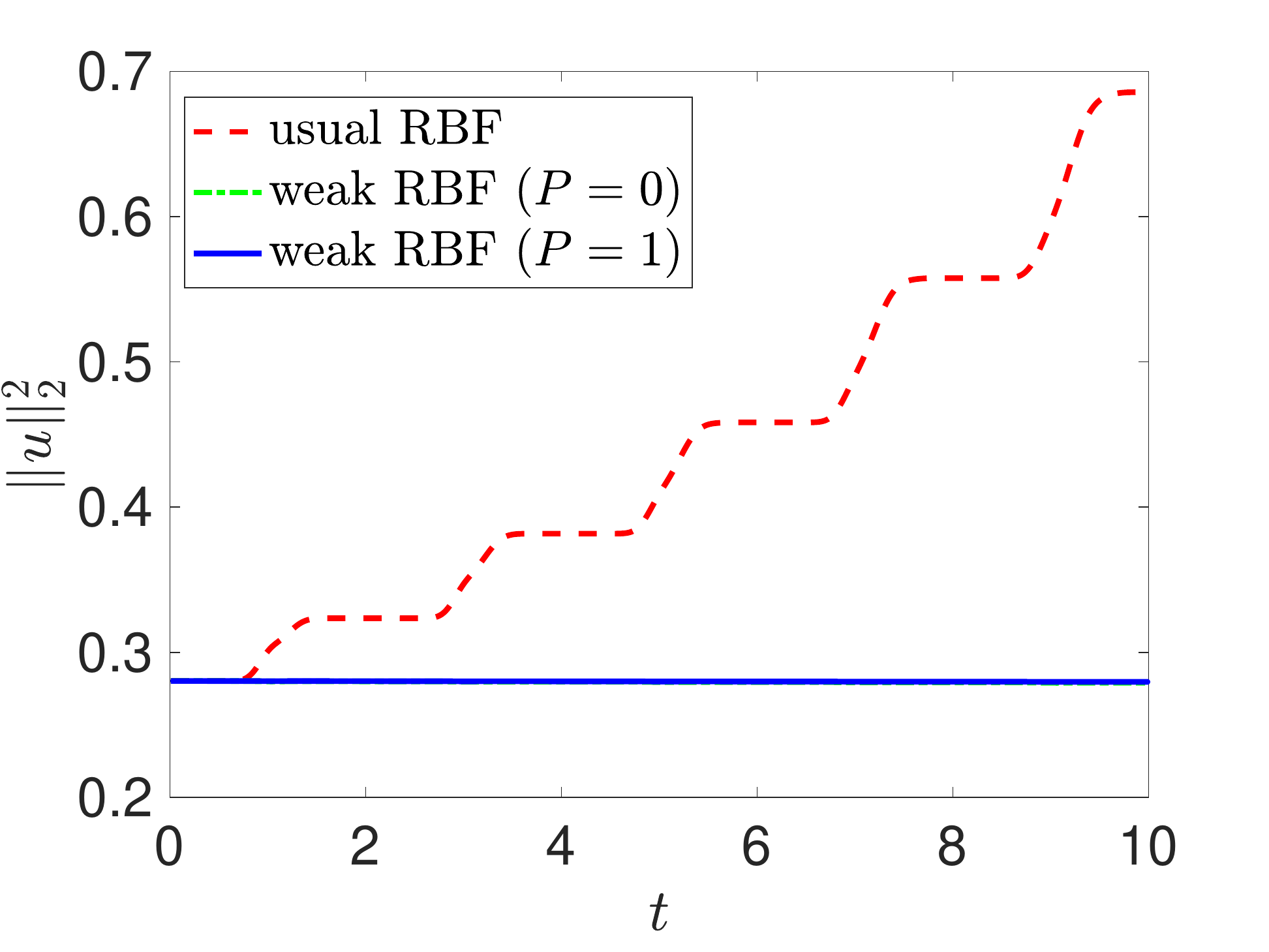}
    		\caption{Energy, Gauss}
    		\label{fig:linear_1d_A_G_en} 
  	\end{subfigure}%
	\\
	\begin{subfigure}[b]{0.33\textwidth}
    		\includegraphics[width=\textwidth]{%
      		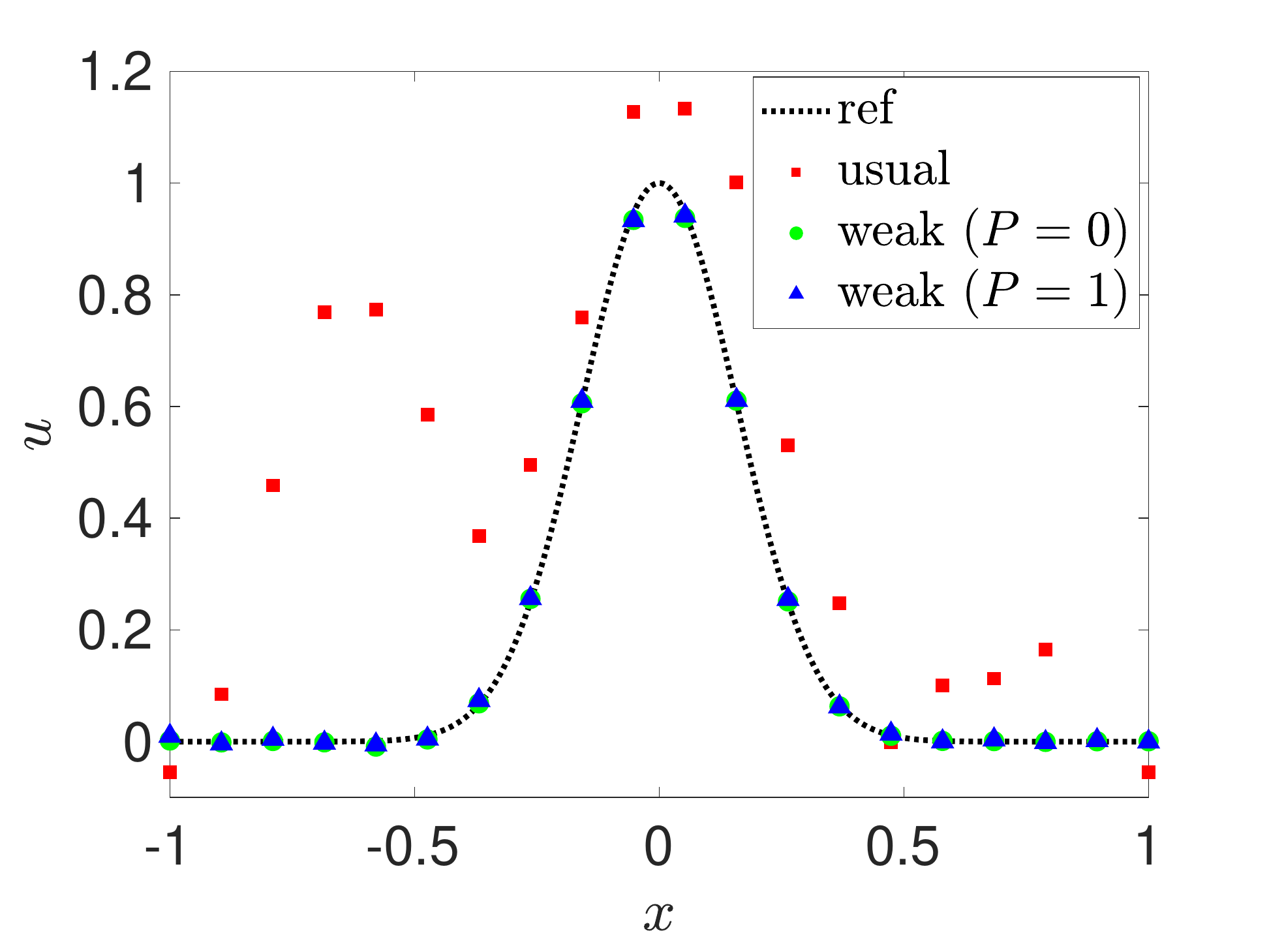}
    		\caption{Solution, IQ}
    		\label{fig:linear_1d_A_IQ_sol} 
  	\end{subfigure}%
  	~ 
  	\begin{subfigure}[b]{0.33\textwidth}
    		\includegraphics[width=\textwidth]{%
      		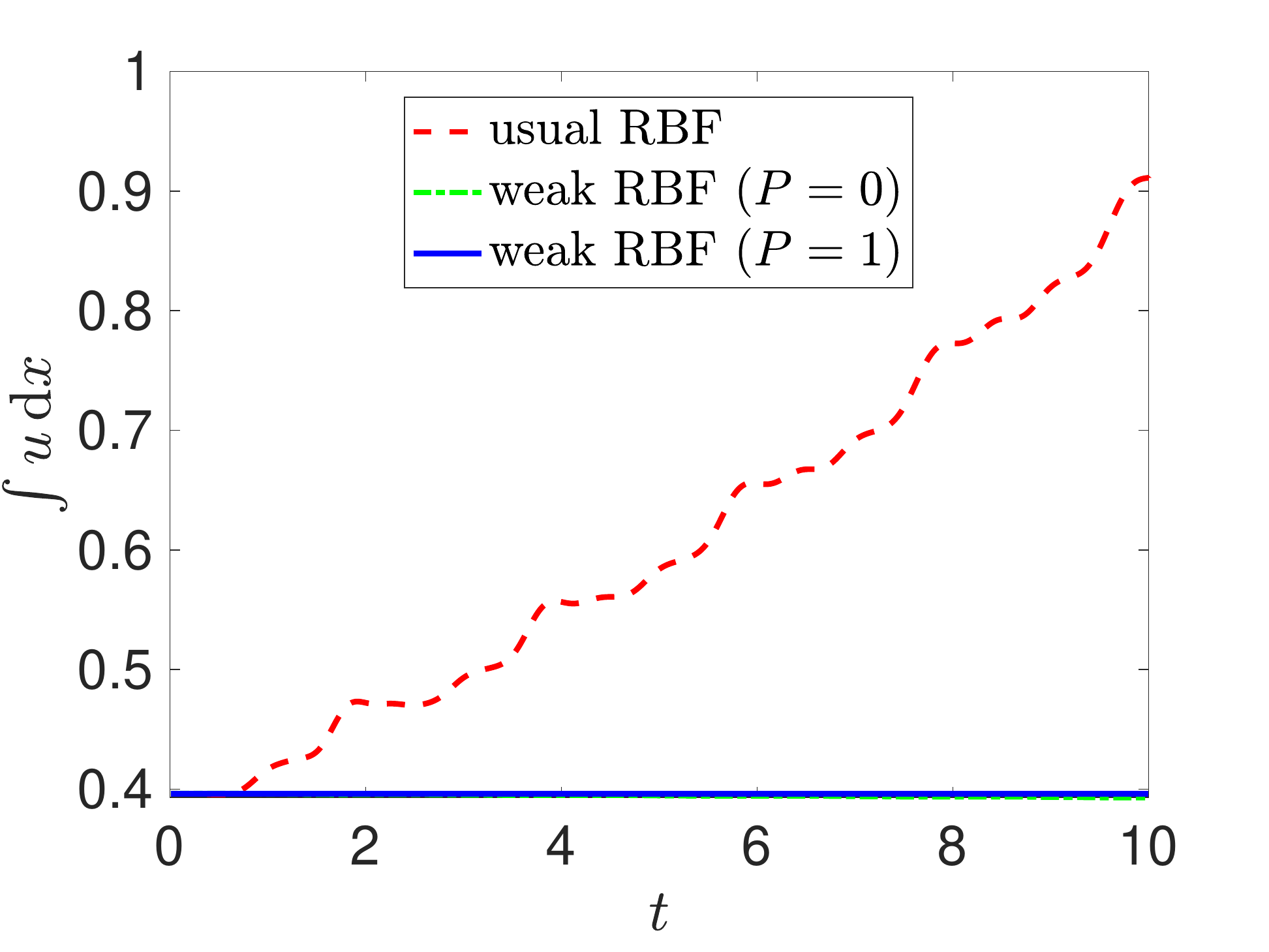}
    		\caption{Momentum, IQ}
    		\label{fig:linear_1d_A_IQ_mom} 
  	\end{subfigure}%
	~ 
  	\begin{subfigure}[b]{0.33\textwidth}
    		\includegraphics[width=\textwidth]{%
      		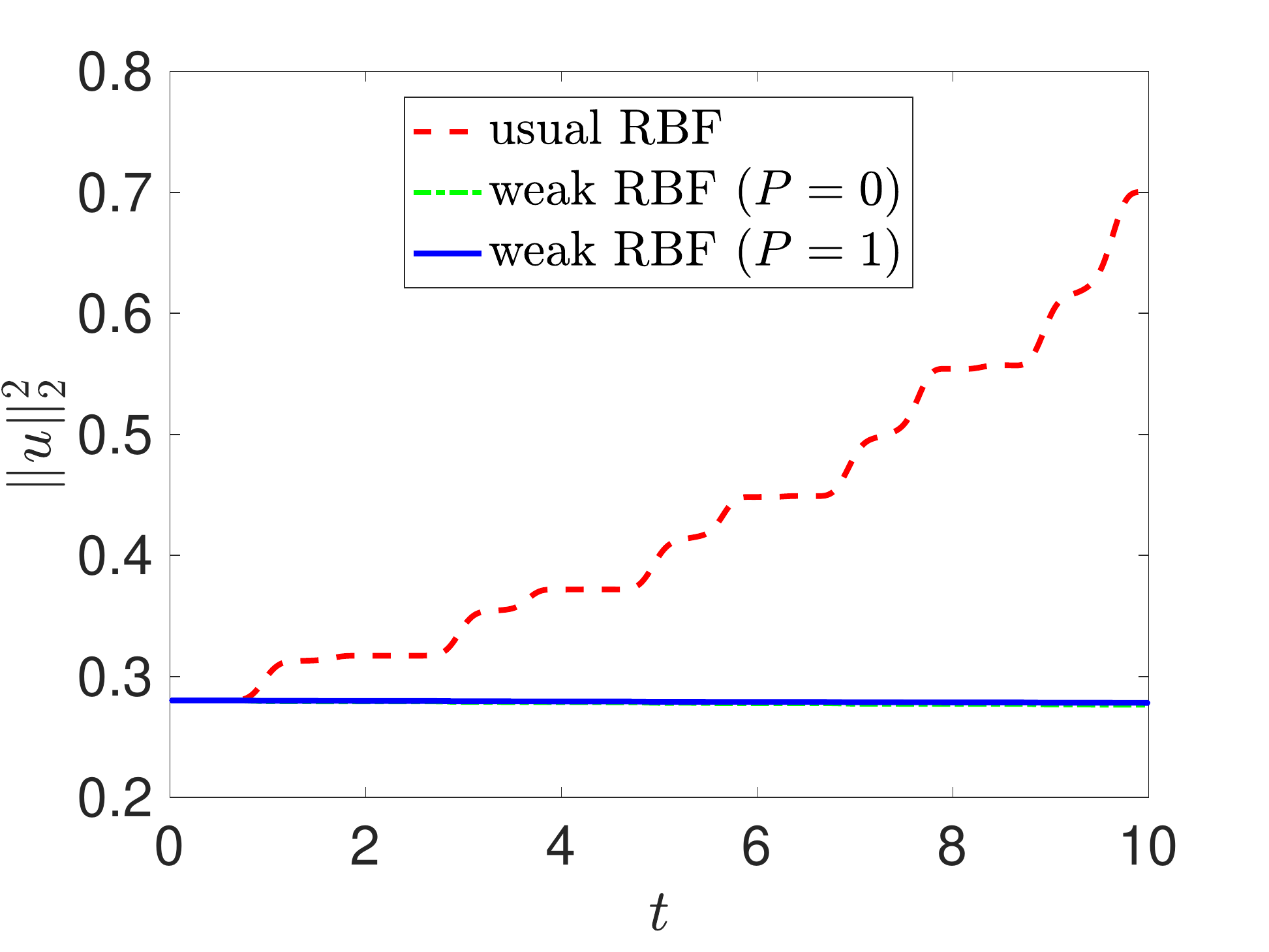}
    		\caption{Energy, IQ}
    		\label{fig:linear_1d_A_IQ_en} 
  	\end{subfigure}%
	\\
	\begin{subfigure}[b]{0.33\textwidth}
    		\includegraphics[width=\textwidth]{%
      		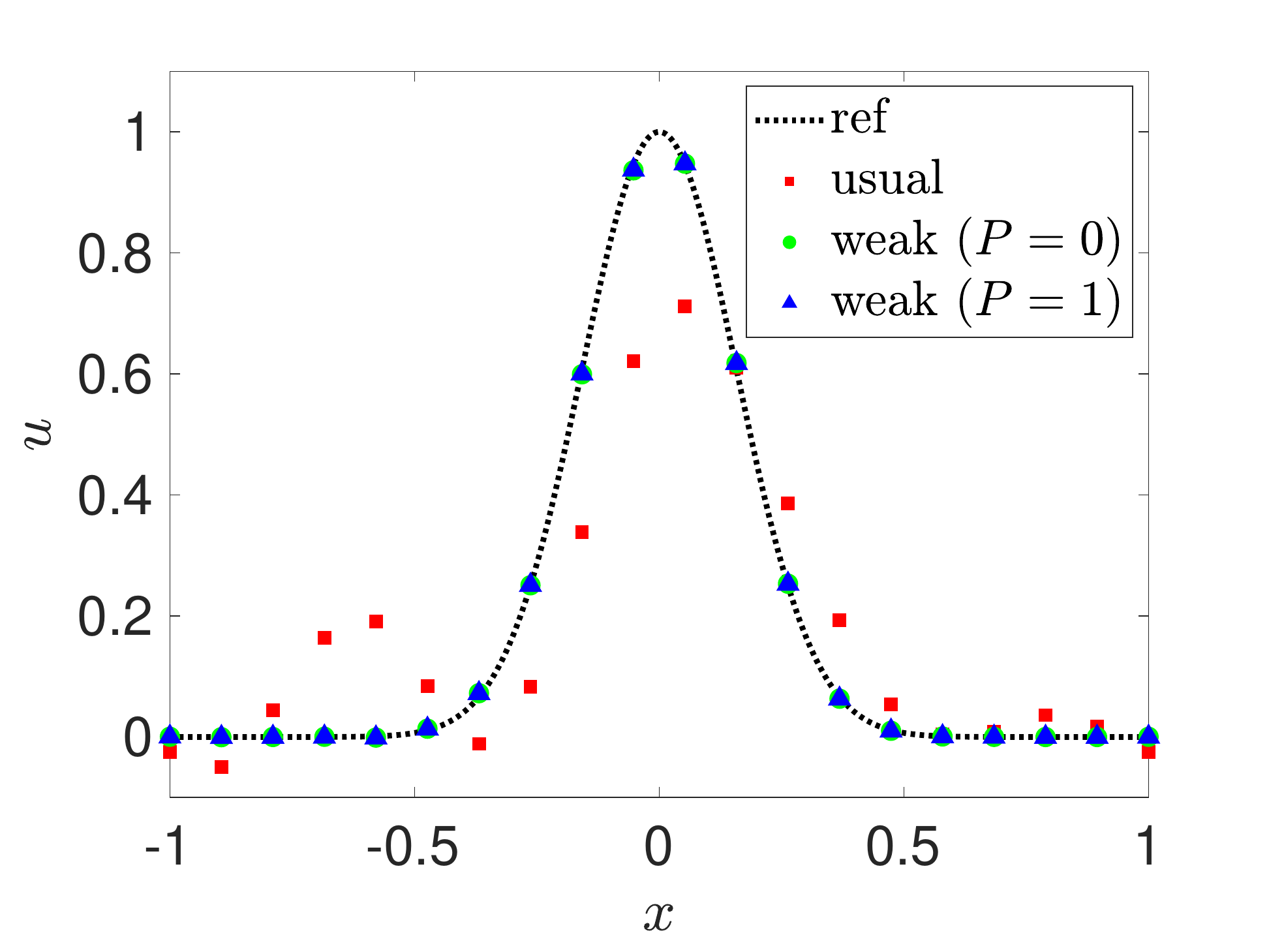}
    		\caption{Solution, multiquadric}
    		\label{fig:linear_1d_A_MQ_sol} 
  	\end{subfigure}%
  	~ 
  	\begin{subfigure}[b]{0.33\textwidth}
    		\includegraphics[width=\textwidth]{%
      		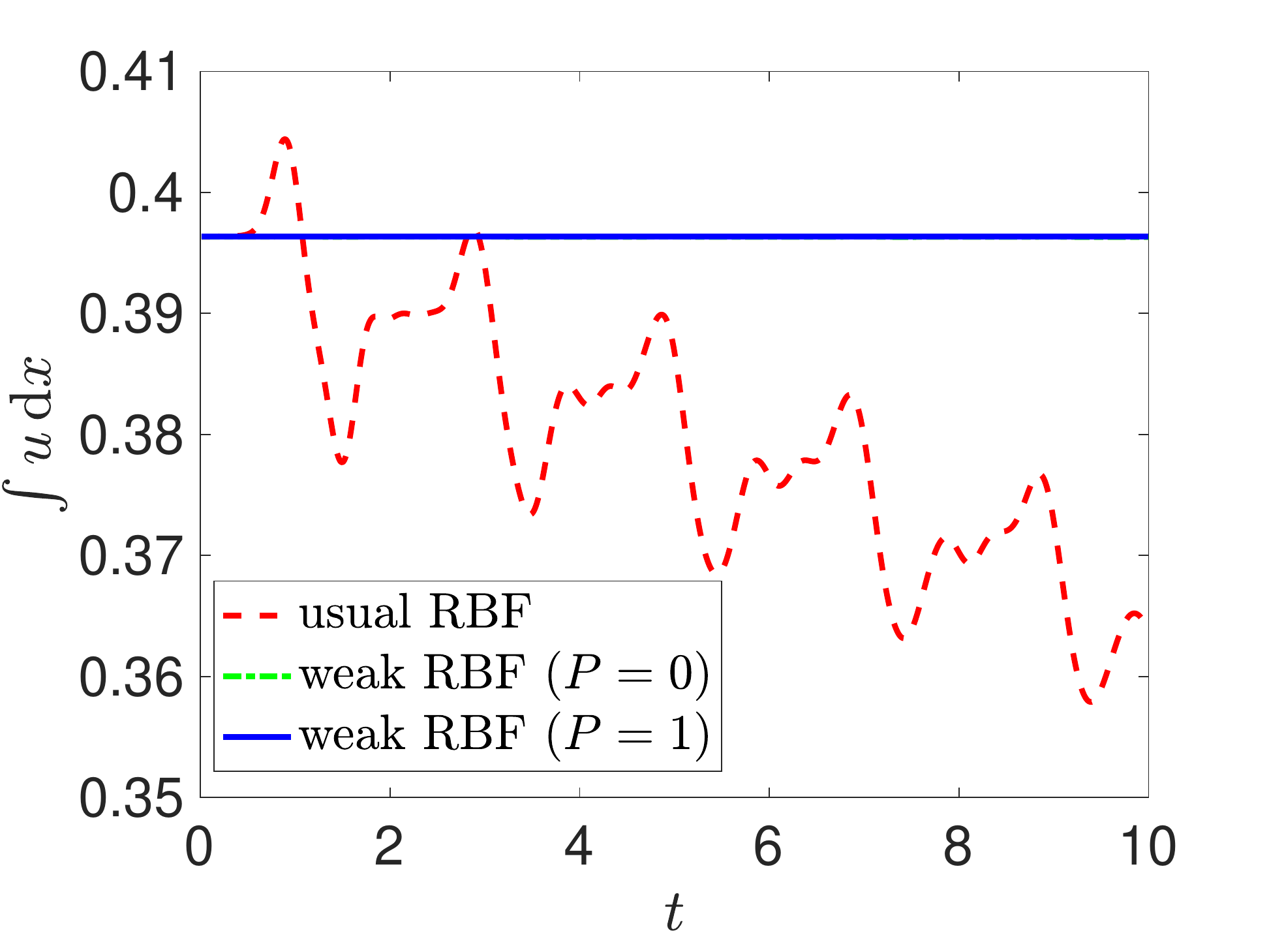}
    		\caption{Momentum, multiquadric}
    		\label{fig:linear_1d_A_MQ_mom} 
  	\end{subfigure}%
	~ 
  	\begin{subfigure}[b]{0.33\textwidth}
    		\includegraphics[width=\textwidth]{%
      		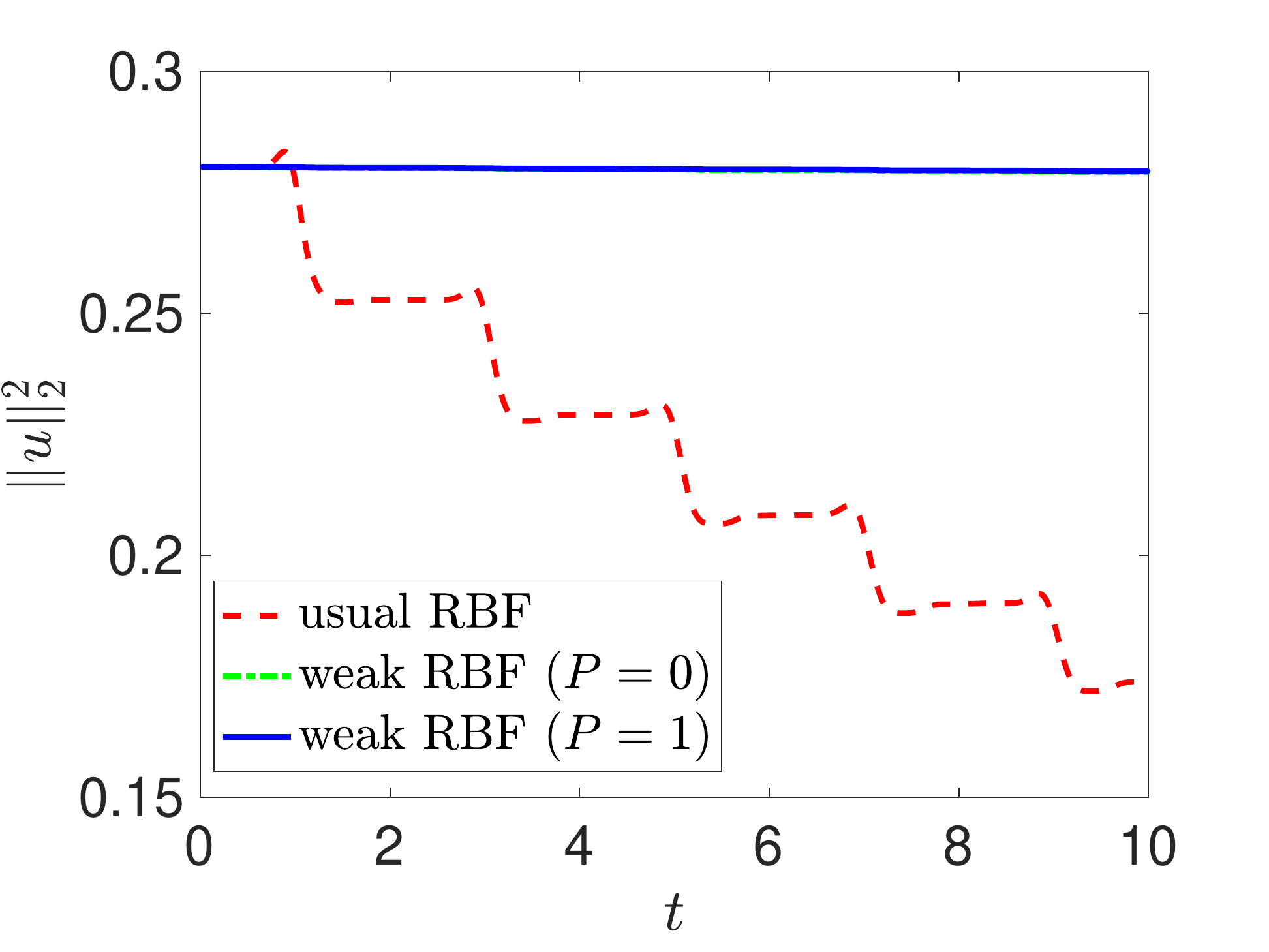}
    		\caption{Energy, multiquadric}
    		\label{fig:linear_1d_A_MQ_en} 
  	\end{subfigure}%
  	\caption{
  		Numerical solutions at $t=10$ (left); their momentum (middle); and energy (right) over time for $u_t + u_x = 0$ with periodic BC \eqref{eq:linear_BC_periodic}. 
  		In all cases, $N = 20$ equidistant nodes and shape parameter $\varepsilon=5$ were used. 
  	}
  	\label{fig:linear_1d_A2}
\end{figure} 

Figures \ref{fig:linear_1d_A} and \ref{fig:linear_1d_A2} illustrate results of the standard as well as the weak RBF method for the linear advection equation with periodic BCs at time $t=10$. 
Different kernels are compared, including the cubic, quintic, Gaussian (G), inverse quadratic (IQ) and multiquadric (MQ) kernel. 
For the latter three a shape parameter of $\varepsilon=5$ was used. 
Furthermore, all tests were performed for $N = 20$ equidistant nodes. 
From Figure \ref{fig:linear_1d_A} it is apparent that in all cases the weak RBF method yields visibly more accurate results than the standard (strong) RBF method. 
In accordance with our previous investigations on conservation and energy stability, we also observe that momentum, $\int u \intd x$, is preserved by the weak RBF method and energy, $\| u \|_2^2$, is nonincreasing.  This is independent of whether $P = 0$ or $1$.
For the the standard RBF method, on the other hand, unphysical profiles for momentum and energy are evident. 
Henceforth we  only focus on the cubic and quintic kernel which allows us to  
 eliminate the potential effects from  poorly chosen shape parameters.

\subsubsection{Error Analysis}

We now provide a more detailed comparison between the standard and weak RBF method for periodic as well as inflow BCs for the cubic and quintic kernel. 

\begin{figure}[!htb]
  	\centering
 	\begin{subfigure}[b]{0.45\textwidth}
    		\includegraphics[width=\textwidth]{%
      		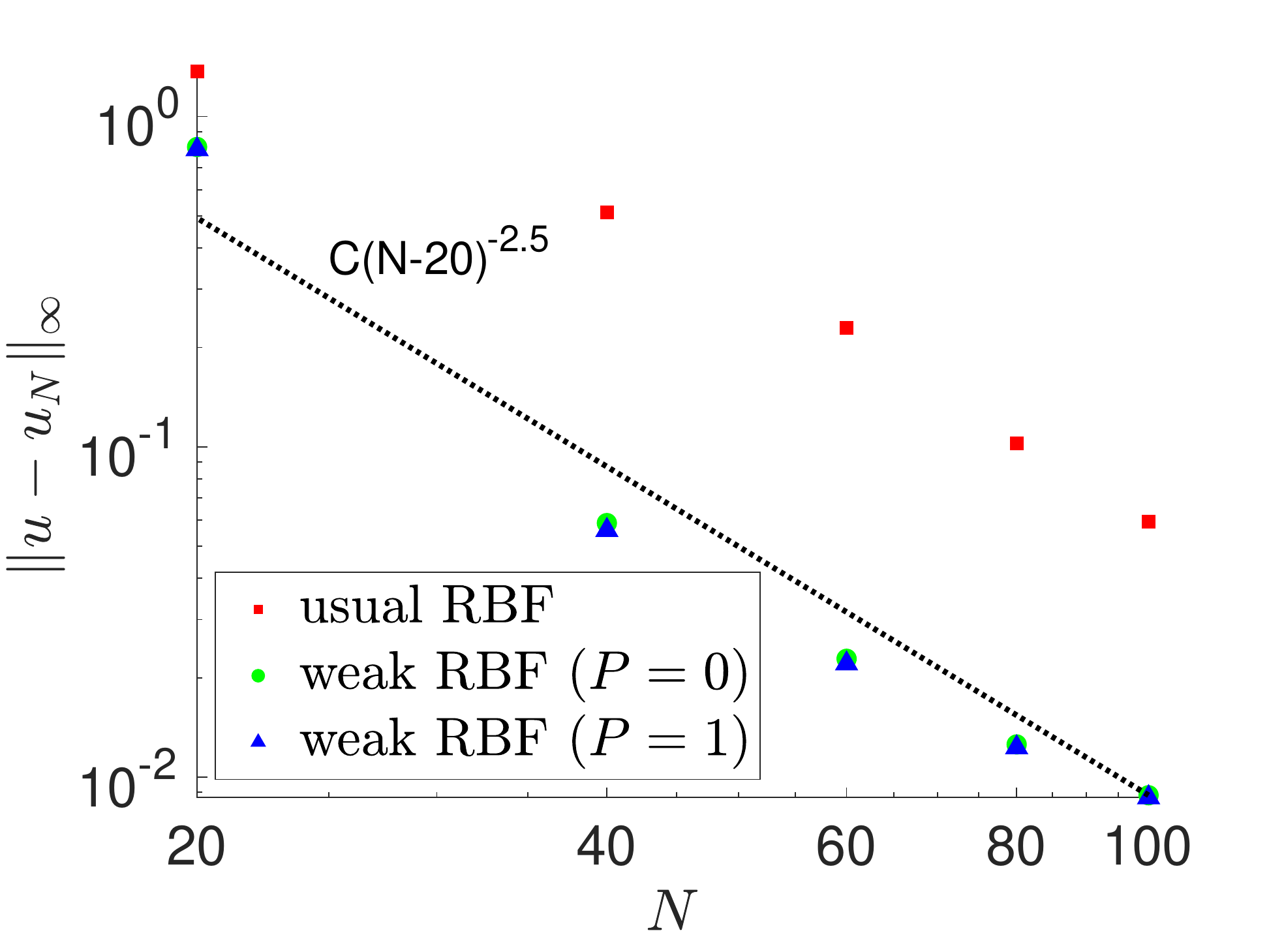}
    		\caption{Cubic, periodic BC}
    		\label{fig:linear_1d_B_max_cubic_periodic} 
  	\end{subfigure}%
  	~ 
  	\begin{subfigure}[b]{0.45\textwidth}
    		\includegraphics[width=\textwidth]{%
      		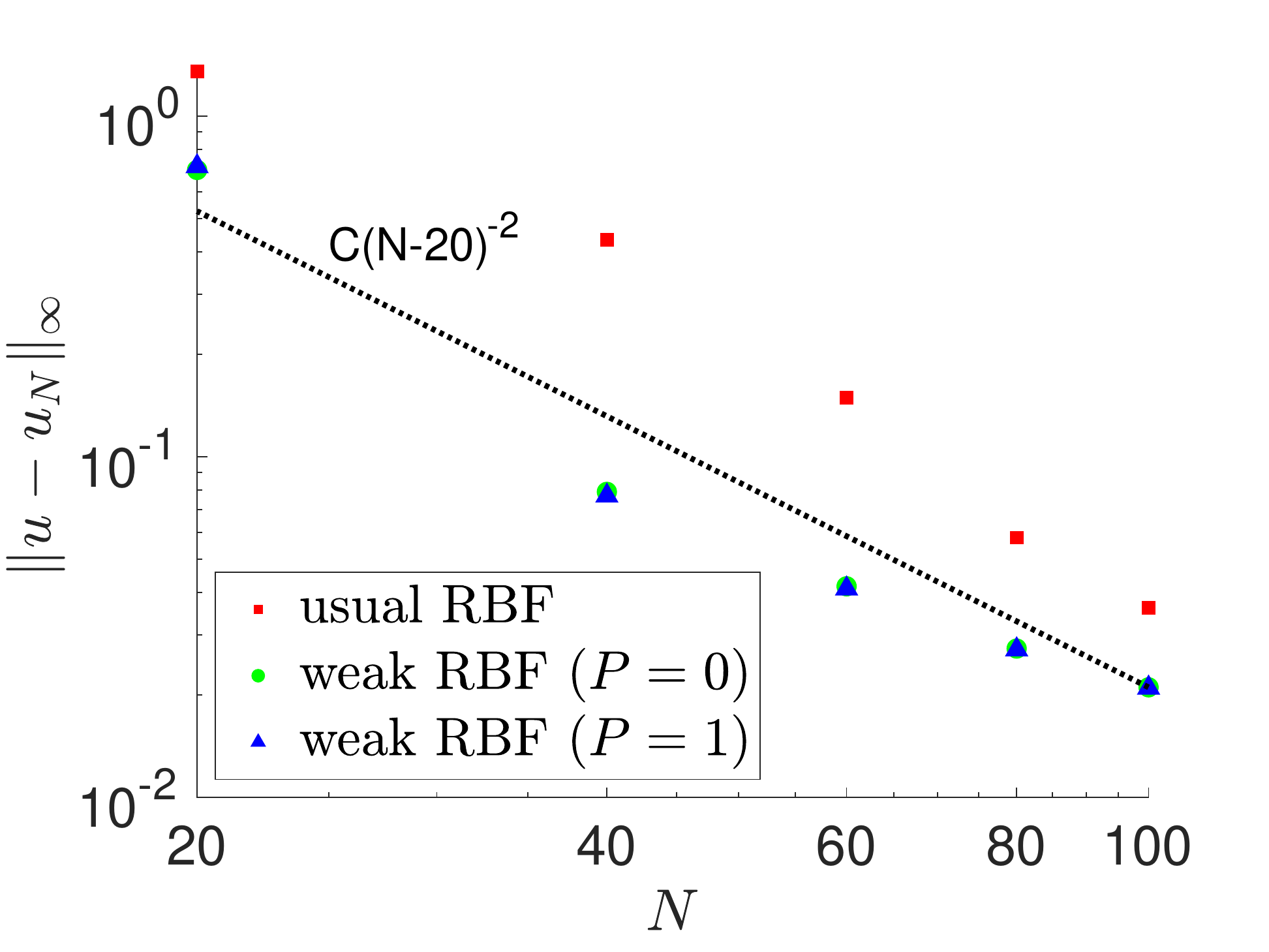}
    		\caption{Cubic, inflow BC}
    		\label{fig:linear_1d_B_max_cubic_inflow}
  	\end{subfigure}%
	\\ 
  	\begin{subfigure}[b]{0.45\textwidth}
    		\includegraphics[width=\textwidth]{%
      		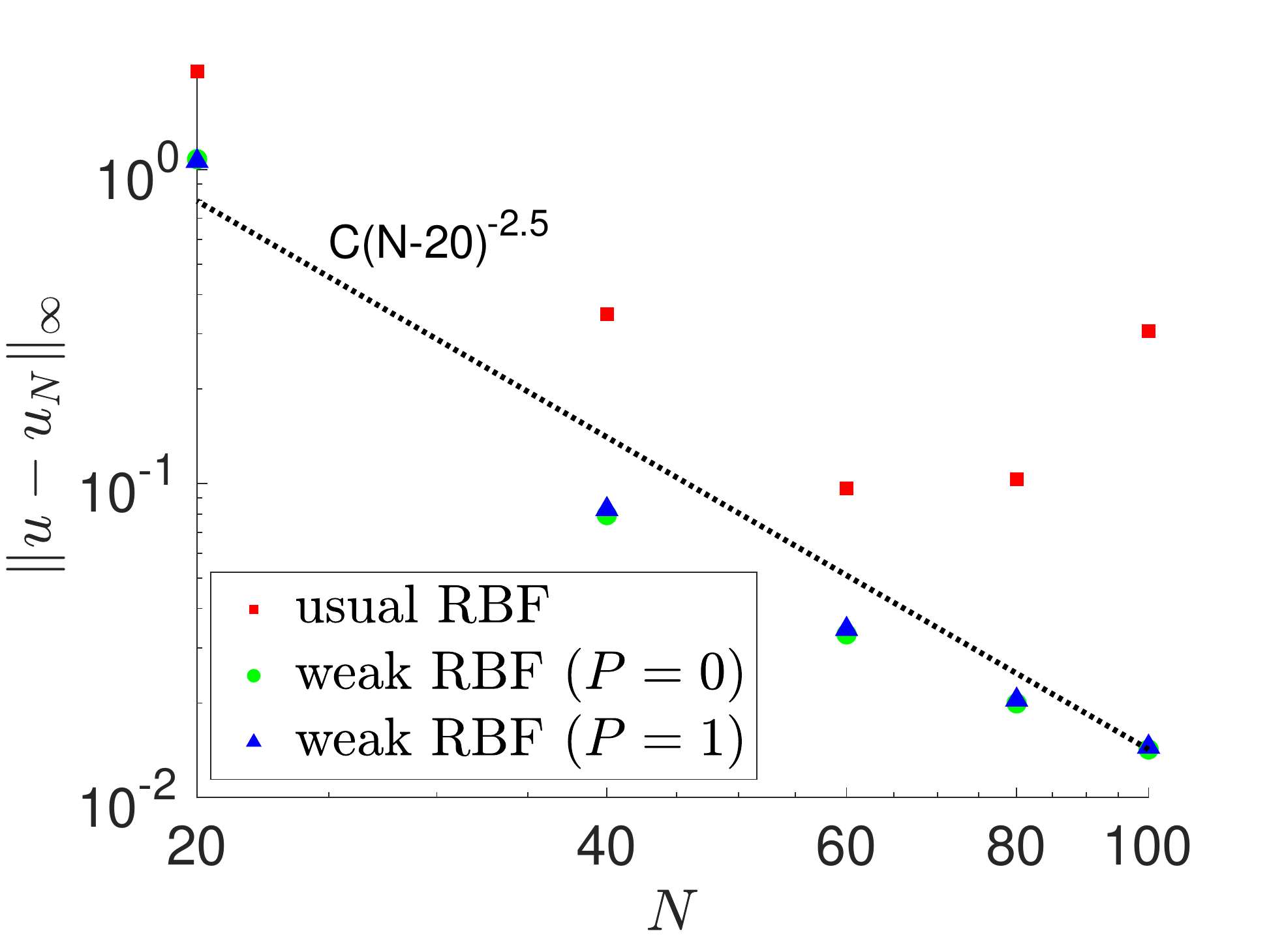}
    		\caption{Quintic, periodic BC}
    		\label{fig:linear_1d_B_max_quintic_periodic} 
  	\end{subfigure}%
  	~ 
  	\begin{subfigure}[b]{0.45\textwidth}
    		\includegraphics[width=\textwidth]{%
      		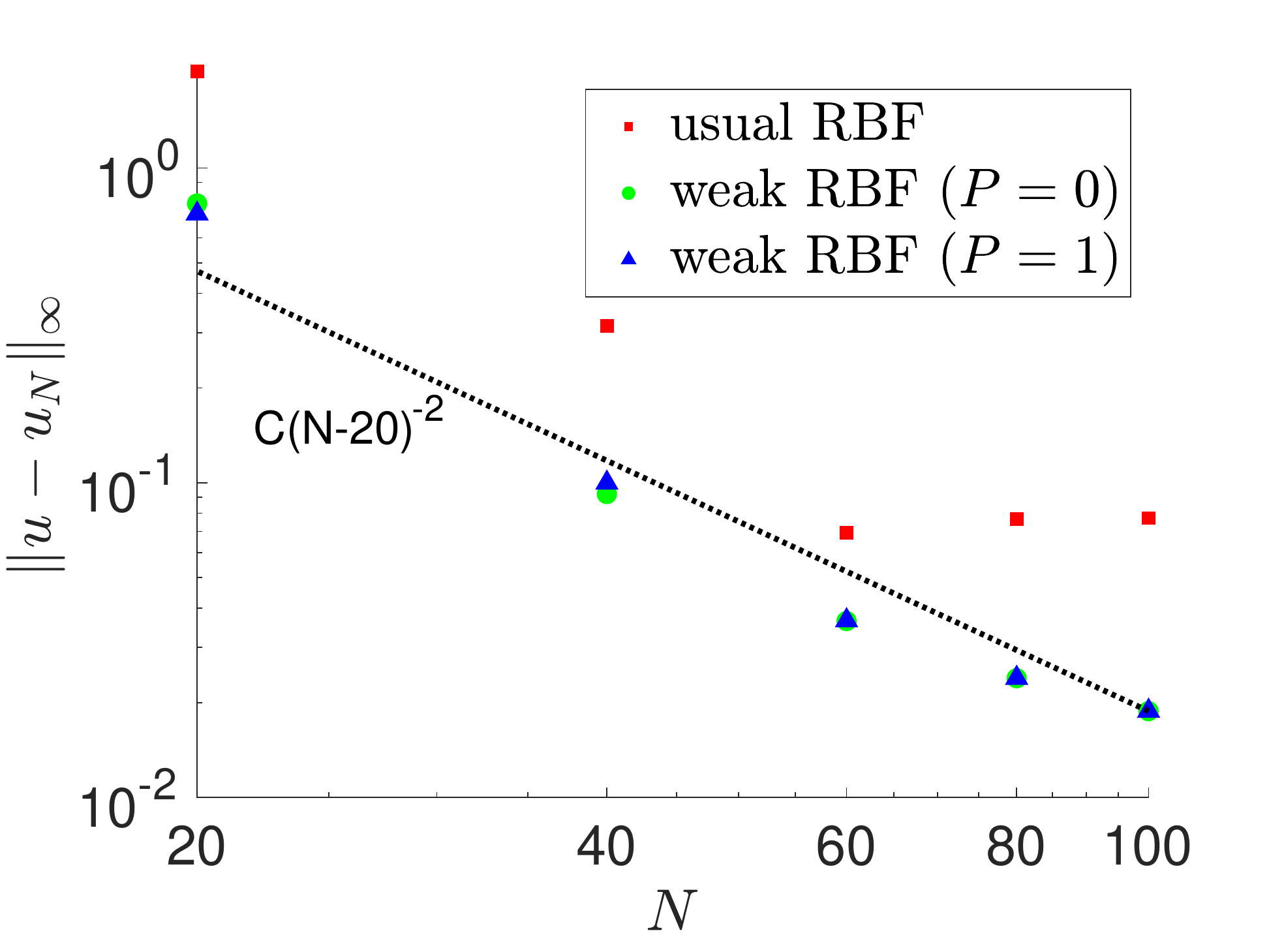}
    		\caption{Quintic, inflow BC}
    		\label{fig:linear_1d_B_max_quintic_inflow}
  	\end{subfigure}%
  	\caption{
  		$\|\cdot\|_{\infty}$-errors of the numerical solutions at $t=10$ for the linear text problem \eqref{eq:linear_1d} with (left) periodic and (right) inflow BC. 
  	}
  	\label{fig:linear_1d_B_max}
\end{figure} 

\begin{figure}[!htb]
  	\centering
 	\begin{subfigure}[b]{0.45\textwidth}
    		\includegraphics[width=\textwidth]{%
      		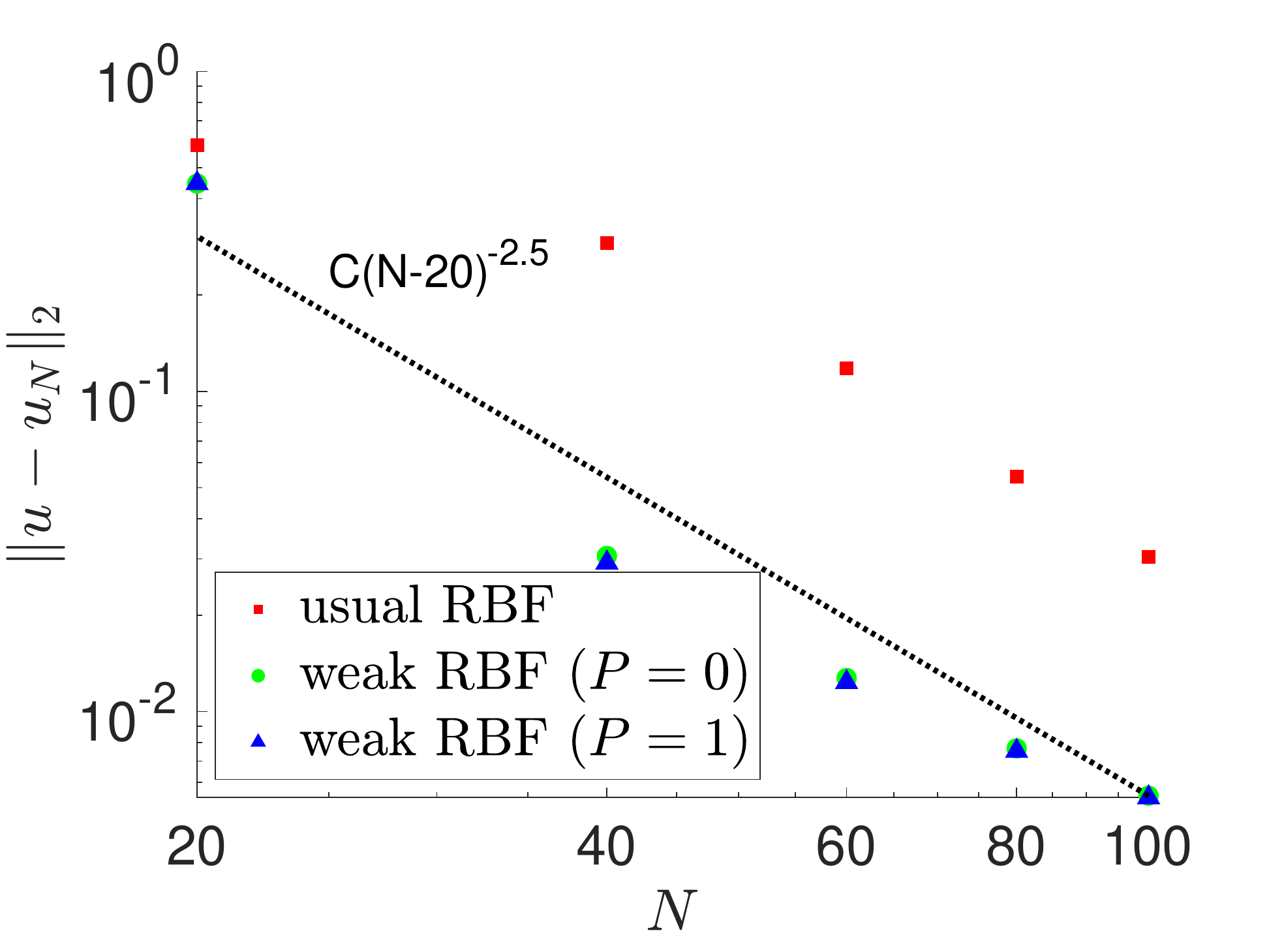}
    		\caption{Cubic, periodic BC}
    		\label{fig:linear_1d_B_L2_cubic_periodic} 
  	\end{subfigure}%
  	~ 
  	\begin{subfigure}[b]{0.45\textwidth}
    		\includegraphics[width=\textwidth]{%
      		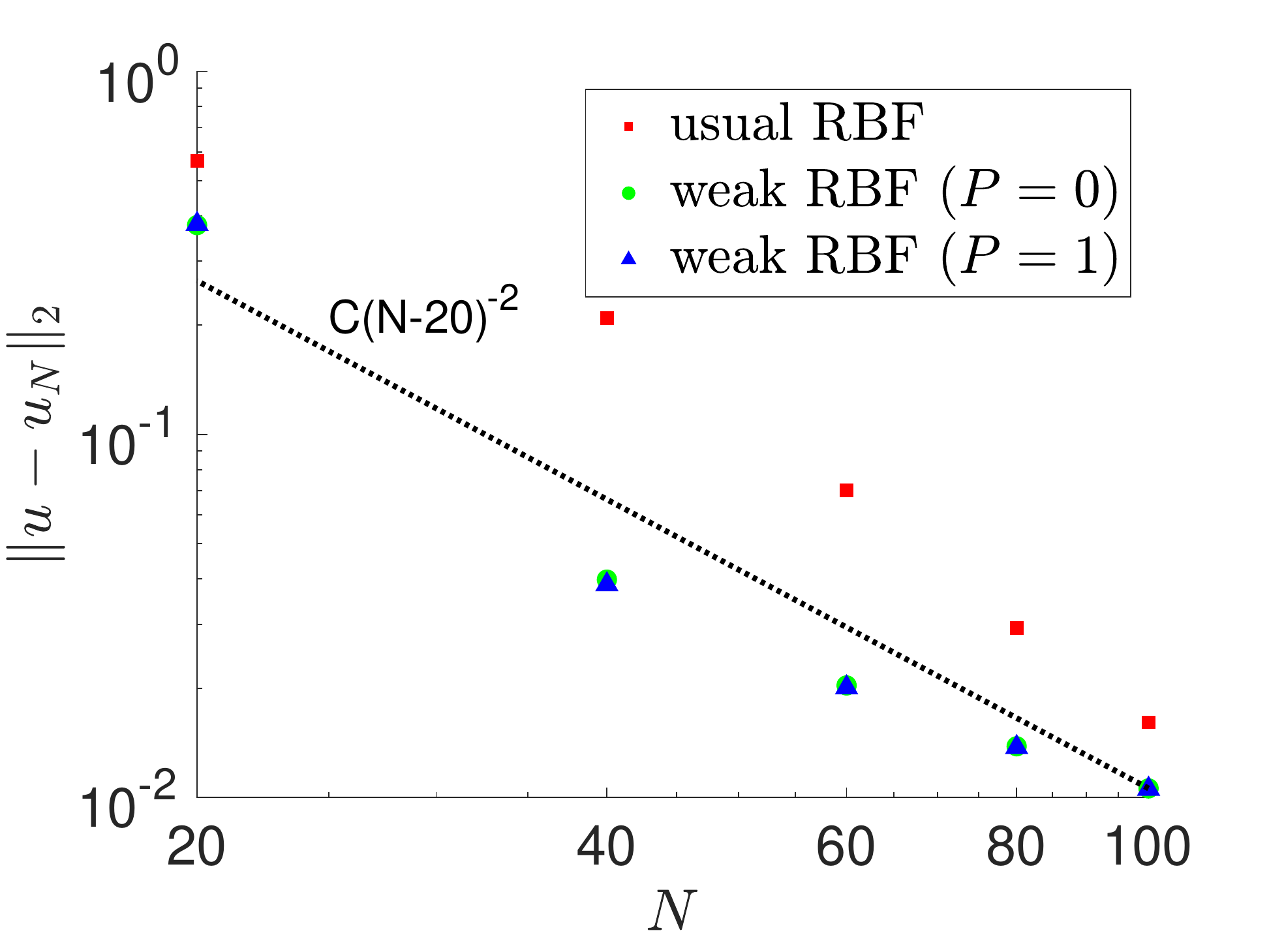}
    		\caption{Cubic, inflow BC}
    		\label{fig:linear_1d_B_L2_cubic_inflow}
  	\end{subfigure}%
	\\ 
  	\begin{subfigure}[b]{0.45\textwidth}
    		\includegraphics[width=\textwidth]{%
      		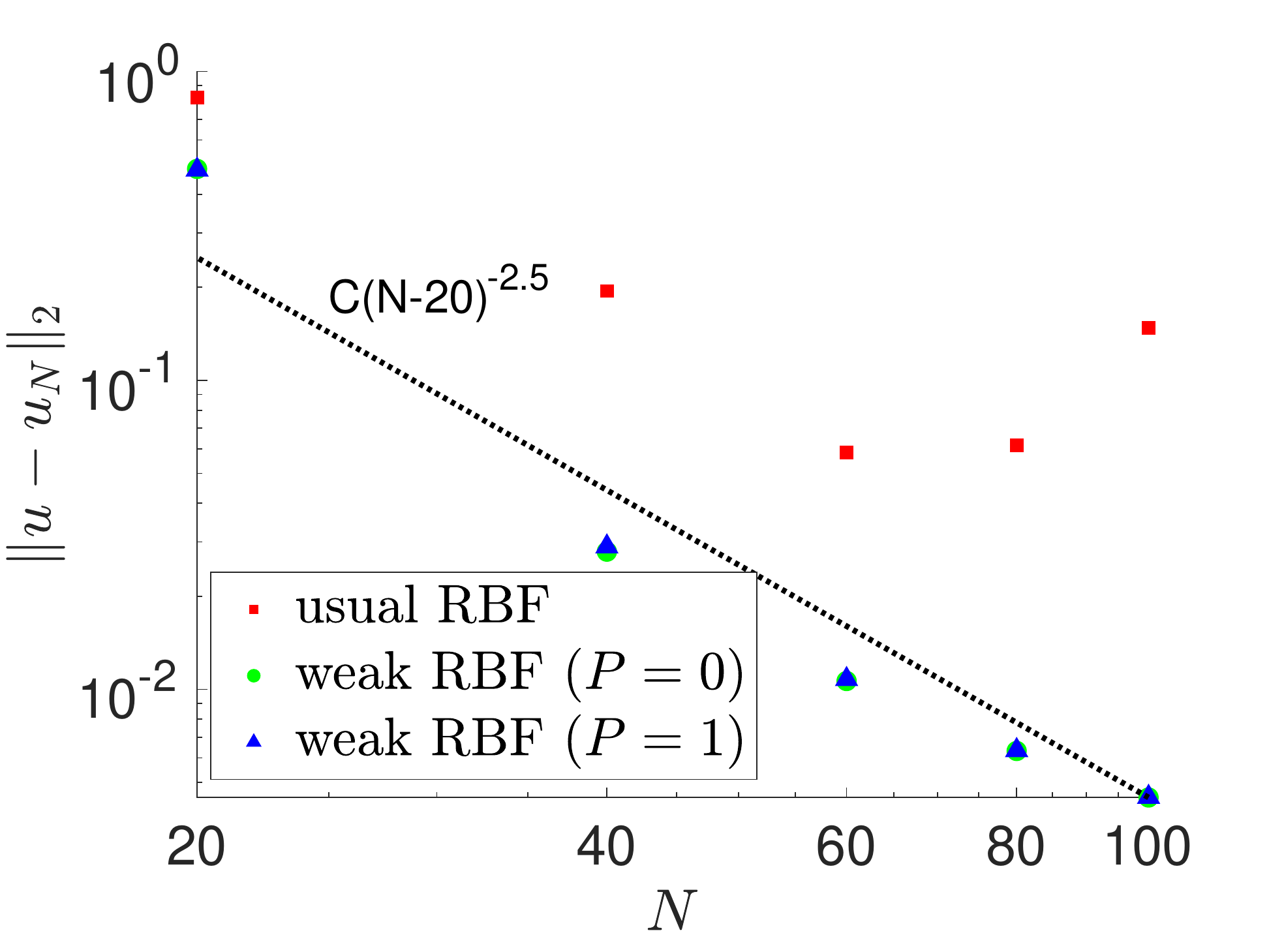}
    		\caption{Quintic, periodic BC}
    		\label{fig:linear_1d_B_L2_quintic_periodic} 
  	\end{subfigure}%
  	~ 
  	\begin{subfigure}[b]{0.45\textwidth}
    		\includegraphics[width=\textwidth]{%
      		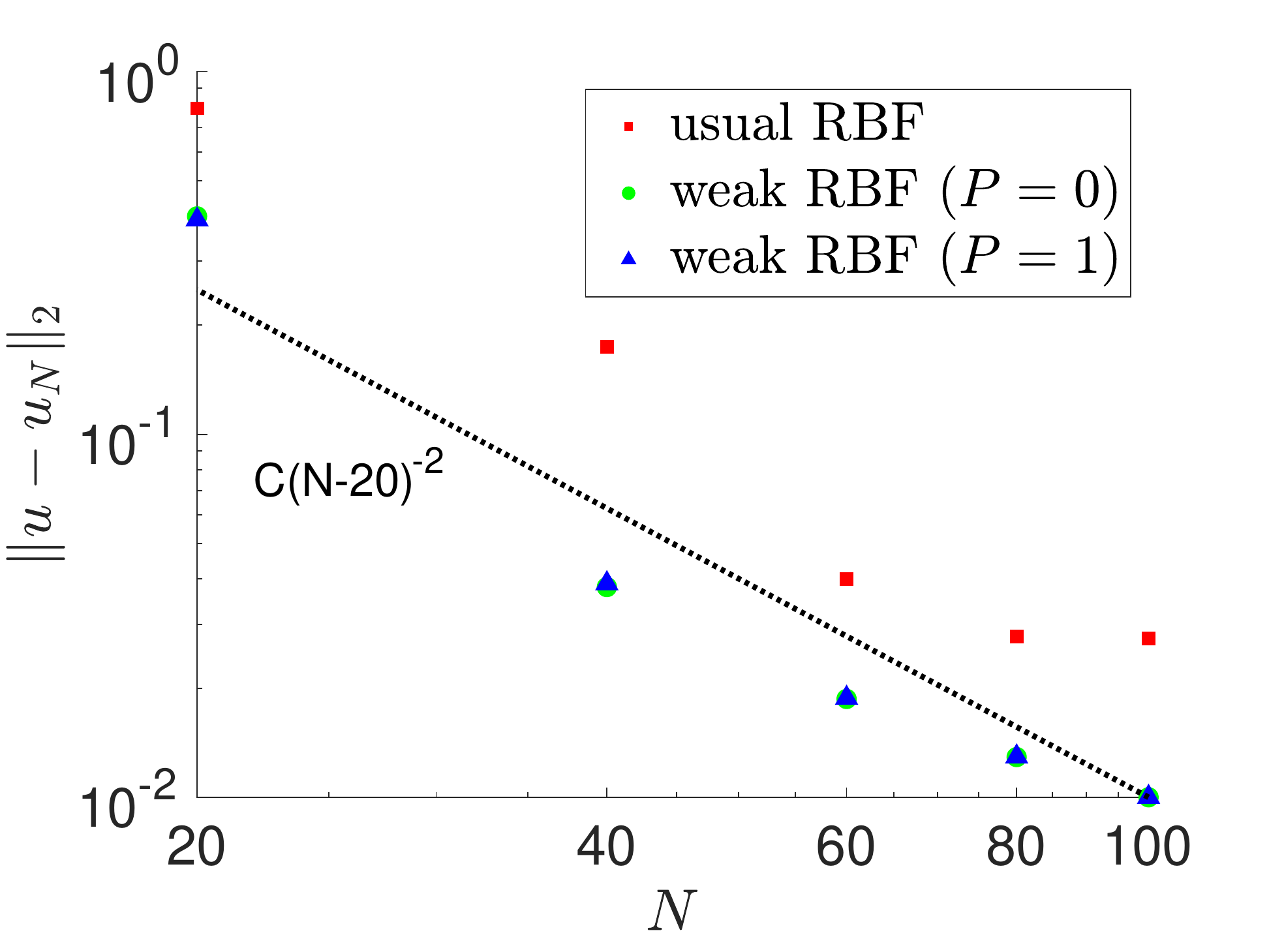}
    		\caption{Quintic, inflow BC}
    		\label{fig:linear_1d_B_L2_quintic_inflow}
  	\end{subfigure}%
  	\caption{
  		$\|\cdot\|_{2}$-errors of the numerical solutions at $t=10$ for the linear text problem \eqref{eq:linear_1d} with (left) periodic and (right) inflow BC}
  	\label{fig:linear_1d_B_L2}
\end{figure} 

Figures \ref{fig:linear_1d_B_max} and \ref{fig:linear_1d_B_L2} illustrate the $\|\cdot\|_{\infty}$- and $\|\cdot\|_{2}$-errors of both methods corresponding to the linear IVP
\begin{equation}\label{eq:linear_1d_error}
	u_t + u_x = 0, \quad 
	u(0,x) = \cos^2 \left( 4 \pi x \right)
\end{equation} 
with $x \in \Omega = [-1,1]$ and periodic as well as inflow BC at $t=2$. 
These error norms are respectively given by 
\begin{equation} 
\begin{aligned}
	\| u - u_N \|_{\infty} = \max_{n=1,\dots,N} | u(x_n) - u_N(x_n) |, \quad
	\| u - u_N \|_{2} = \sqrt{ \sum_{n=1}^N | u(x_n) - u_N(x_n) |^2 }, 
\end{aligned}
\end{equation}
where $u$ denotes the exact solution, $u_N$ the numerical solution, and $x_1,\dots,x_N$ are the nodes. 
Figures \ref{fig:linear_1d_B_max} and \ref{fig:linear_1d_B_L2} consider the errors using  equidistant nodes. 
It is clear that the weak RBF method yields more accurate results than the standard RBF method in all cases, and 
that the standard RBF method does not even seem to converge for the quintic kernel case. 
This may be due to rising instability in combination with the resulting numerical artifacts never leaving the computational domain in case of periodic BCs. 
The weak RBF method, on the other hand, is observed to have a convergence rate of $2.5$ in the periodic case, regardless of whether the cubic or quintic kernel is used. 
We note that the local approximation orders of the cubic and quintic kernel are respectively $2$ and $3$, \cite{iske2003radial,iske2020ten}. 
For the inflow BC, the rate of convergence of the weak RBF method is observed to decrease to $2$ for both kernels. 
Moreover, for the inflow BC, the standard RBF method displays a similar rate of convergence. 
It might be that this increase of stability (and therefore accuracy) for the standard RBF method is related to numerical artifacts being allowed to leave the computational domain while only exact information (due to the inflow BC) flows into the computational domain from the left hand side.
This behavior will be considered more in future investigations.

\subsubsection{Equidistant vs Nonequidistant Points}

As demonstrated in \S \ref{sub:ana-RBF-method} and \S \ref{sub:coll-RBF-method}, neither conservation nor energy stability of the weak RBF method depend on the choice of the nodes.  
However accuracy of the weak RBF method might suffer from poor distributions of the nodes.  Below we further investigate the potential implication of different nodal distributions.

\begin{figure}[!htb]
  	\centering
 	\begin{subfigure}[b]{0.33\textwidth}
    		\includegraphics[width=\textwidth]{%
      		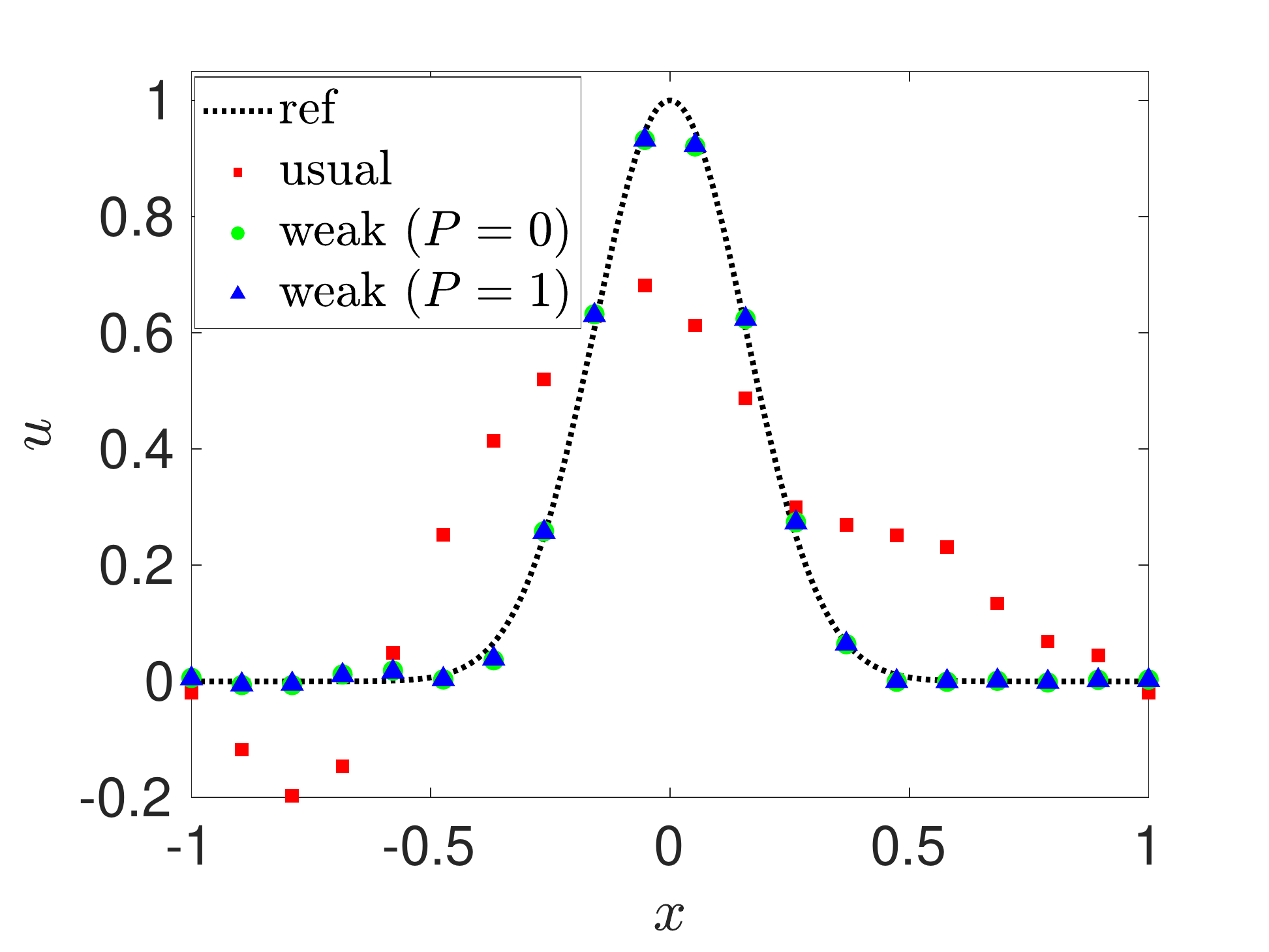}
    		\caption{Solution, equidistant}
    		\label{fig:linear_1d_C_sol_equid_cubic} 
  	\end{subfigure}%
  	~ 
  	\begin{subfigure}[b]{0.33\textwidth}
    		\includegraphics[width=\textwidth]{%
      		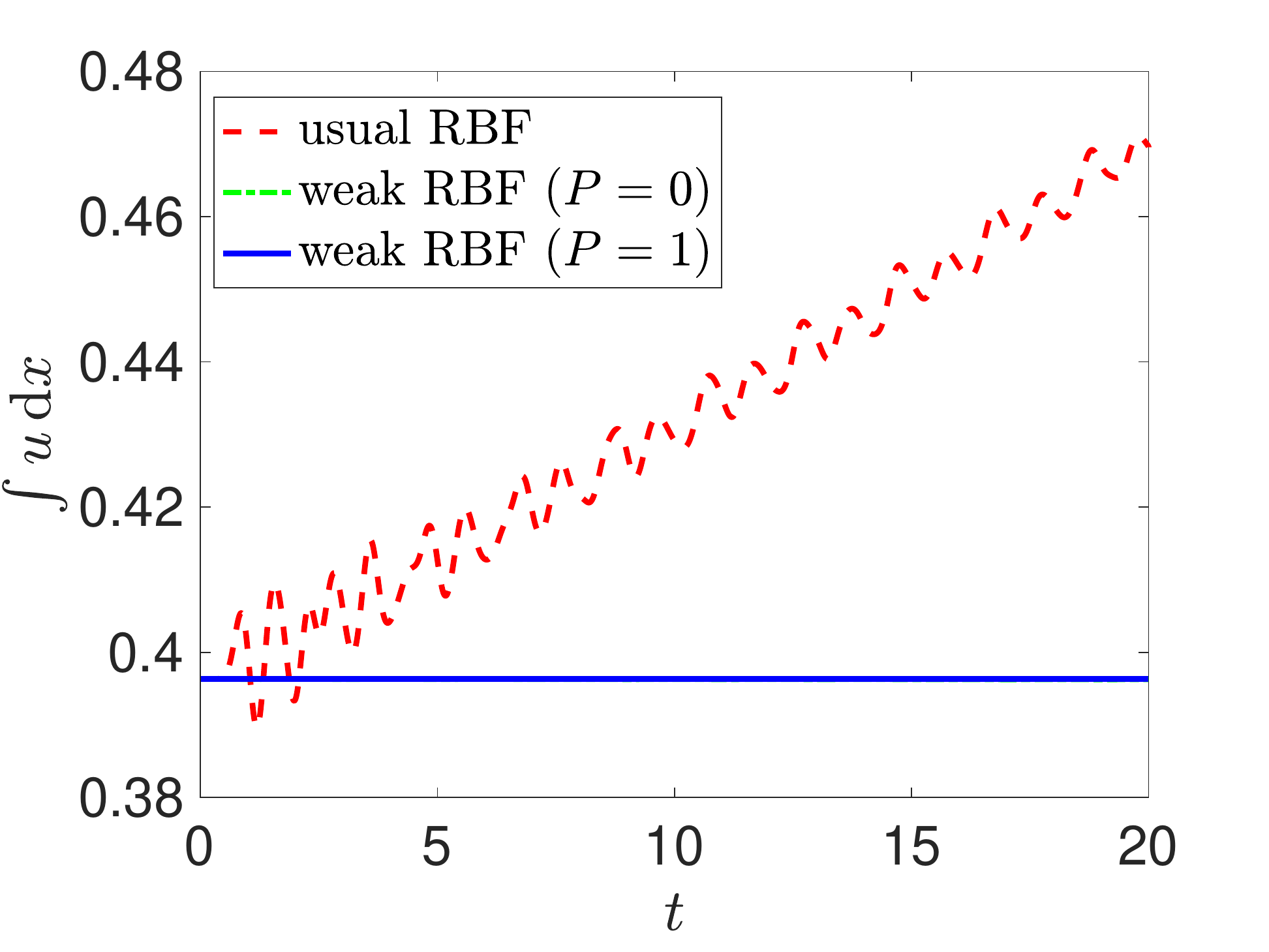}
    		\caption{Momentum, equidistant}
    		\label{fig:linear_1d_C_momentum_equid_cubic} 
  	\end{subfigure}%
	~ 
  	\begin{subfigure}[b]{0.33\textwidth}
    		\includegraphics[width=\textwidth]{%
      		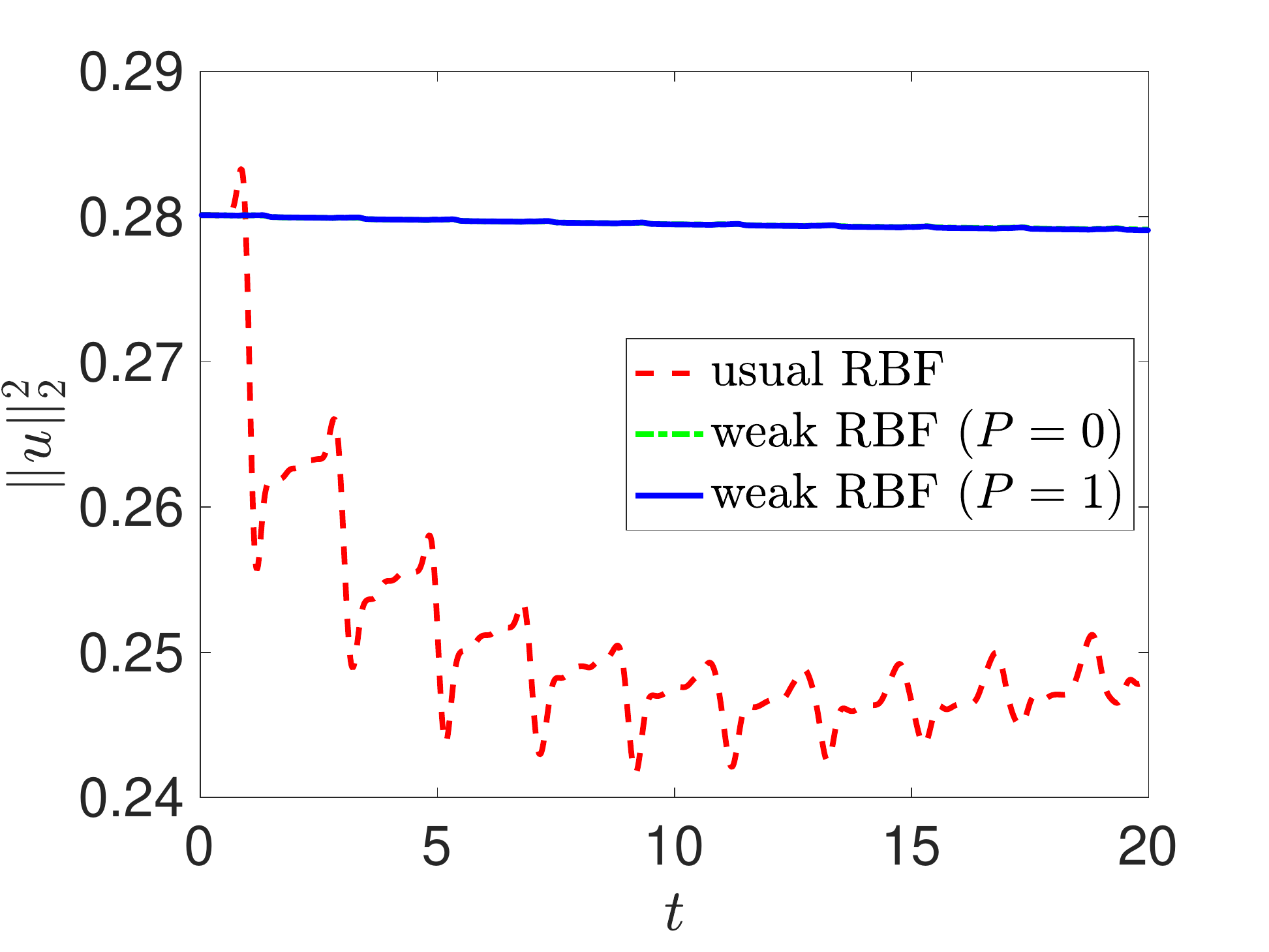}
    		\caption{Energy, equidistant}
    		\label{fig:linear_1d_C_energy_equid_cubic} 
  	\end{subfigure}%
	\\ 
  	\begin{subfigure}[b]{0.33\textwidth}
    		\includegraphics[width=\textwidth]{%
      		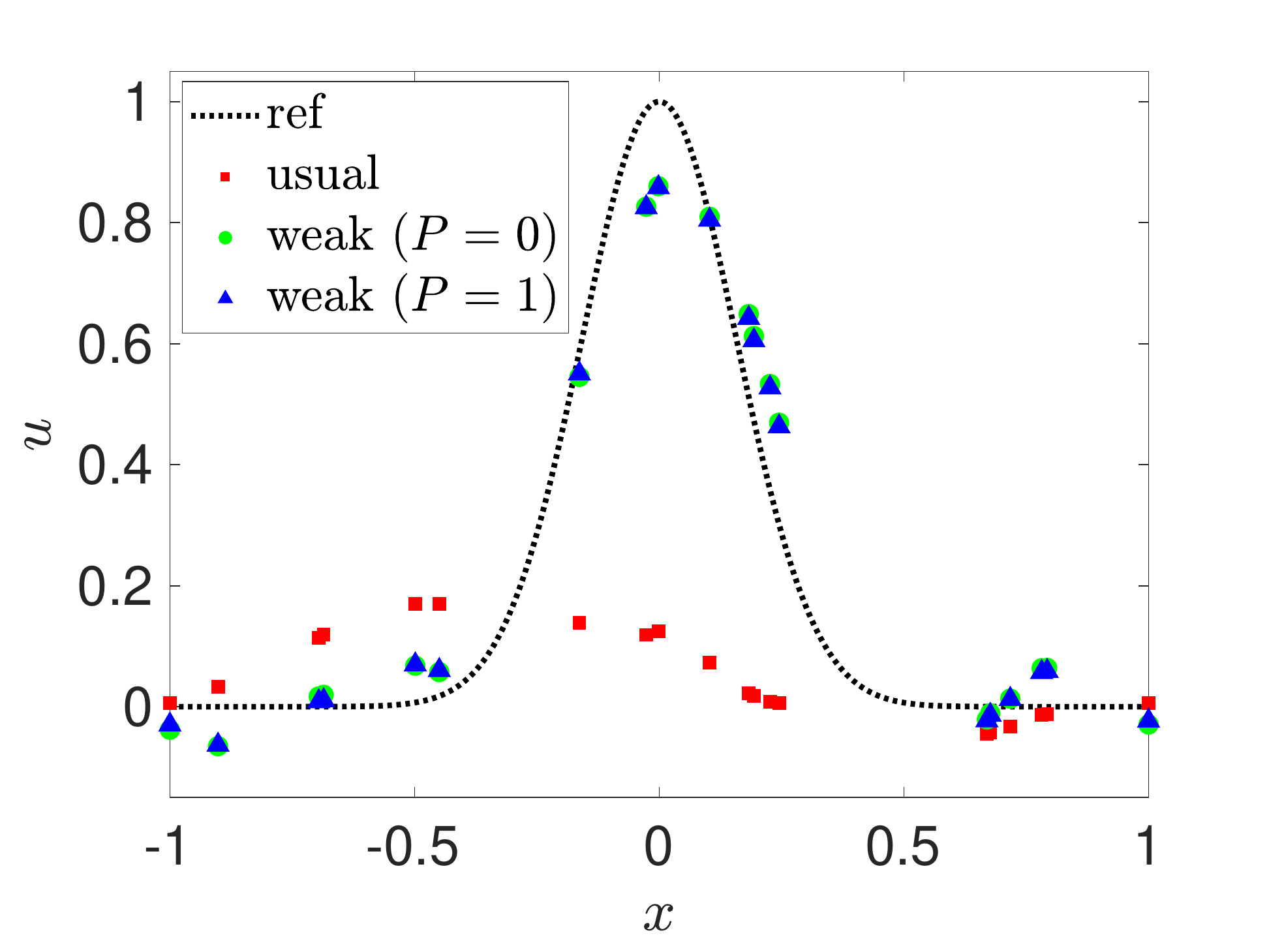}
    		\caption{Solution, random}
    		\label{fig:linear_1d_C_sol_random_cubic} 
  	\end{subfigure}%
  	~ 
  	\begin{subfigure}[b]{0.33\textwidth}
    		\includegraphics[width=\textwidth]{%
      		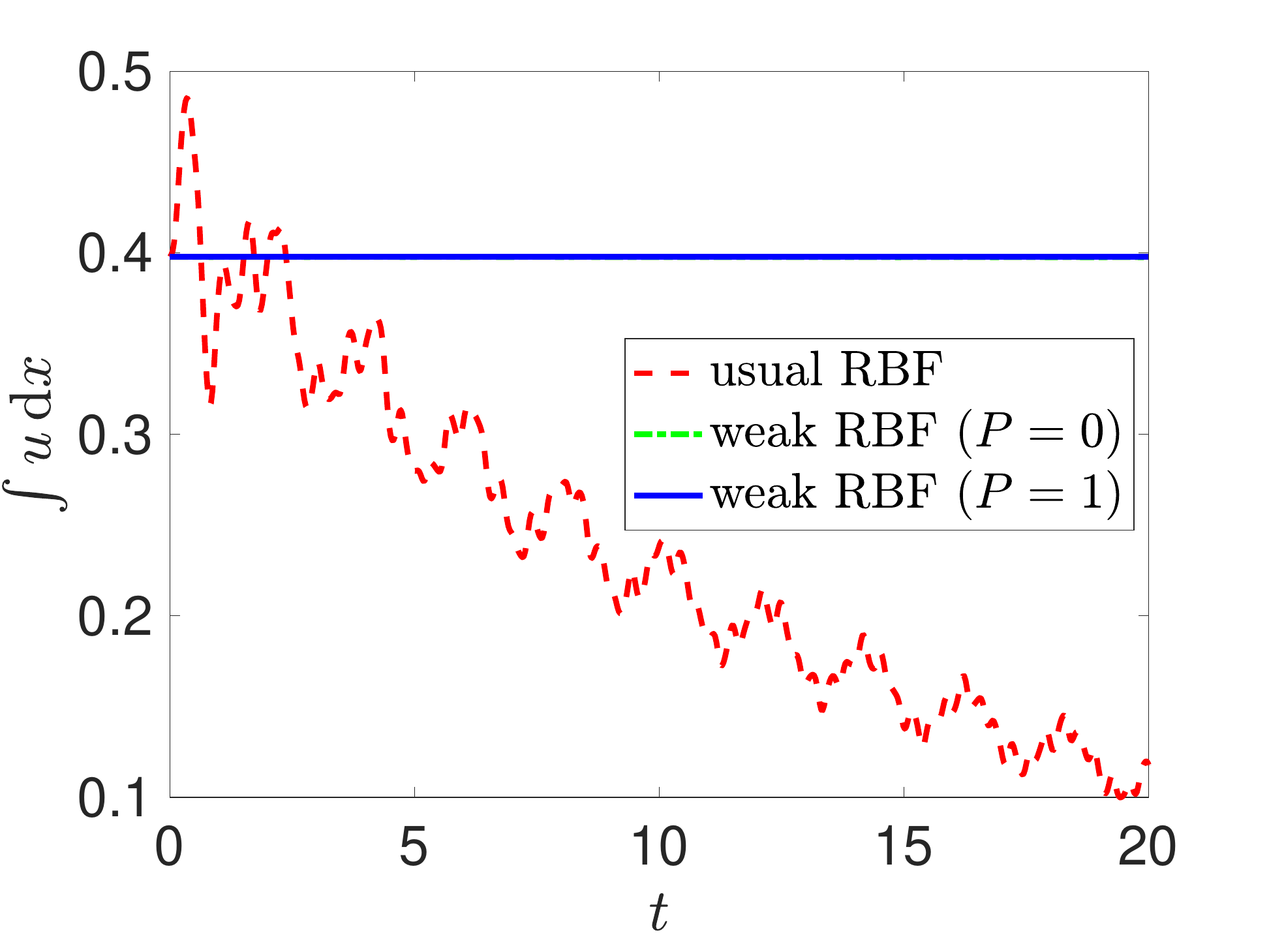}
    		\caption{Momentum, random}
    		\label{fig:linear_1d_C_momentum_random_cubic} 
  	\end{subfigure}%
	~ 
  	\begin{subfigure}[b]{0.33\textwidth}
    		\includegraphics[width=\textwidth]{%
      		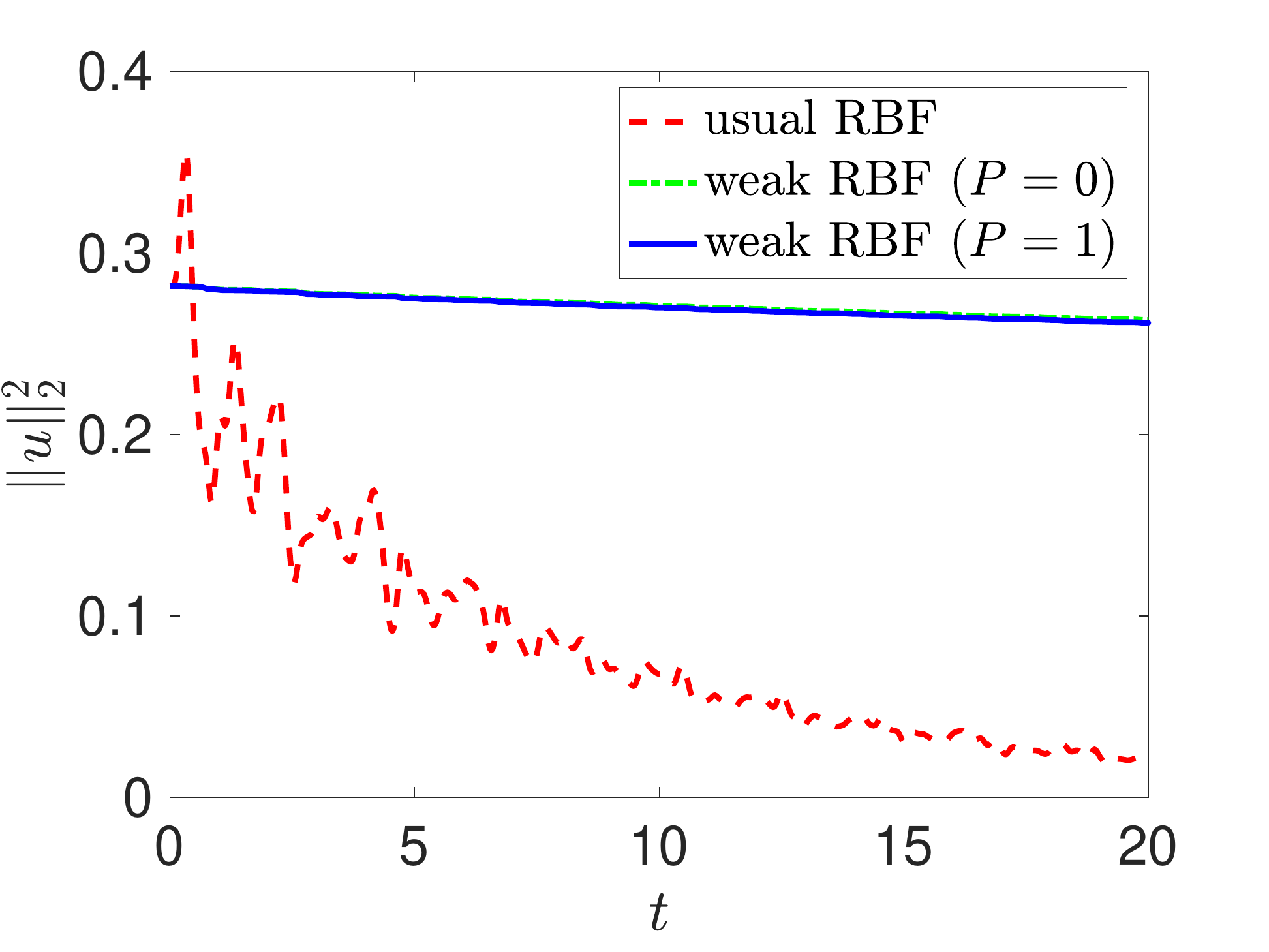}
    		\caption{Energy, random}
    		\label{fig:linear_1d_C_energy_random_cubic} 
  	\end{subfigure}%
  	\caption{
  		Numerical solutions at $t=20$ (left); momentum (middle); and energy (right) over time for $u_t + u_x = 0$ with periodic BCs. 
  		$N = 20$ equidistant (top) and randomly uniformly distributed (bottom) nodes are compared for a cubic kernel.
  	}
  	\label{fig:linear_1d_C_cubic}
\end{figure} 

\begin{figure}[!htb]
  	\centering
 	\begin{subfigure}[b]{0.33\textwidth}
    		\includegraphics[width=\textwidth]{%
      		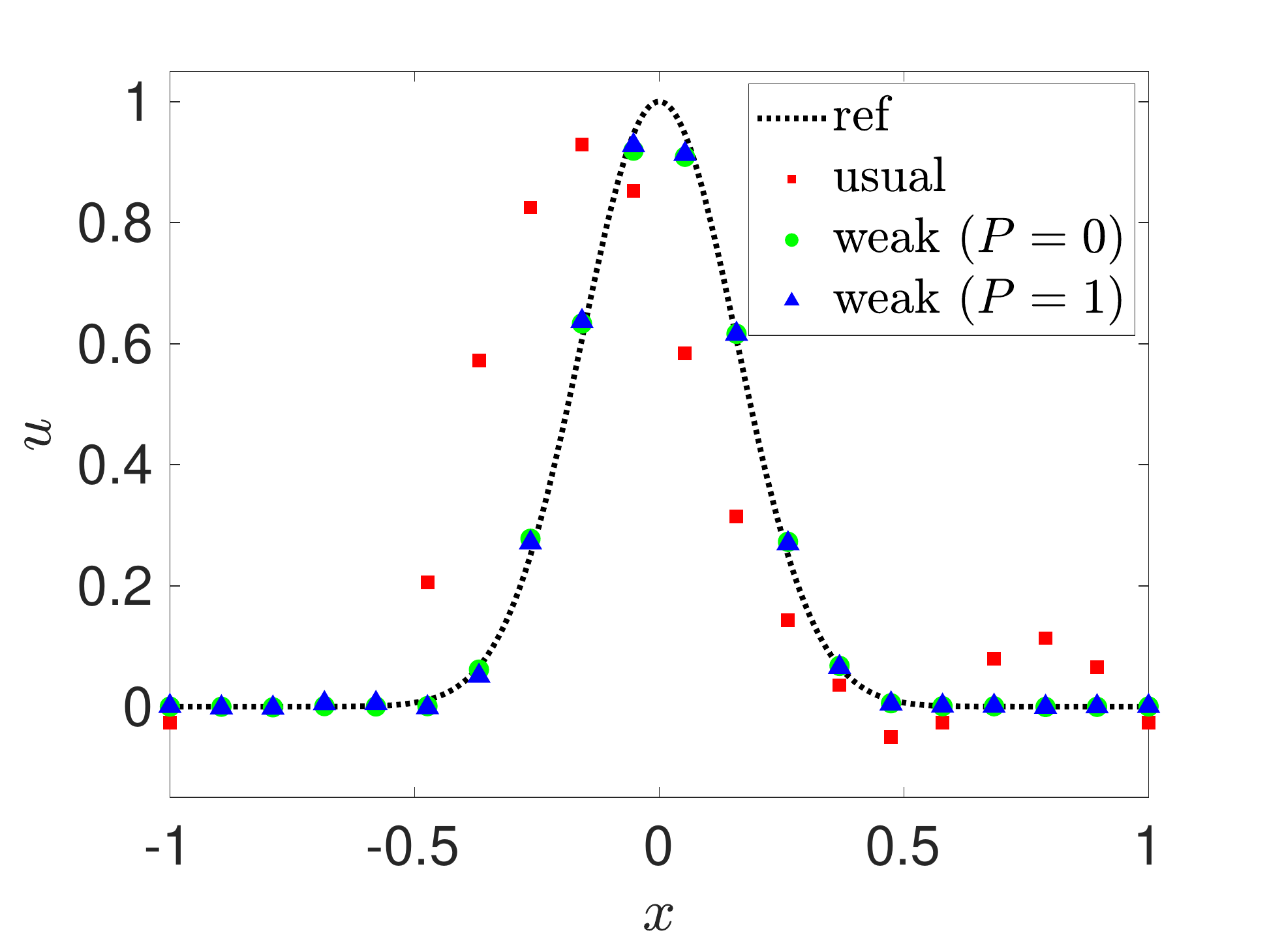}
    		\caption{Solution, equidistant}
    		\label{fig:linear_1d_C_sol_equid_quintic} 
  	\end{subfigure}%
  	~ 
  	\begin{subfigure}[b]{0.33\textwidth}
    		\includegraphics[width=\textwidth]{%
      		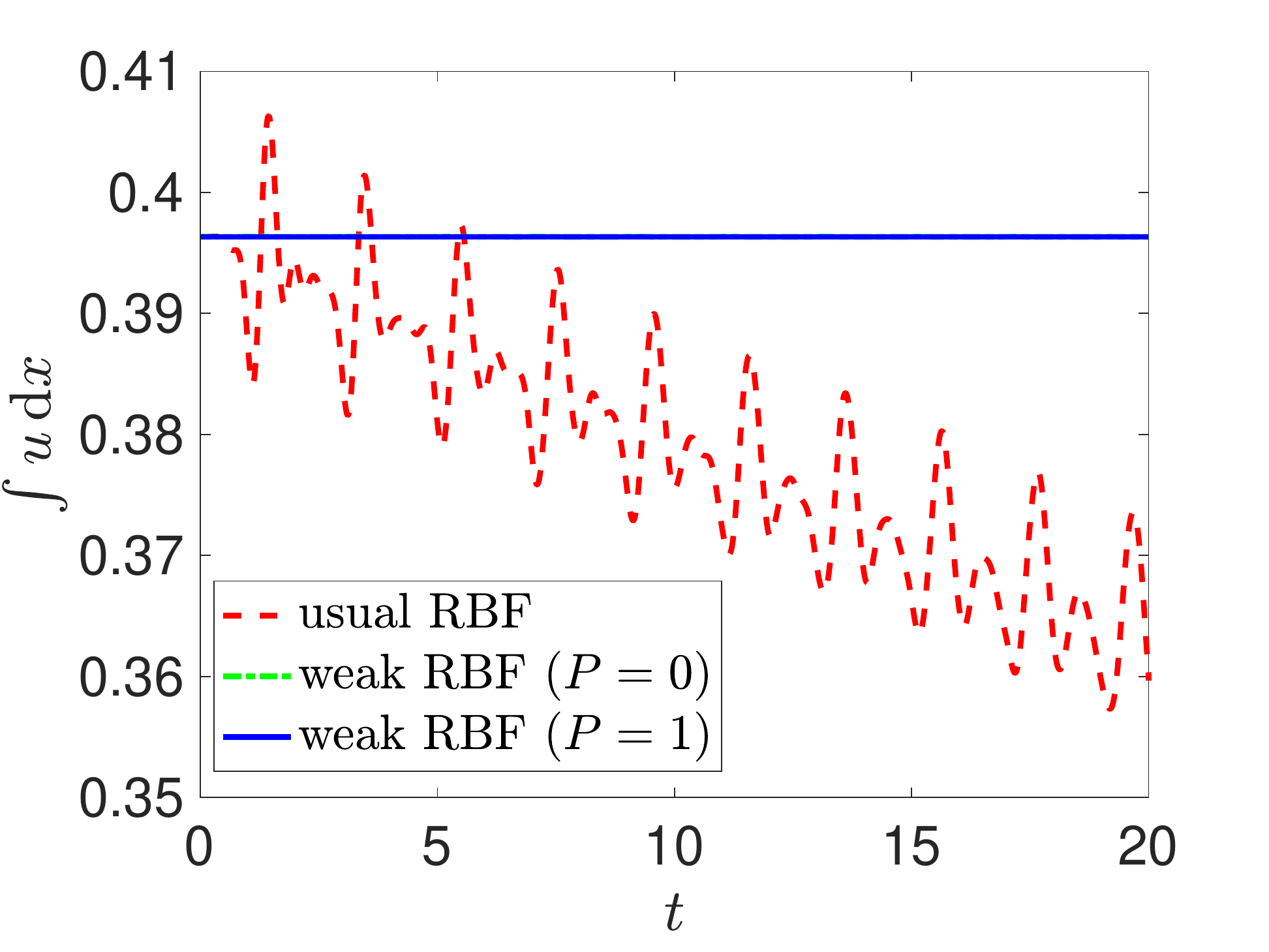}
    		\caption{Momentum, equidistant}
    		\label{fig:linear_1d_C_momentum_equid_quintic} 
  	\end{subfigure}%
	~ 
  	\begin{subfigure}[b]{0.33\textwidth}
    		\includegraphics[width=\textwidth]{%
      		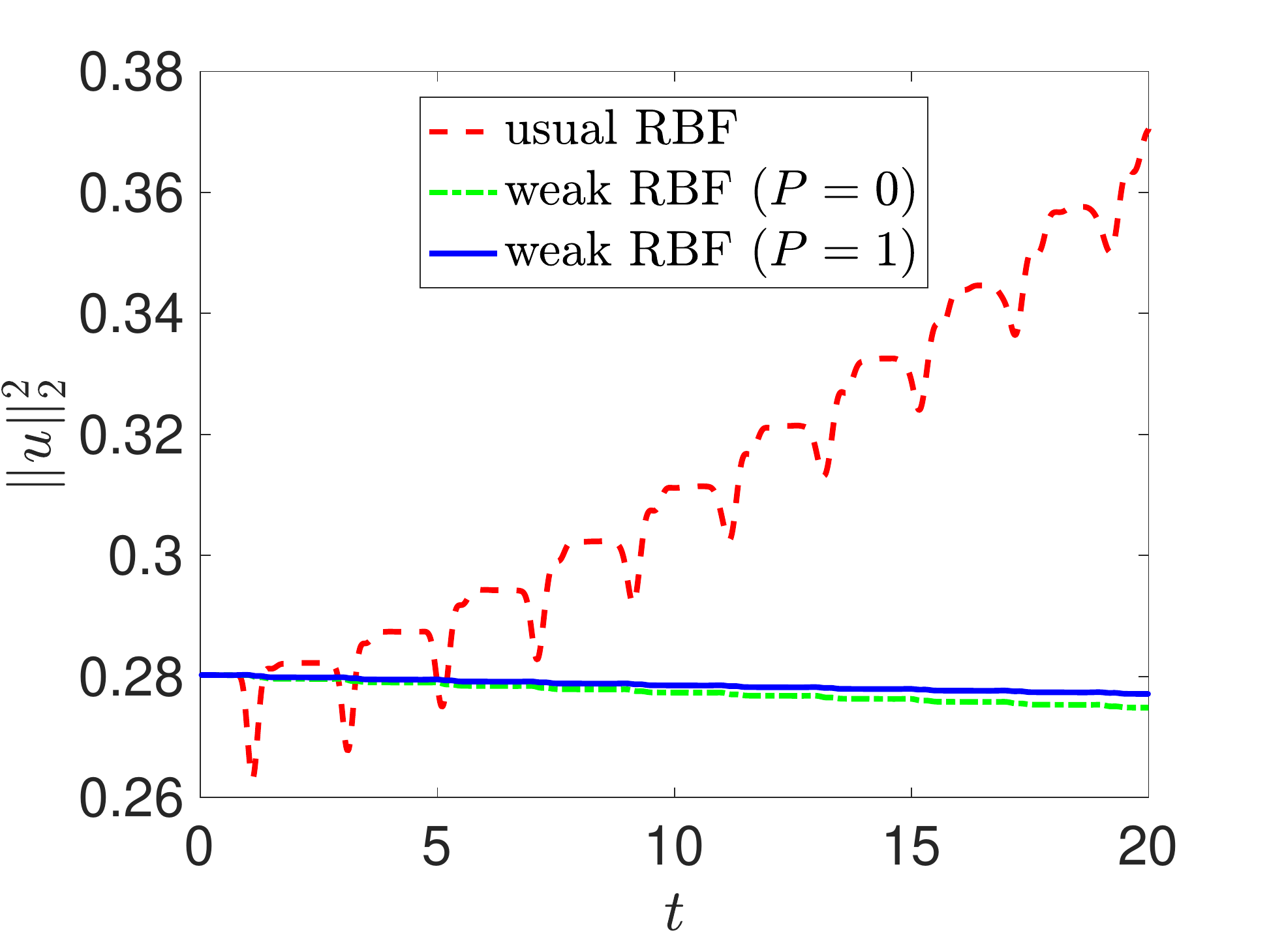}
    		\caption{Energy, equidistant}
    		\label{fig:linear_1d_C_energy_equid_quintic} 
  	\end{subfigure}%
	\\ 
  	\begin{subfigure}[b]{0.33\textwidth}
    		\includegraphics[width=\textwidth]{%
      		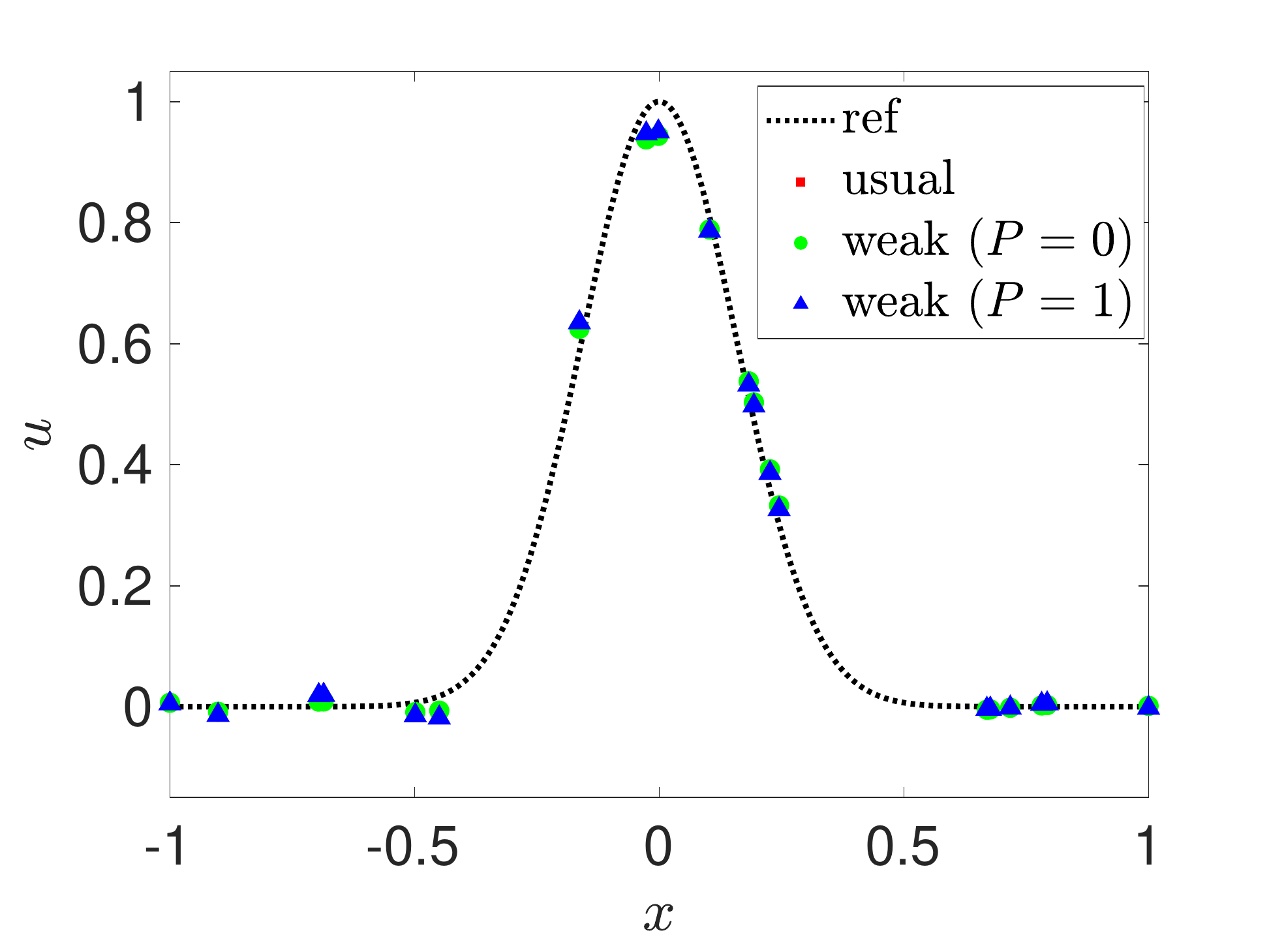}
    		\caption{Solution, random}
    		\label{fig:linear_1d_C_sol_random_quintic} 
  	\end{subfigure}%
  	~ 
  	\begin{subfigure}[b]{0.33\textwidth}
    		\includegraphics[width=\textwidth]{%
      		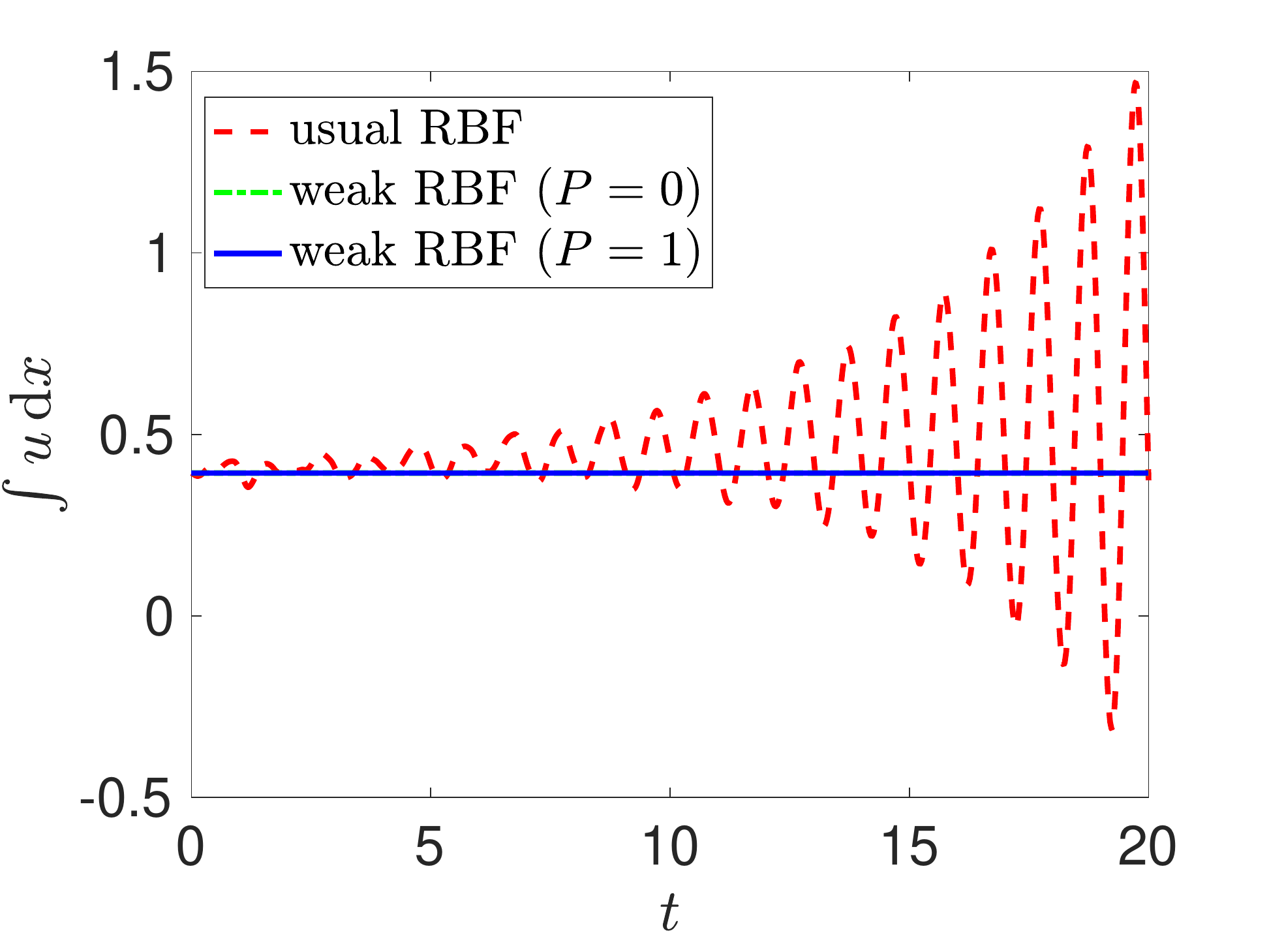}
    		\caption{Momentum, random}
    		\label{fig:linear_1d_C_momentum_random_quintic} 
  	\end{subfigure}%
	~ 
  	\begin{subfigure}[b]{0.33\textwidth}
    		\includegraphics[width=\textwidth]{%
      		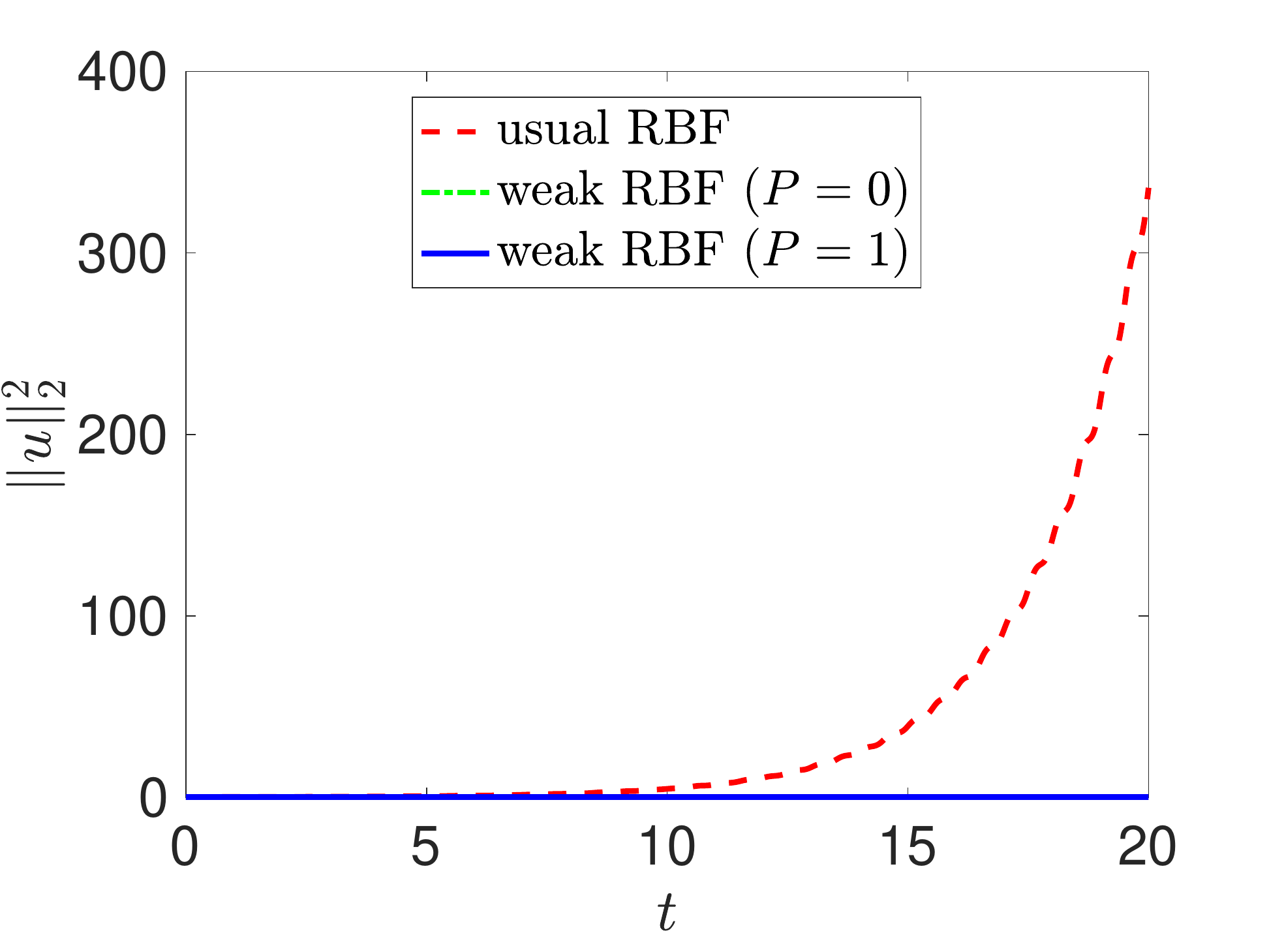}
    		\caption{Energy, random}
    		\label{fig:linear_1d_C_energy_random_quintic} 
  	\end{subfigure}%
  	\caption{
  		Numerical solutions at $t=20$ (left); momentum (middle); and energy (right) over time for $u_t + u_x = 0$ with periodic BCs. 
  		$N = 20$ equidistant (top) and randomly uniformly distributed (bottom) nodes are compared for a quintic kernel.
  	}
  	\label{fig:linear_1d_C_quintic}
\end{figure} 

Figures \ref{fig:linear_1d_C_cubic} and \ref{fig:linear_1d_C_quintic} illustrate this potential decrease in accuracy but preserved conservation and stability properties for the weak RBF method. 
The results for $N=20$ equidistant and randomly (uniformly distributed) nodes are compared for the cubic and quintic kernel. 
In all cases, the linear IVP \eqref{eq:linear_1d} with periodic BC is considered at time $t=20$. 
While the weak RBF method yields consistent results for all cases, the standard RBF method varies considerably, and essentially blew up when the quintic kernel was employed before the final time was reached (see the energy profile).

\subsubsection{Exact vs Numerical Integration} 
\label{subsub:num_int}

As discussed in \S \ref{sub:num-int}, to increase efficiency and reduce runtimes, an exact integration is often replaced by a numerical approximation.\footnote{In our implementation we are using  the  MATLAB function \emph{integral} for their computation so that strictly speaking, none of our integration is exact. This MATLAB function uses global adaptive quadrature and certain (default) error tolerances.}

\begin{figure}[!htb]
  	\centering
 	\begin{subfigure}[b]{0.33\textwidth}
    		\includegraphics[width=\textwidth]{%
      		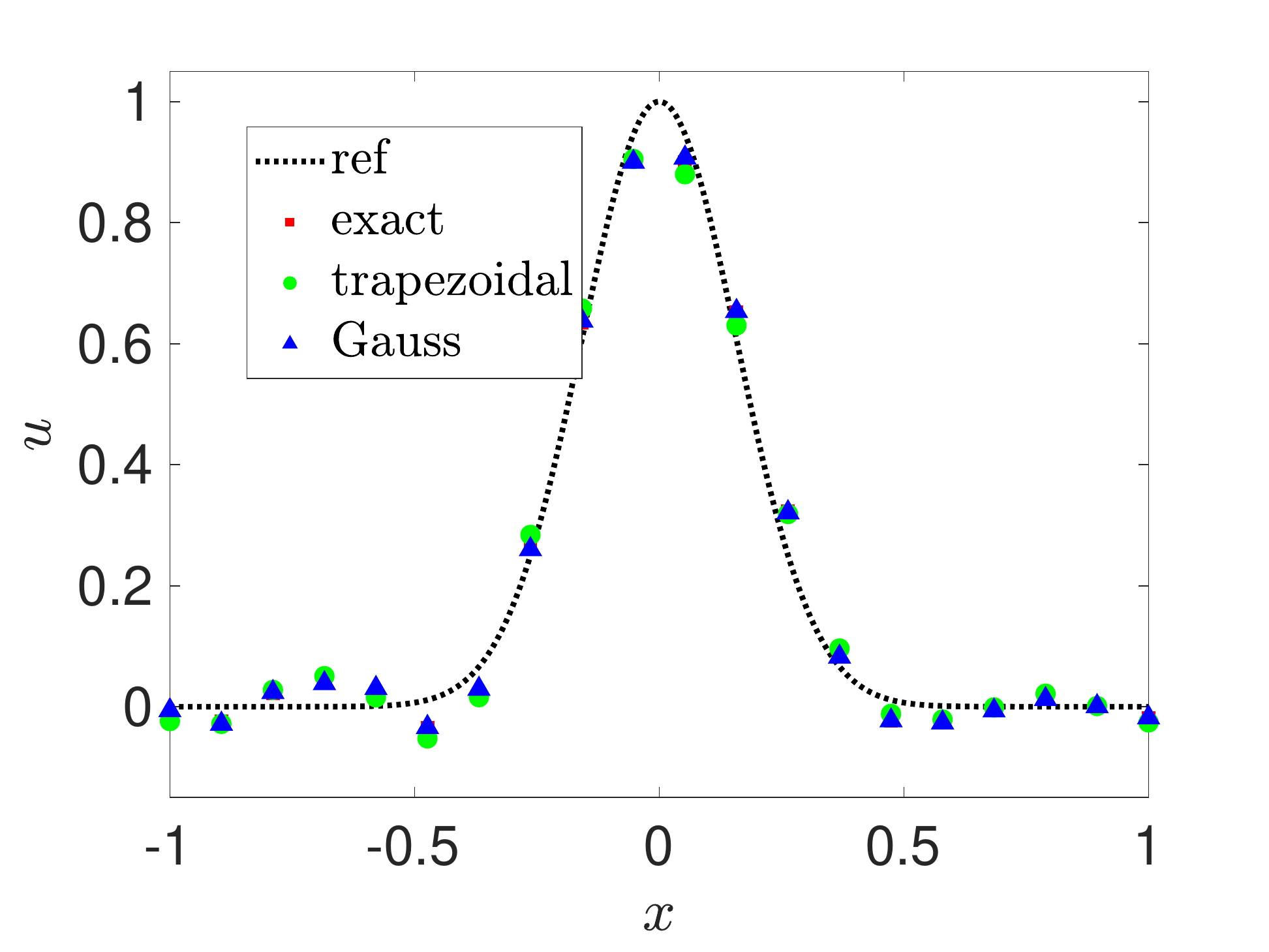}
    		\caption{Solution, cubic}
    		\label{fig:linear_1d_D_sol_cubic} 
  	\end{subfigure}%
  	~ 
  	\begin{subfigure}[b]{0.33\textwidth}
    		\includegraphics[width=\textwidth]{%
      		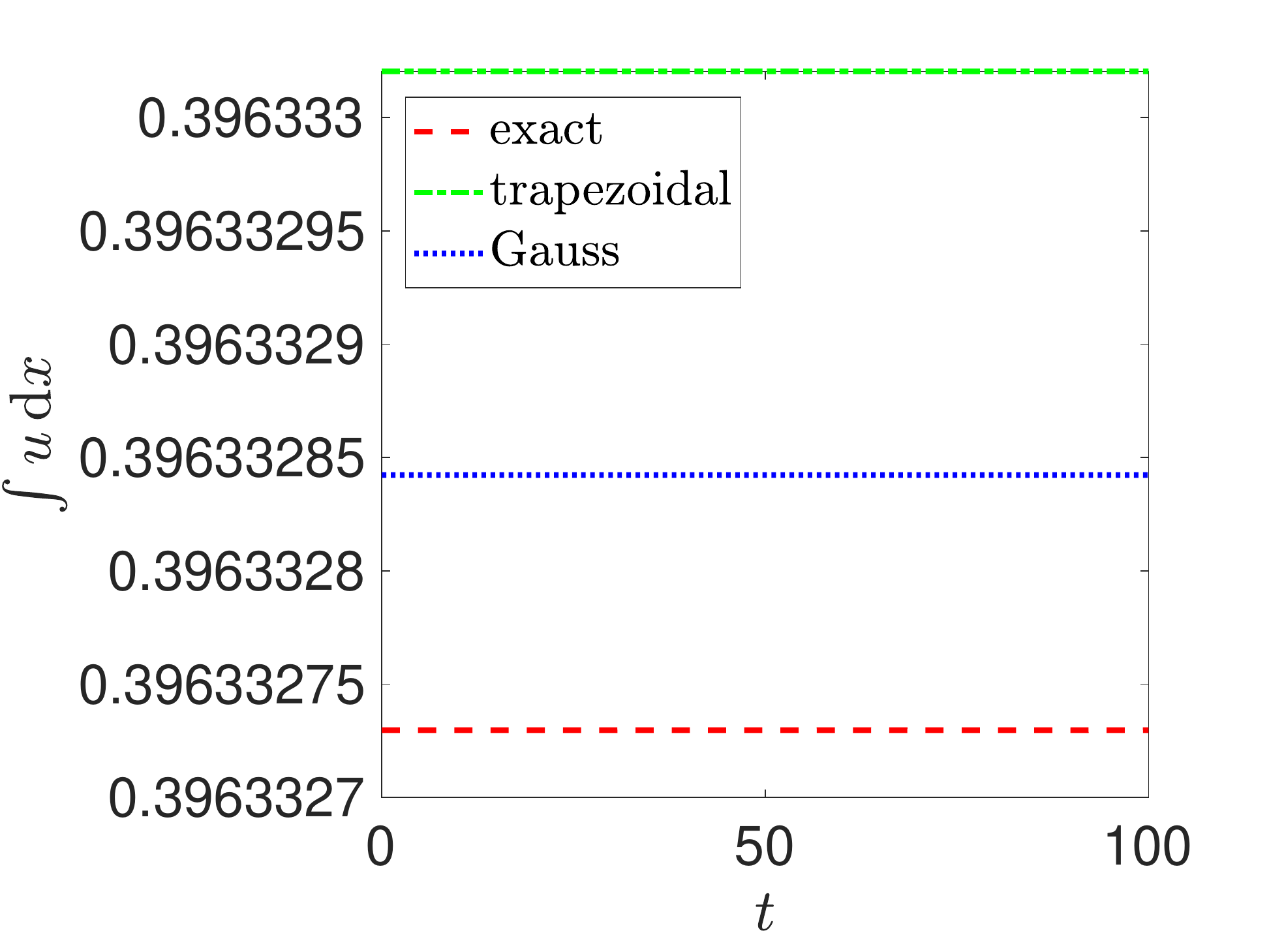}
    		\caption{Momentum, cubic}
    		\label{fig:linear_1d_D_momentum_cubic} 
  	\end{subfigure}%
	~ 
  	\begin{subfigure}[b]{0.33\textwidth}
    		\includegraphics[width=\textwidth]{%
      		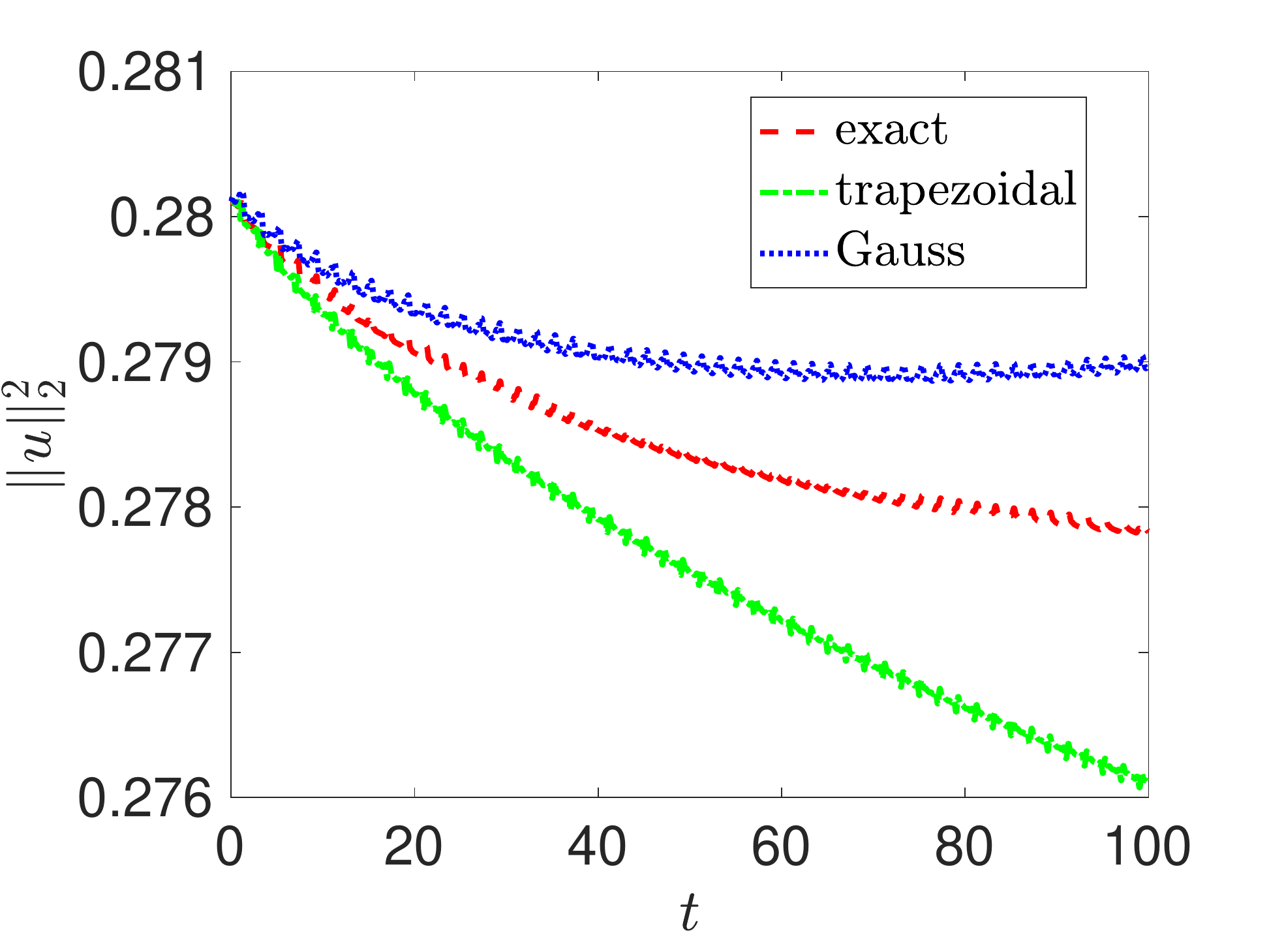}
    		\caption{Energy, cubic}
    		\label{fig:linear_1d_D_energy_cubic} 
  	\end{subfigure}%
	\\ 
  	\begin{subfigure}[b]{0.33\textwidth}
    		\includegraphics[width=\textwidth]{%
      		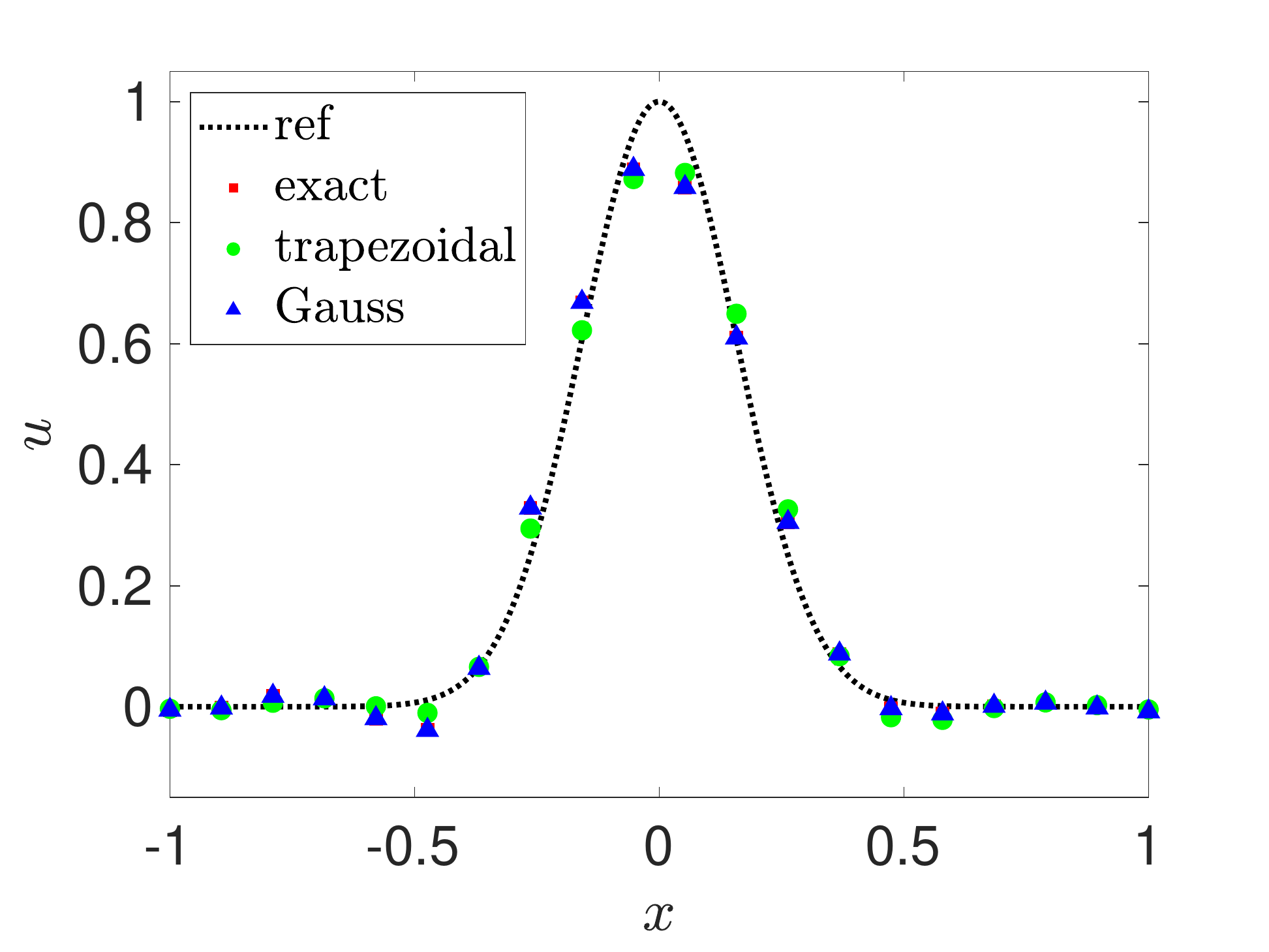}
    		\caption{Solution, quintic}
    		\label{fig:linear_1d_D_sol_quintic} 
  	\end{subfigure}%
  	~ 
  	\begin{subfigure}[b]{0.33\textwidth}
    		\includegraphics[width=\textwidth]{%
      		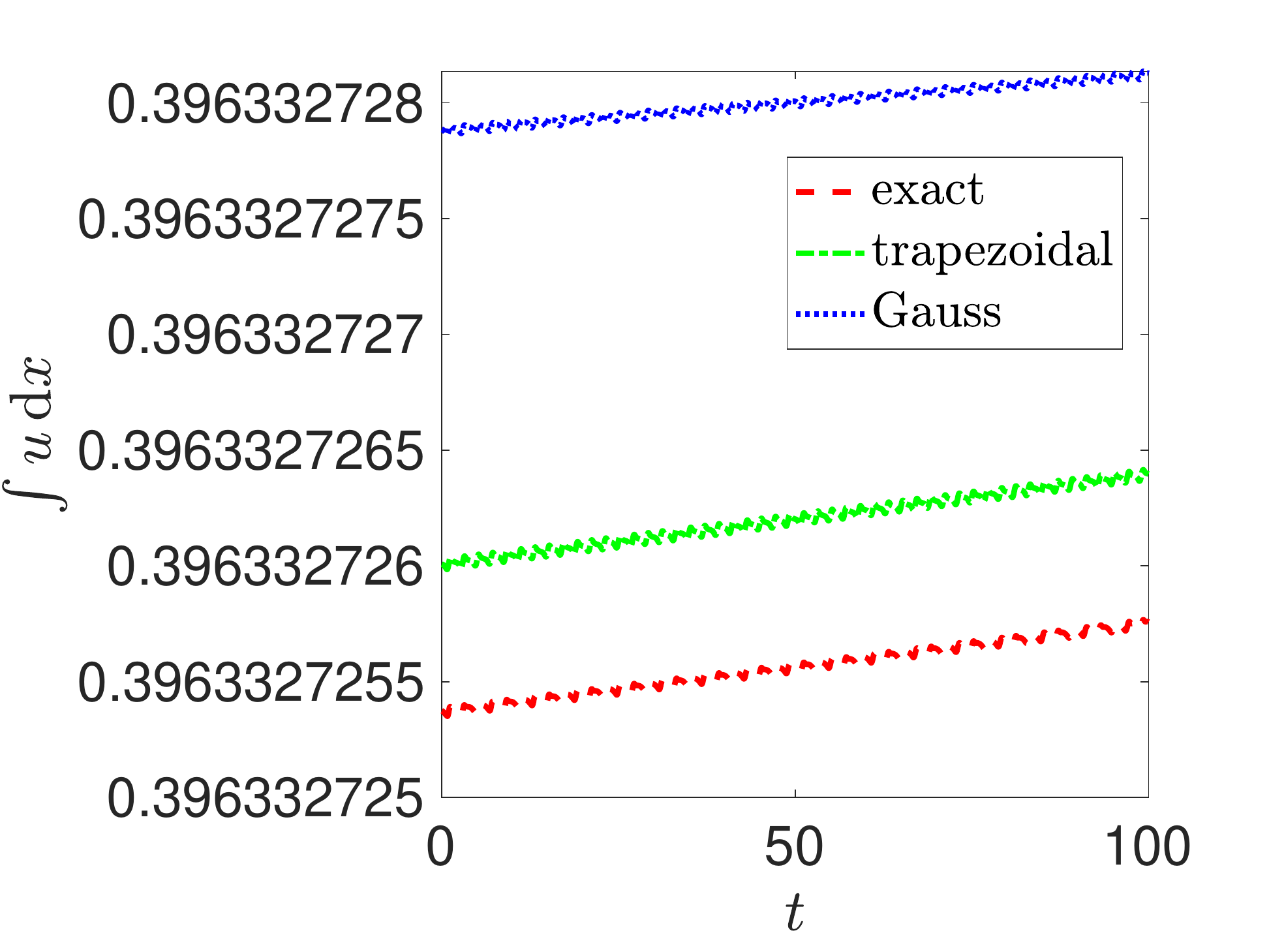}
    		\caption{Momentum, quintic}
    		\label{fig:linear_1d_D_momentum_random_quintic} 
  	\end{subfigure}%
	~ 
  	\begin{subfigure}[b]{0.33\textwidth}
    		\includegraphics[width=\textwidth]{%
      		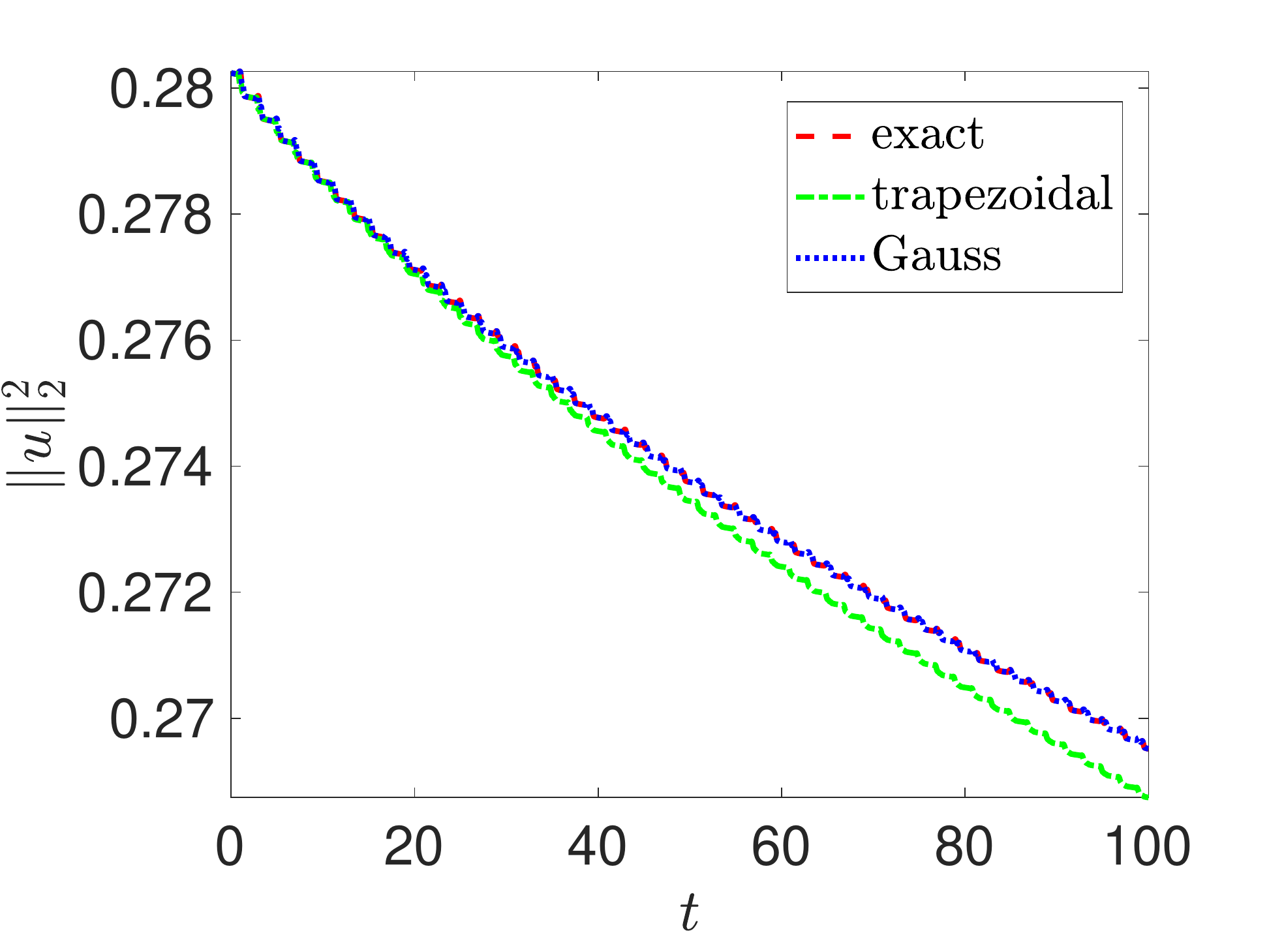}
    		\caption{Energy, quintic}
    		\label{fig:linear_1d_D_energy_quintic} 
  	\end{subfigure}%
  	\caption{
  		Numerical solutions at $t=100$ (left); their momentum (middle); and energy (right) over time for $u_t + u_x = 0$ with periodic BCs. 
  		Different integration techniques, all using $J=100$ quadrature points, are compared for a cubic and quintic kernel as well as $N = 20$ equidistant nodes.
  	}
  	\label{fig:linear_1d_D}
\end{figure} 

Figure \ref{fig:linear_1d_D} illustrates the numerical solution by the weak RBF method ($P=1$) together with the corresponding momentum and energy over time for the linear IVP \eqref{eq:linear_1d} with periodic BC and end time $t=100$. 
Here we compare ``exact'' integration (employing the MATLAB function \emph{integral}) with simple trapezoidal and Gauss(--Legendre) quadratures.  
The results demonstrate that even when only $J=100$ quadrature points are used, the numerical solution as well as the momentum are essentially the same for all integration techniques. This is also true for the energy 
in case of the quintic kernel. 
For the cubic kernel,  there are noticeable  differences in the energy for the different integration techniques, however. 
It is possible to overcome such discrepancies by increasing the number of quadrature points $J$ so that it is proportional to the number of nodes $N$. 
It is interesting to note that the trapezoidal rule yields more dissipation (lower energy profiles) than both the Gauss rule and ``exact'' integration. While the reasons for this should be investigated further, for now we simply note that the trapezoidal rule allows an efficient implementation of the weak RBF method while still preserving energy stability.

\subsection{Euler Equations} 
\label{sub:Euler} 

We now address the extension of the weak RBF method to systems of nonlinear hyperbolic CLs. 
To this end, we consider the one-dimensional Euler equations given by
\begin{equation}
  	\boldsymbol{U}_t + \boldsymbol{F}(\boldsymbol{U})_x = 0
\end{equation}
for ${x \in \Omega = [-1,1]}$, where $\boldsymbol{U}$ and $\boldsymbol{F}(\boldsymbol{U})$ respectively are the vector of conserved variables and fluxes:  
\begin{equation}
  \boldsymbol{U} = 
  \begin{pmatrix}
    u_1 \\ u_2 \\ u_3 
  \end{pmatrix}
  = 
  \begin{pmatrix} 
    \rho \\ \rho u \\ E
  \end{pmatrix}, \quad 
  \boldsymbol{F} = 
  \begin{pmatrix}
    f_1 \\ f_2 \\ f_3 
  \end{pmatrix}
  = 
  \begin{pmatrix} 
    \rho u \\ \rho u^2 + p \\ u(E+p)
  \end{pmatrix}. 
\end{equation}
Here, $\rho$ is the {density}, $u$ is the {velocity}, $p$ is the {pressure}, and $E$ is the {total energy per unit volume}. 
The Euler equations are completed by addition of an equation of state (EOS) with general form 
\begin{equation}
  p = p(\rho,e),
\end{equation} 
where $e = E/\rho - u^2/2$ is the {specific internal energy}.  
For the case of ideal gases the EOS is given by
\begin{equation}
  p = (\gamma - 1) \rho e
\end{equation}
with $\gamma$ denoting the {ratio of specific heats}. 
For the subsequent numerical tests, we set $\gamma = 3$ and consider a smooth isentropic flow resulting from the Euler equations with smooth ICs 
\begin{equation}
\begin{aligned}
  \rho(0,x) = 1 +\frac{1}{2} \sin(\pi x), \quad 
  u(0,x) = 0, \quad 
  p(0,x) = \rho(0,x)^{\gamma},
\end{aligned} 
\end{equation}
and periodic BCs. 
A similar test problem has been proposed in \cite{cheng2014positivity} as well as  
in \cite{abgrall2019high} in the context of (positivity-preserving) high-order methods.
Utilizing the method of characteristics, the exact density $\rho$ and velocity $u$ are given by 
\begin{equation}
\begin{aligned}
  \rho(t,x) = \frac{1}{2} \left[ \rho_0(x_1) + \rho_0(x_2) \right], \quad 
  u(t,x) = \sqrt{3} \left[ \rho(t,x) - \rho_0(x_1) \right],
\end{aligned}
\end{equation}
where $x_1 = x_1(t,x)$ and $x_2 = x_2(t,x)$ are solutions of the nonlinear equations 
\begin{equation}
  x + \sqrt{3} \rho_0(x_1) t - x_1 = 0, \quad 
  x - \sqrt{3} \rho_0(x_2) t - x_2 = 0.
\end{equation} 
Finally, the exact pressure $p$ can be computed by the isentropic law $p = C \rho^{\gamma}$ for smooth flows \cite[Chapter 3.1]{toro2013riemann}.

\begin{figure}[!htb]
  	\centering
 	\begin{subfigure}[b]{0.33\textwidth}
    		\includegraphics[width=\textwidth]{%
      		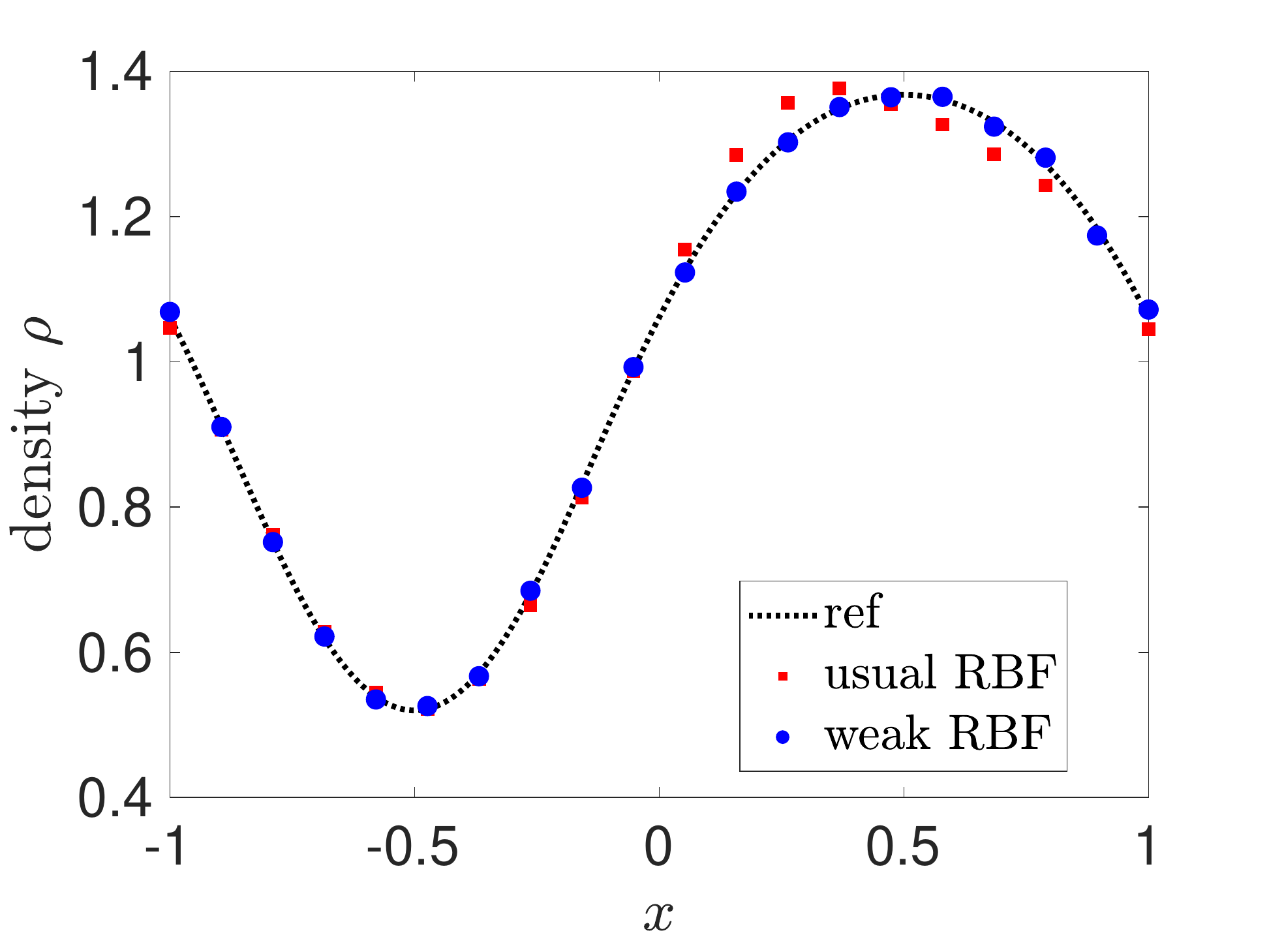}
    		\caption{Density, cubic}
    		\label{fig:Euler_1d_density_cubic} 
  	\end{subfigure}%
  	~ 
  	\begin{subfigure}[b]{0.33\textwidth}
    		\includegraphics[width=\textwidth]{%
      		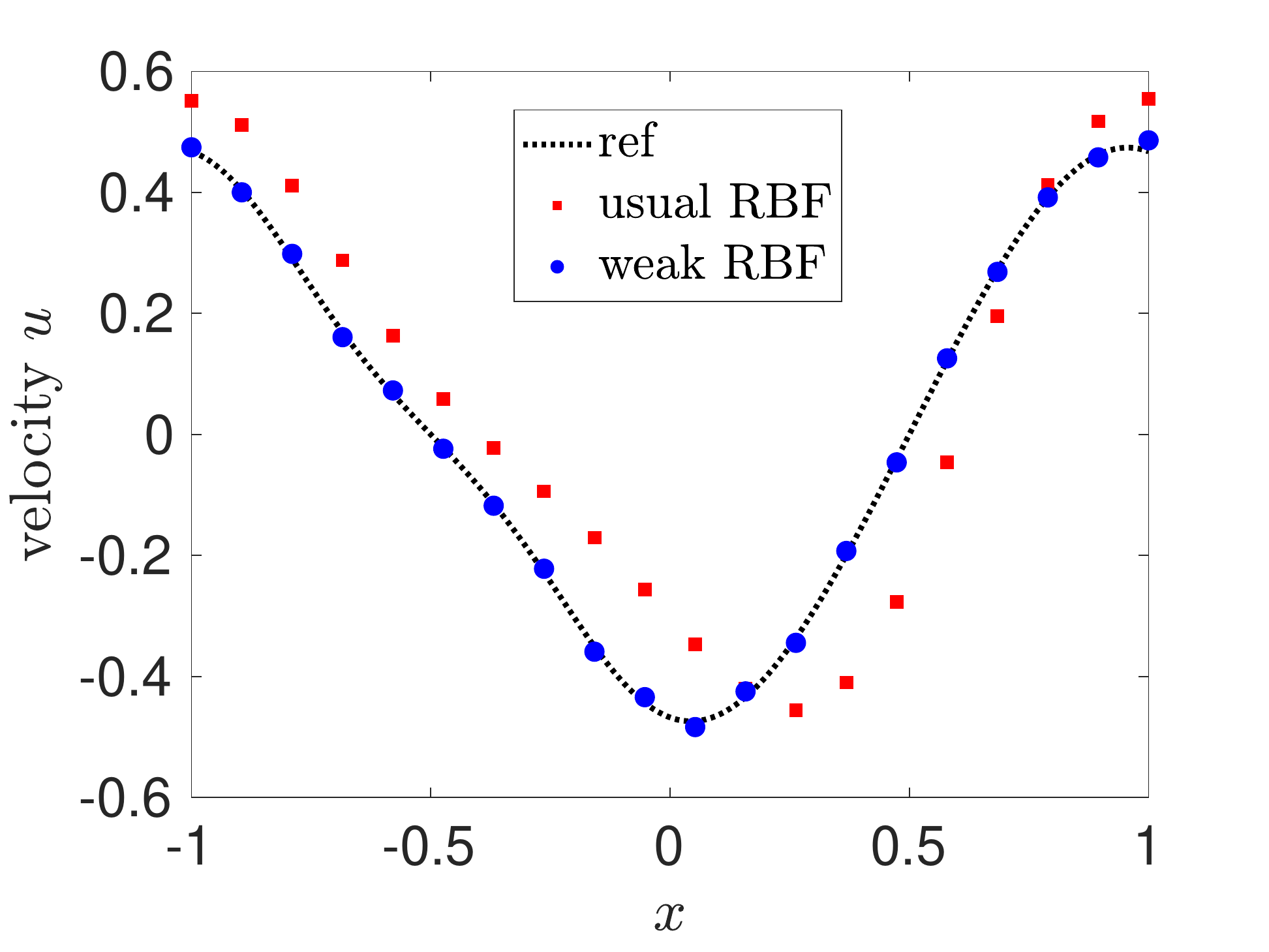}
    		\caption{Velocity, cubic}
    		\label{fig:Euler_1d_velocity_cubic} 
  	\end{subfigure}%
	~ 
  	\begin{subfigure}[b]{0.33\textwidth}
    		\includegraphics[width=\textwidth]{%
      		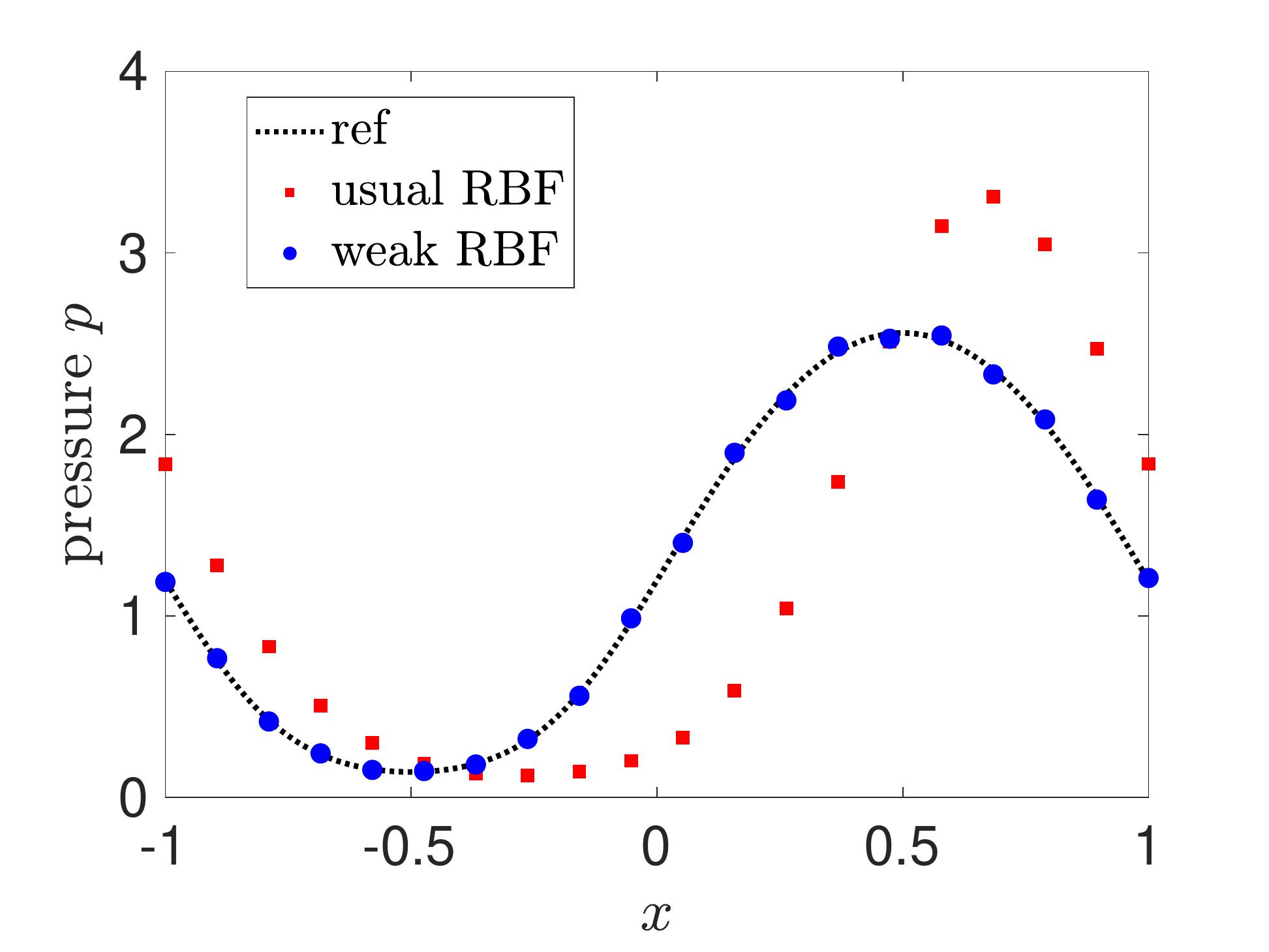}
    		\caption{Pressure, cubic}
    		\label{fig:Euler_1d_pressure_cubic} 
  	\end{subfigure}%
	\\ 
  	\begin{subfigure}[b]{0.33\textwidth}
    		\includegraphics[width=\textwidth]{%
      		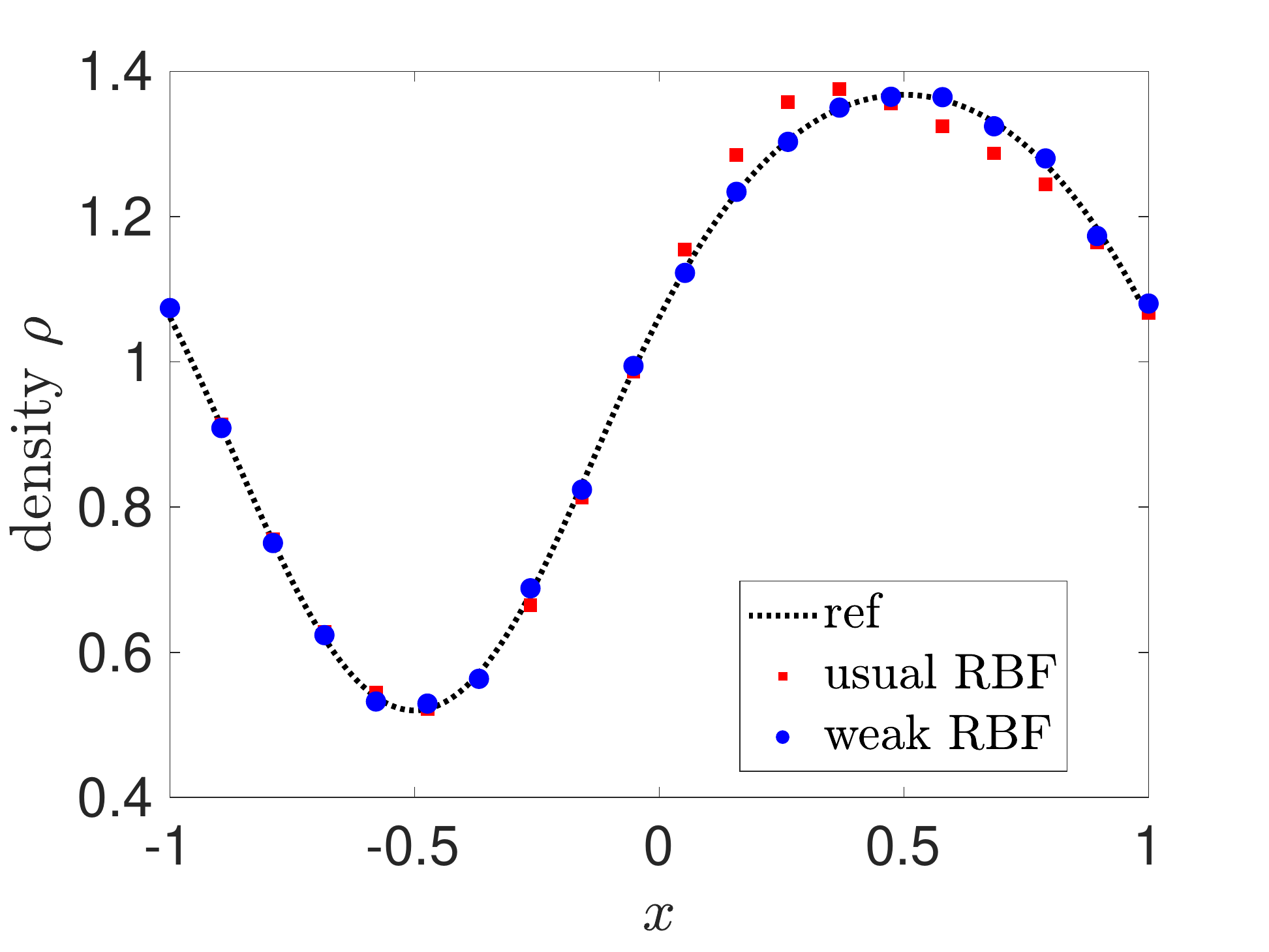}
    		\caption{Density, quintic}
    		\label{fig:Euler_1d_density_quintic} 
  	\end{subfigure}%
  	~ 
  	\begin{subfigure}[b]{0.33\textwidth}
    		\includegraphics[width=\textwidth]{%
      		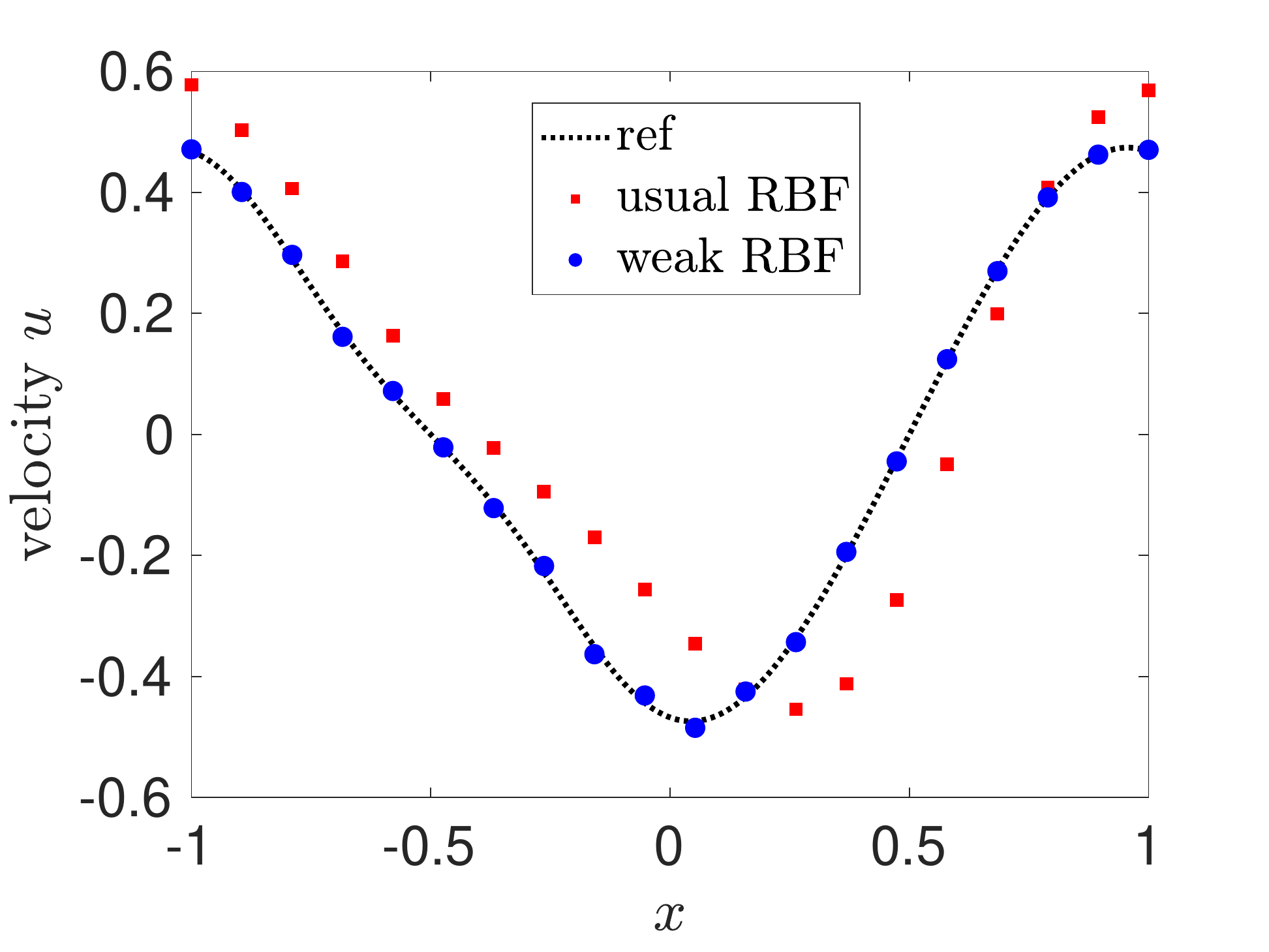}
    		\caption{Velocity, quintic}
    		\label{fig:Euler_1d_velocity_quintic} 
  	\end{subfigure}%
	~ 
  	\begin{subfigure}[b]{0.33\textwidth}
    		\includegraphics[width=\textwidth]{%
      		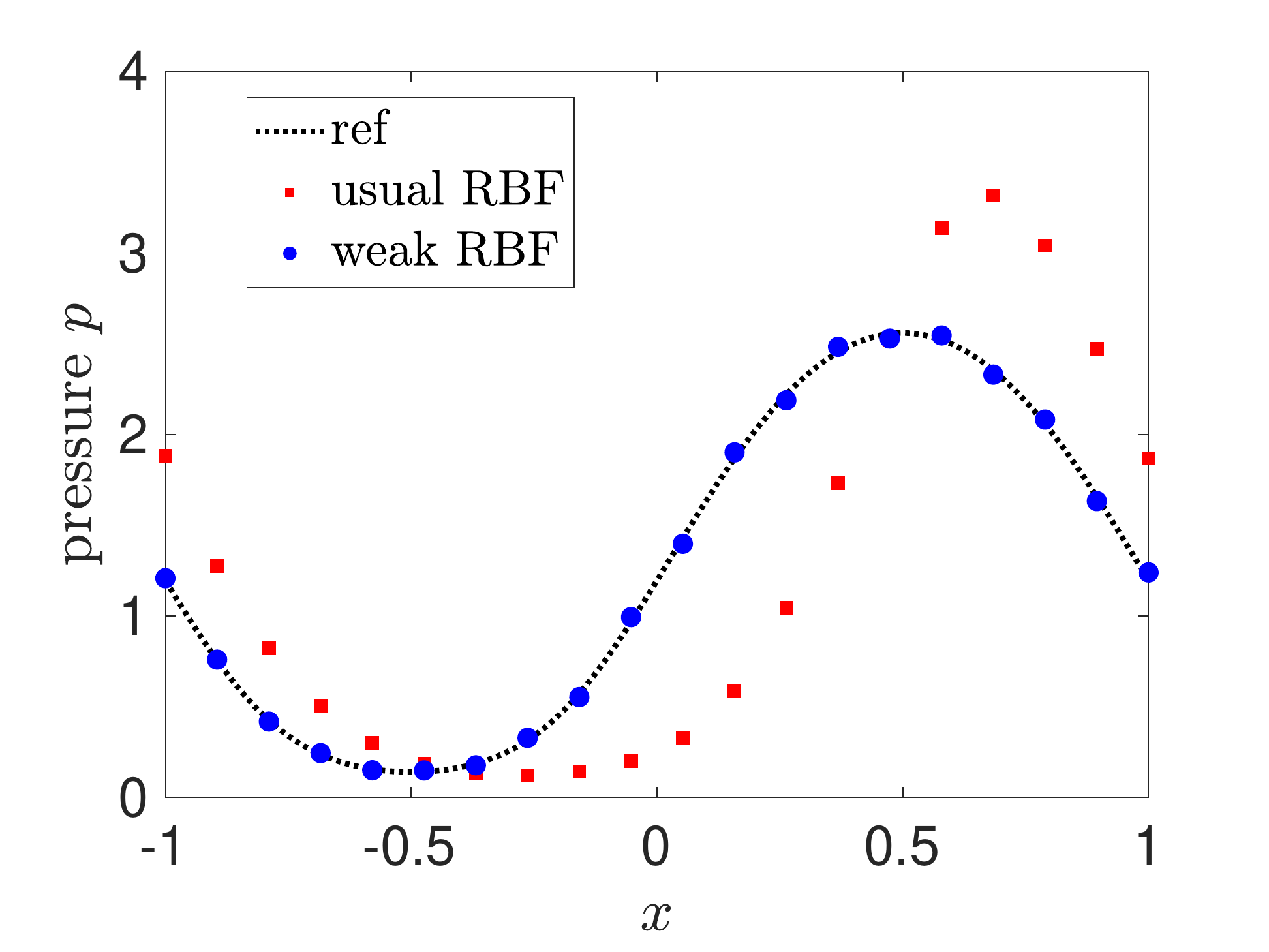}
    		\caption{Pressure, quintic}
    		\label{fig:Euler_1d_pressure_quintic} 
  	\end{subfigure}%
  	\caption{
  		Numerical results (density, velocity, and pressure at the final time $t=0.1$) for the Euler 
equations. 
		The cubic and quintic kernel with $N = 20$ equidistant nodes were used. 
		The weak RBF method includes constants ($P=1$).
  	}
  	\label{fig:Euler_1d}
\end{figure} 

Figure \ref{fig:Euler_1d} illustrates the numerical results at time $t=0.1$ comparing the strong and weak RBF collocation method using the cubic and quintic kernel. 
For the weak RBF method, constants have been included ($P=1$).
As in the case for linear advection, we observe that the weak RBF collocation method for is more accurate than the strong RBF.

\subsection{{Extension to Two Dimensions}}
\label{sub:2d-num} 

To conclude our numerical experiments we apply the weak RBF method to a two-dimensional problem and consider
\begin{equation}\label{eq:2d-problem}
\begin{aligned}
  	u_t + u_x  & = 0, \\ 
	u(0,x,y) & = \sin(2\pi x) \left( \frac{1}{2} \sin(2\pi y) - 1 \right), \\
\end{aligned}
\end{equation}
on ${\Omega=[-1,1]^2 \subset \R^2}$ with periodic BCs.
This test is designed to demonstrate the validity of conservation and energy stability of the weak RBF method in higher dimensions, as discussed in \S \ref{sub:2d-stability}.
In addition, it is once more illustrated that these properties are not affected by using a nonequidistant distribution of nodes, in this case random uniformly distributed. 
Finally, this example also illustrates the limitations of the proposed weak RBF methods for long time simulations. 

\begin{figure}[!htb]
  	\centering
 	\begin{subfigure}[b]{0.33\textwidth}
    		\includegraphics[width=\textwidth]{%
      		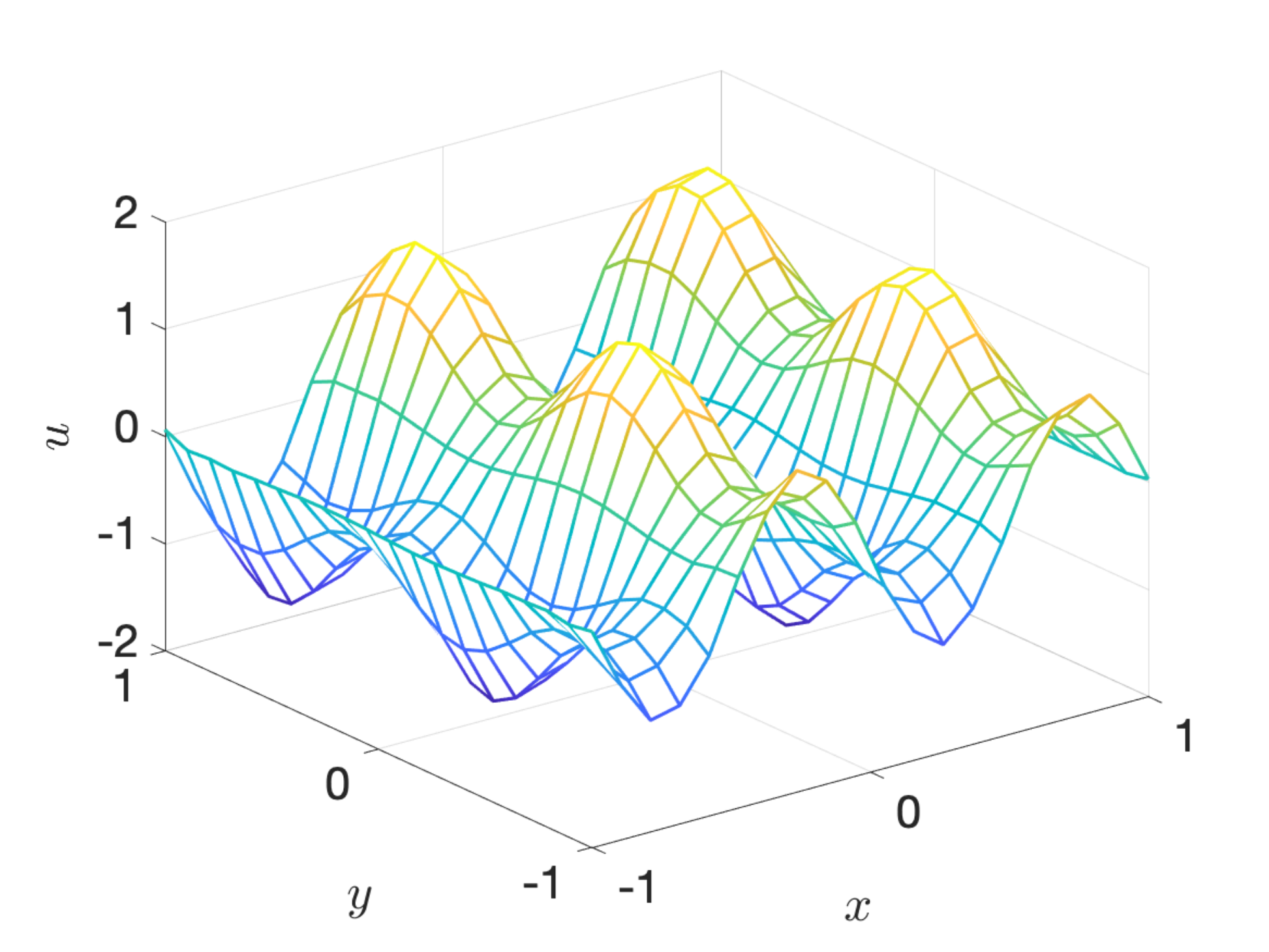}
    		\caption{Reference, $t=10$}
    		\label{fig:linear_2d_cubic_equid_ref_t=10} 
  	\end{subfigure}%
  	~ 
  	\begin{subfigure}[b]{0.33\textwidth}
    		\includegraphics[width=\textwidth]{%
      		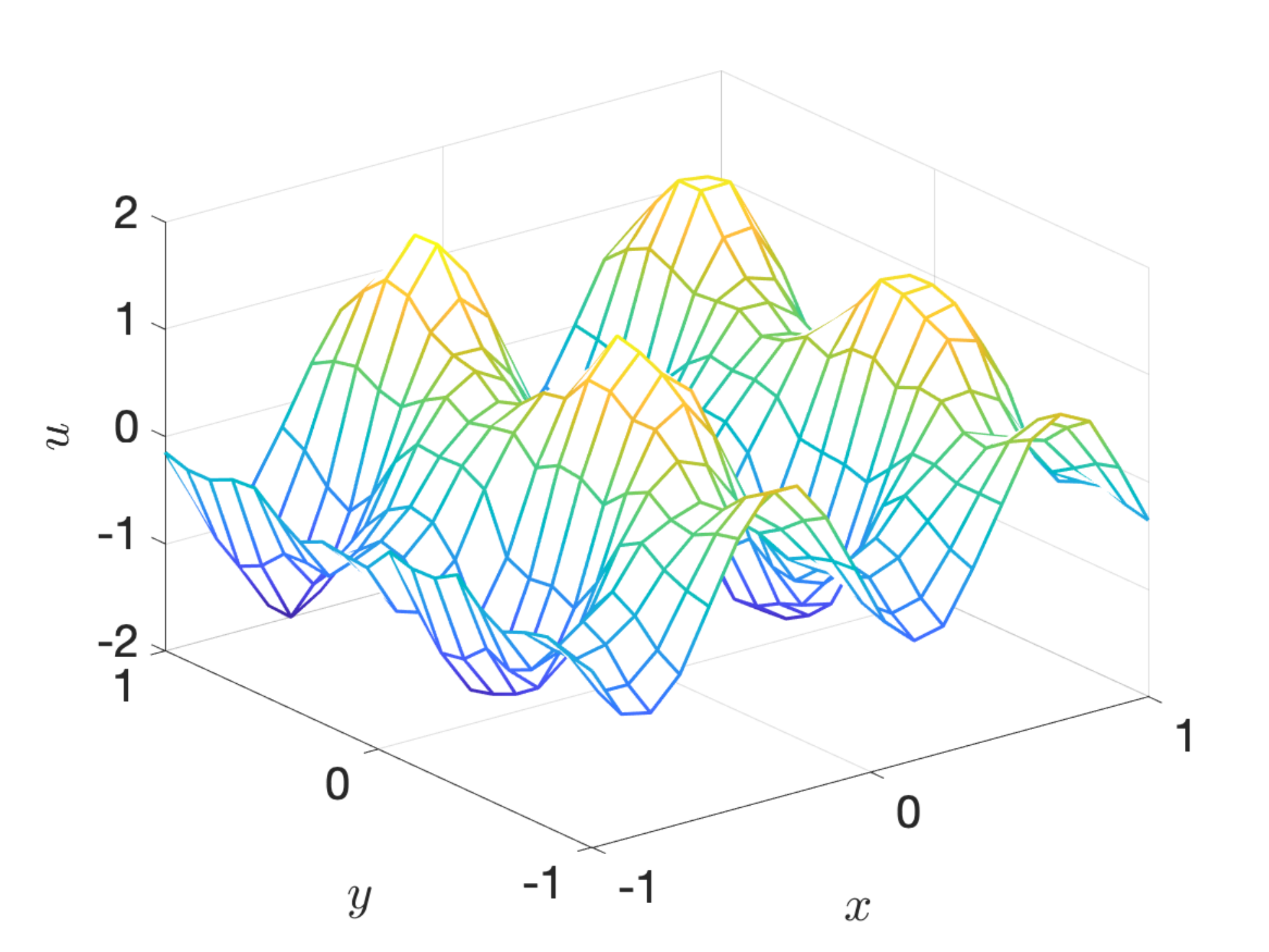}
    		\caption{{Standard} RBF, $t=10$}
    		\label{fig:linear_2d_cubic_equid_usual_t=10} 
  	\end{subfigure}%
	~ 
  	\begin{subfigure}[b]{0.33\textwidth}
    		\includegraphics[width=\textwidth]{%
      		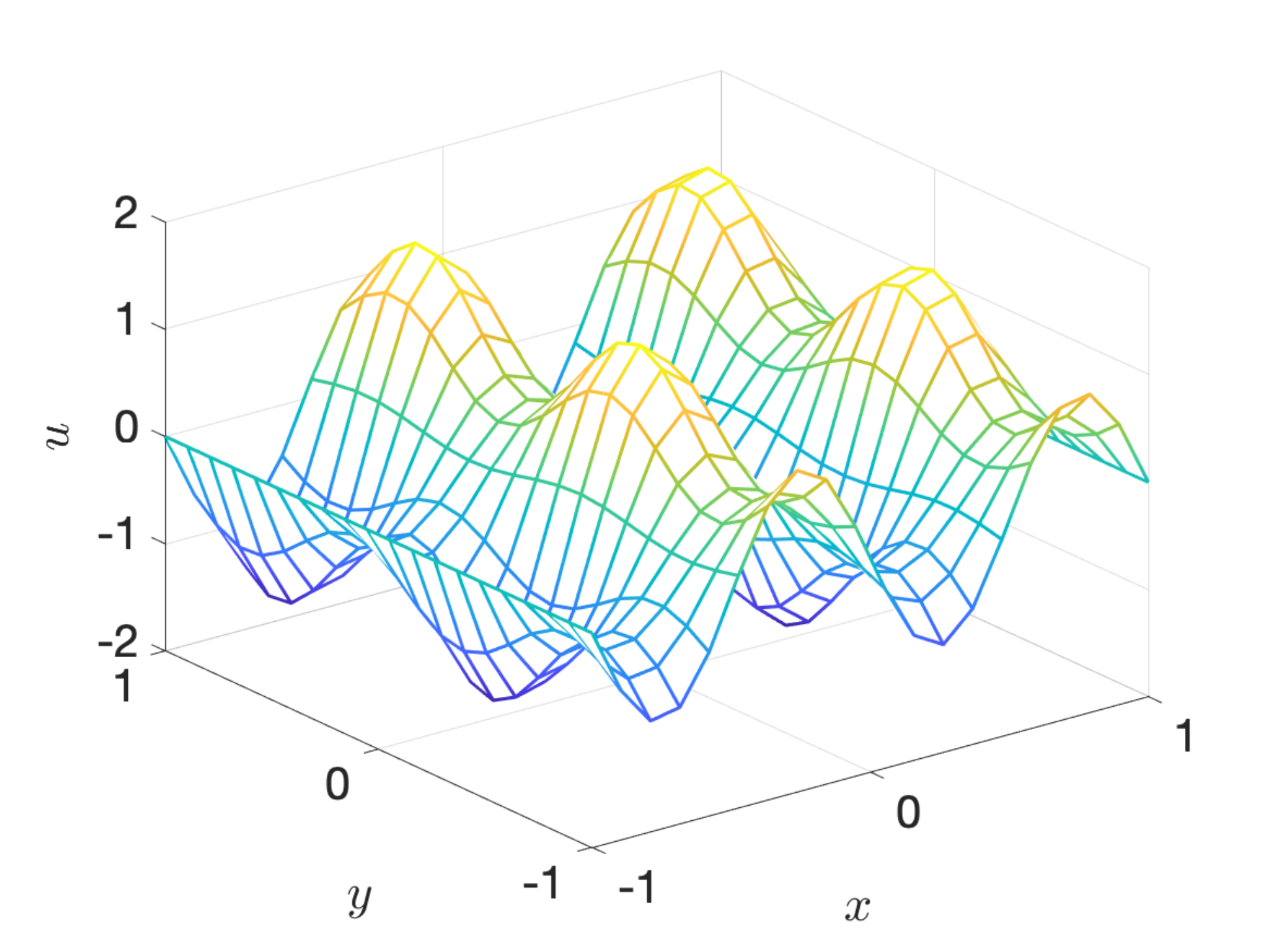}
    		\caption{Weak RBF, $t=10$}
    		\label{fig:linear_2d_cubic_equid_weak_t=10} 
  	\end{subfigure}%
	\\ 
  	\begin{subfigure}[b]{0.33\textwidth}
    		\includegraphics[width=\textwidth]{%
      		plots/linear_2d_ref}
    		\caption{Reference, $t=400$}
    		\label{fig:linear_2d_cubic_equid_ref_t=400} 
  	\end{subfigure}%
  	~ 
  	\begin{subfigure}[b]{0.33\textwidth}
    		\includegraphics[width=\textwidth]{%
      		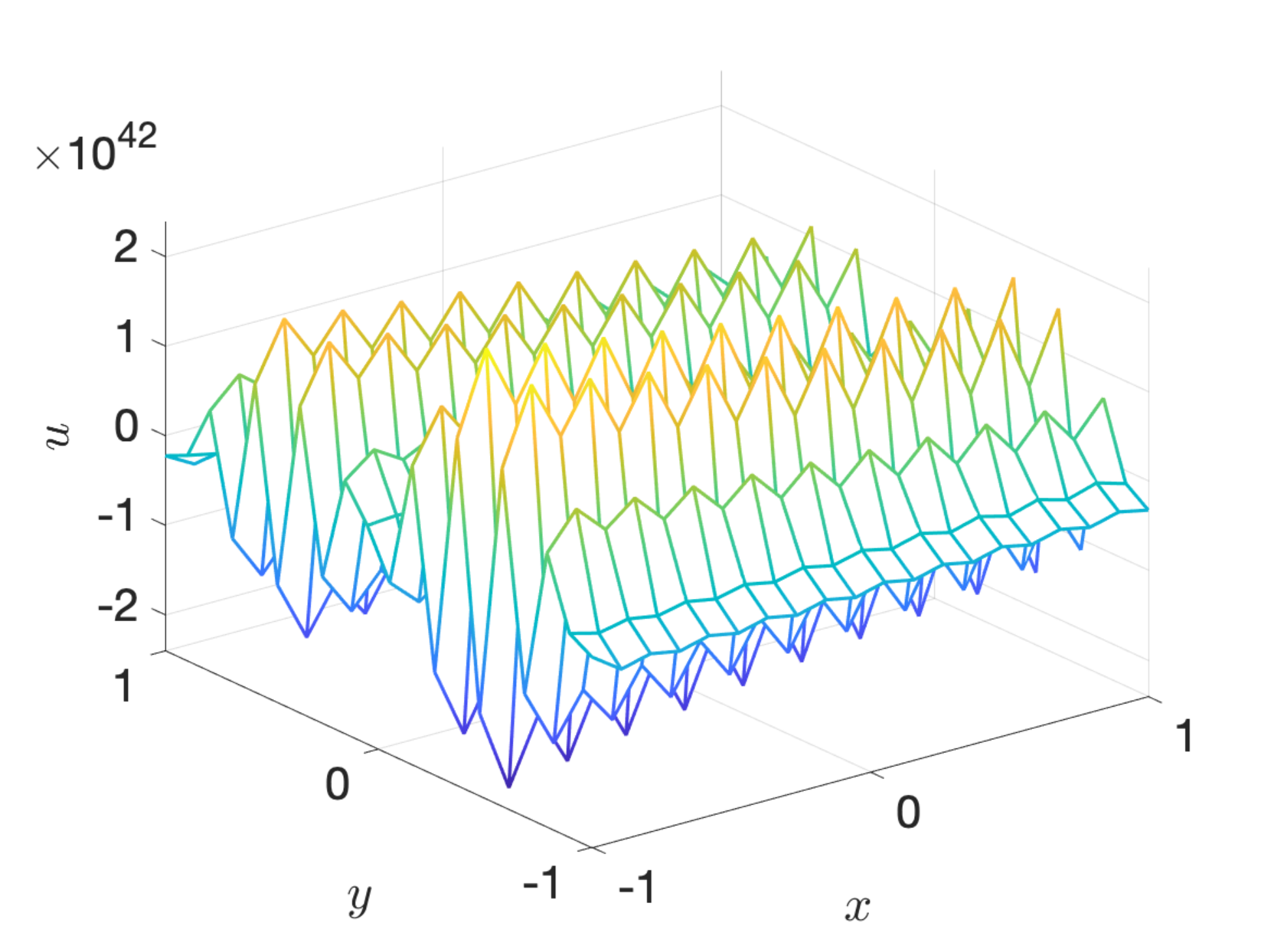}
    		\caption{{Standard} RBF, $t=400$}
    		\label{fig:linear_2d_cubic_equid_usual_t=400} 
  	\end{subfigure}%
	~ 
  	\begin{subfigure}[b]{0.33\textwidth}
    		\includegraphics[width=\textwidth]{%
      		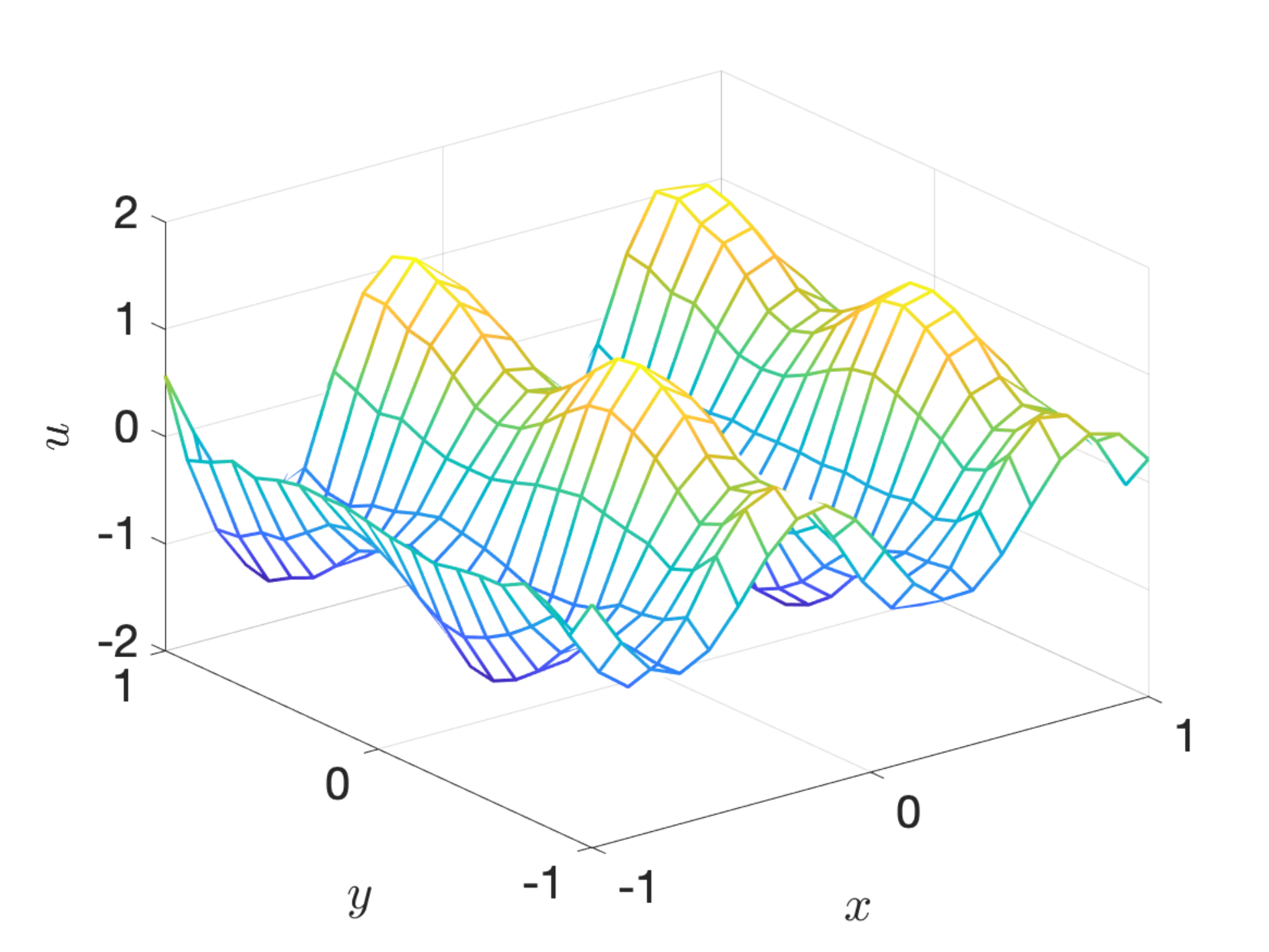}
    		\caption{Weak RBF, $t=400$}
    		\label{fig:linear_2d_cubic_equid_weak_t=400} 
  	\end{subfigure}%
  	\caption{
  		Numerical results for the two-dimensional linear IVP with periodic BCs. 
		The cubic kernel with $N = 400$ equidistant nodes was used.   	
	}
  	\label{fig:linear_2d_cubic_equid}
\end{figure} 

\begin{figure}[!htb]
  	\centering
 	\begin{subfigure}[b]{0.33\textwidth}
    		\includegraphics[width=\textwidth]{%
      		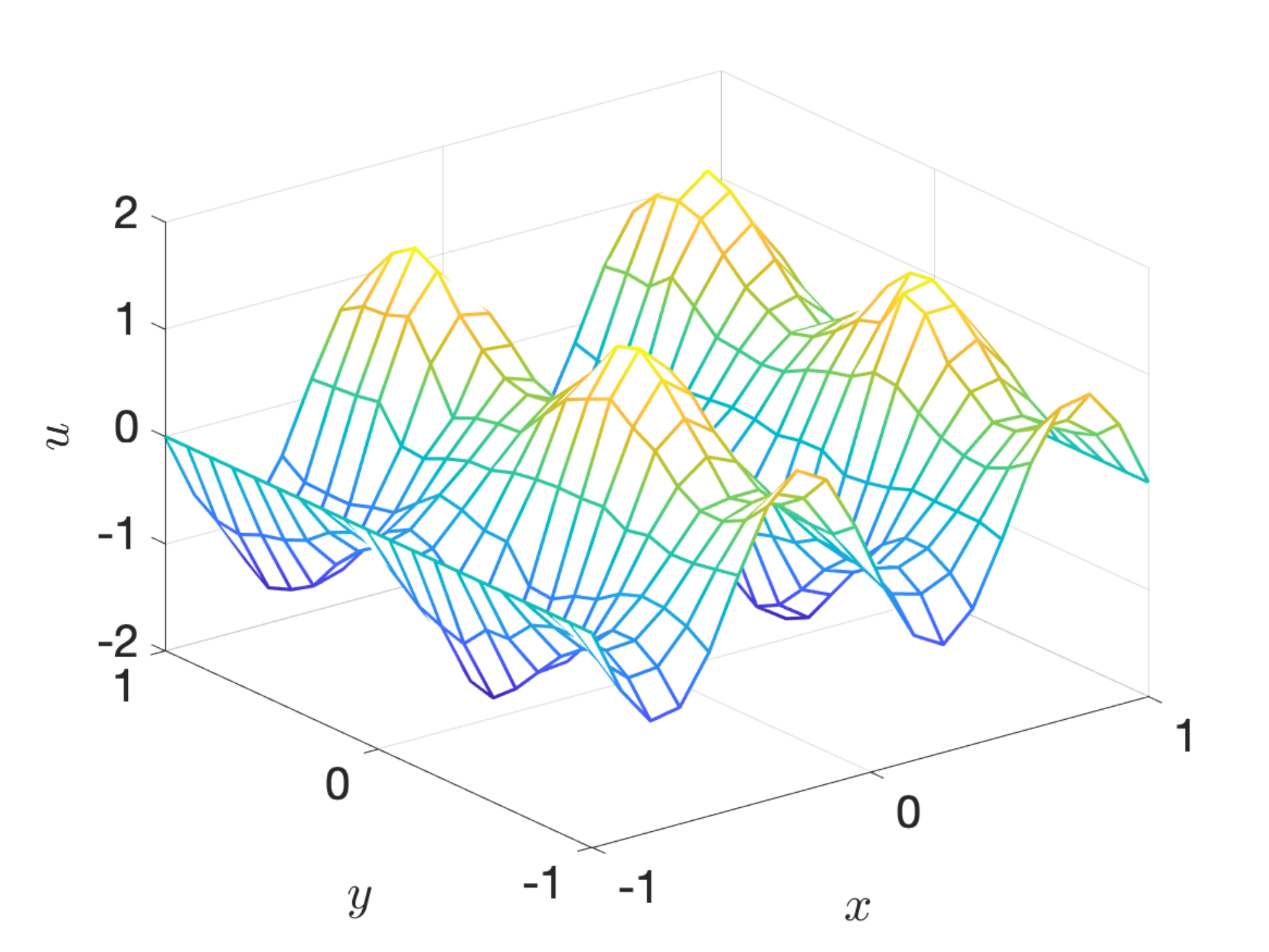}
    		\caption{Reference, $t=10$}
    		\label{fig:linear_2d_cubic_random_ref_t=10} 
  	\end{subfigure}%
  	~ 
  	\begin{subfigure}[b]{0.33\textwidth}
    		\includegraphics[width=\textwidth]{%
      		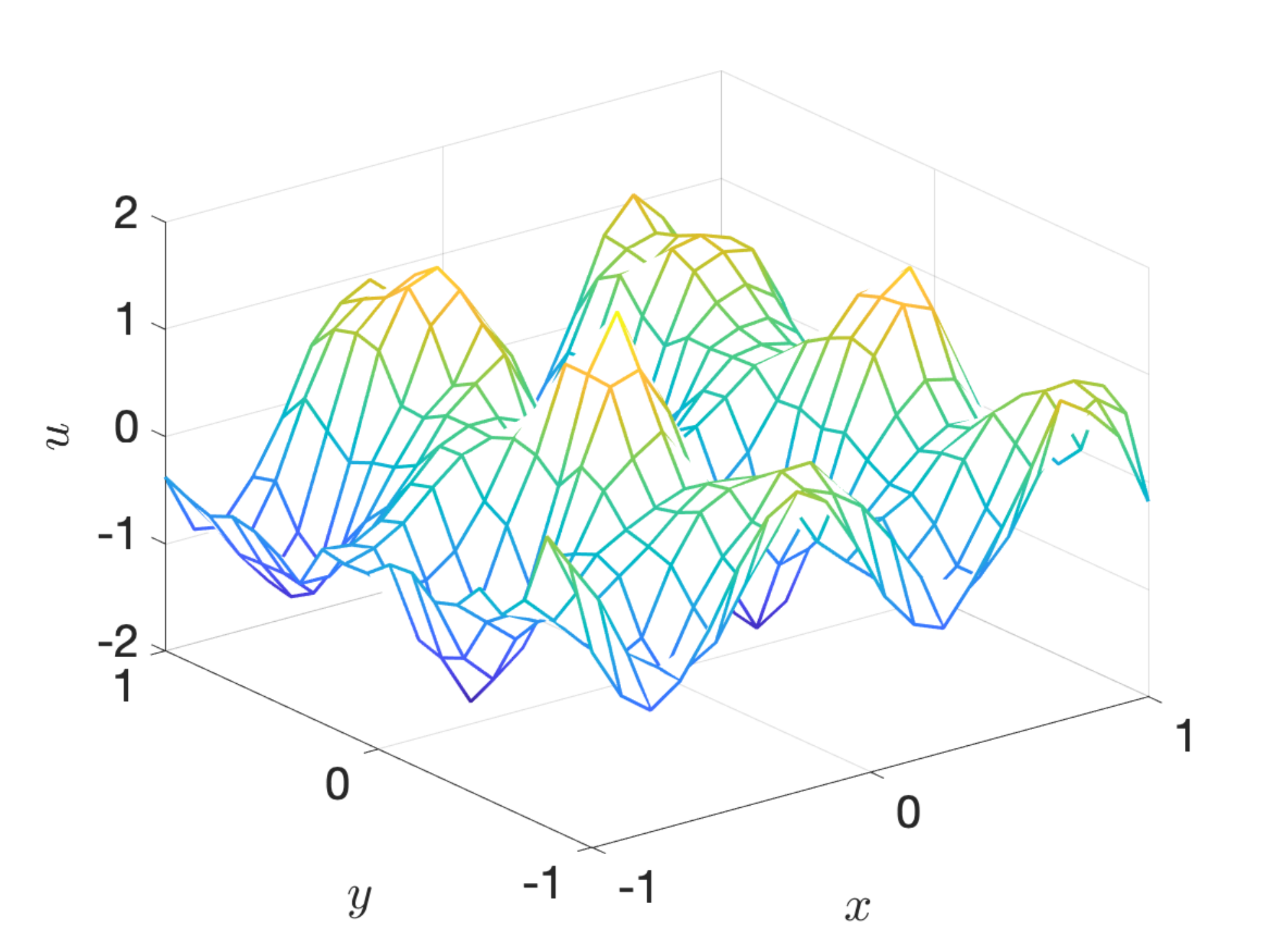}
    		\caption{{Standard} RBF, $t=10$}
    		\label{fig:linear_2d_cubic_random_usual_t=10} 
  	\end{subfigure}%
	~ 
  	\begin{subfigure}[b]{0.33\textwidth}
    		\includegraphics[width=\textwidth]{%
      		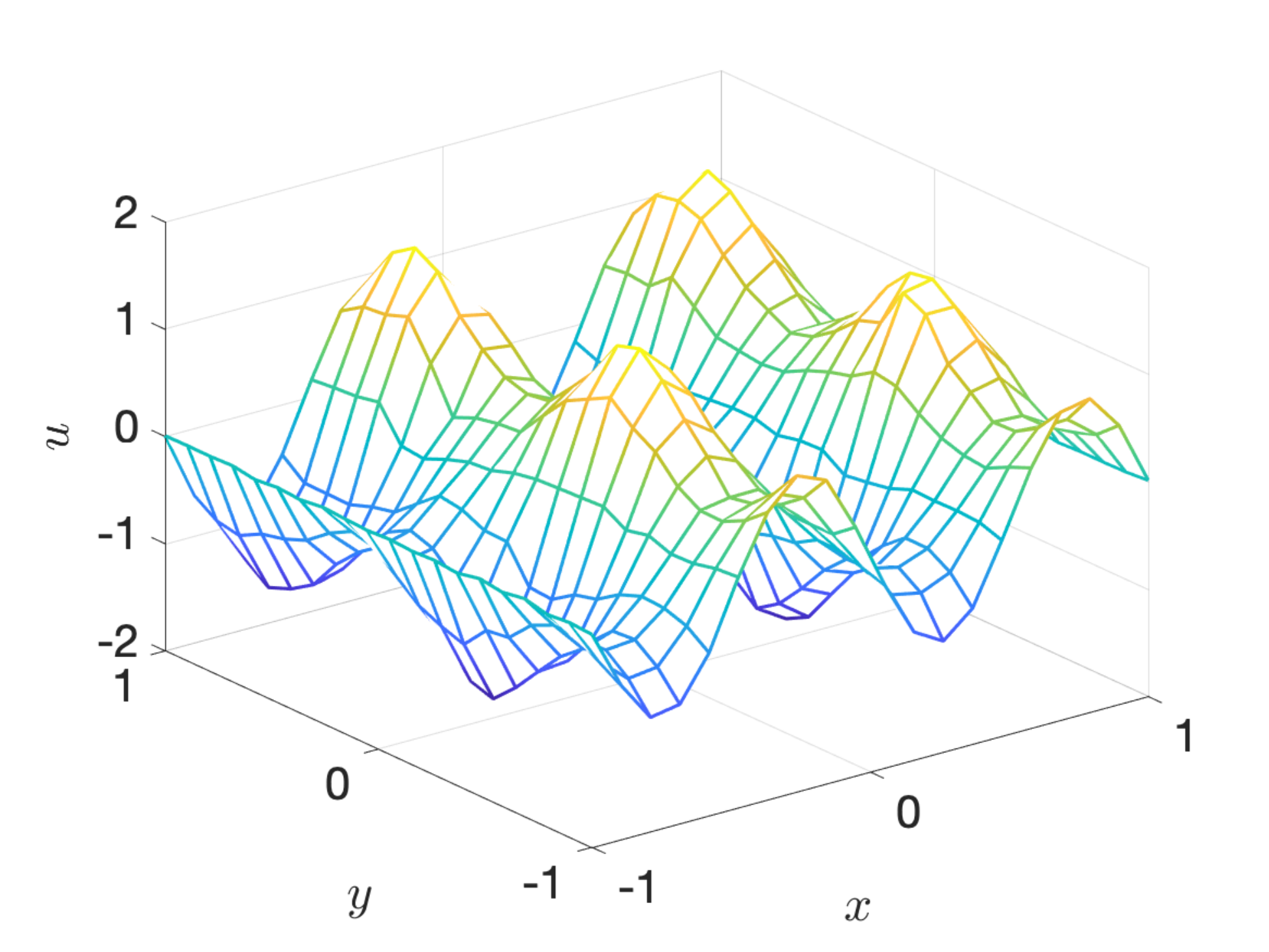}
    		\caption{Weak RBF, $t=10$}
    		\label{fig:linear_2d_cubic_random_weak_t=10} 
  	\end{subfigure}%
	\\ 
  	\begin{subfigure}[b]{0.33\textwidth}
    		\includegraphics[width=\textwidth]{%
      		plots/linear_2d_ref_random}
    		\caption{Reference, $t=400$}
    		\label{fig:linear_2d_cubic_random_ref_t=400} 
  	\end{subfigure}%
  	~ 
  	\begin{subfigure}[b]{0.33\textwidth}
    		\includegraphics[width=\textwidth]{%
      		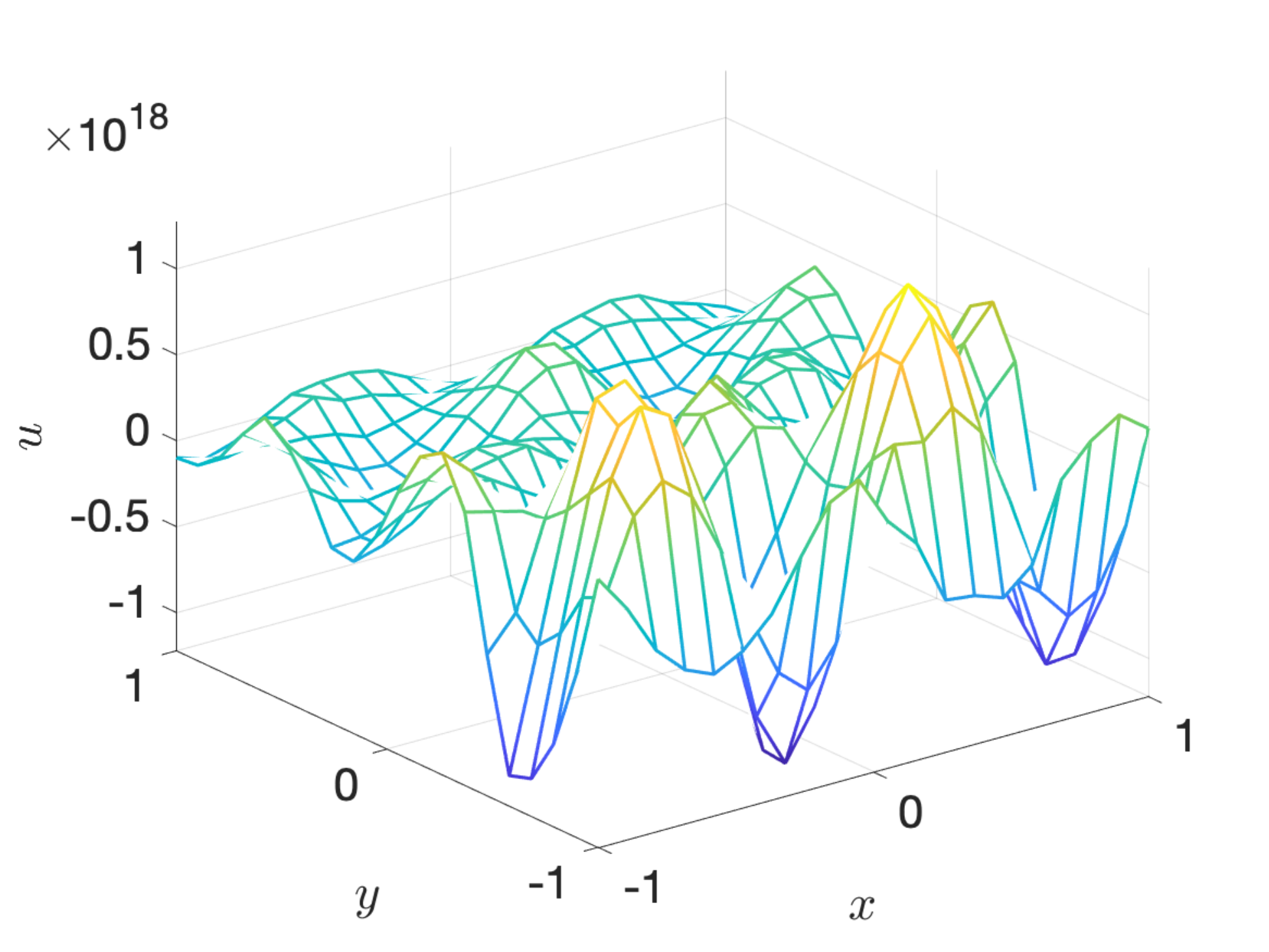}
    		\caption{{Standard} RBF, $t=400$}
    		\label{fig:linear_2d_cubic_random_usual_t=400} 
  	\end{subfigure}%
	~ 
  	\begin{subfigure}[b]{0.33\textwidth}
    		\includegraphics[width=\textwidth]{%
      		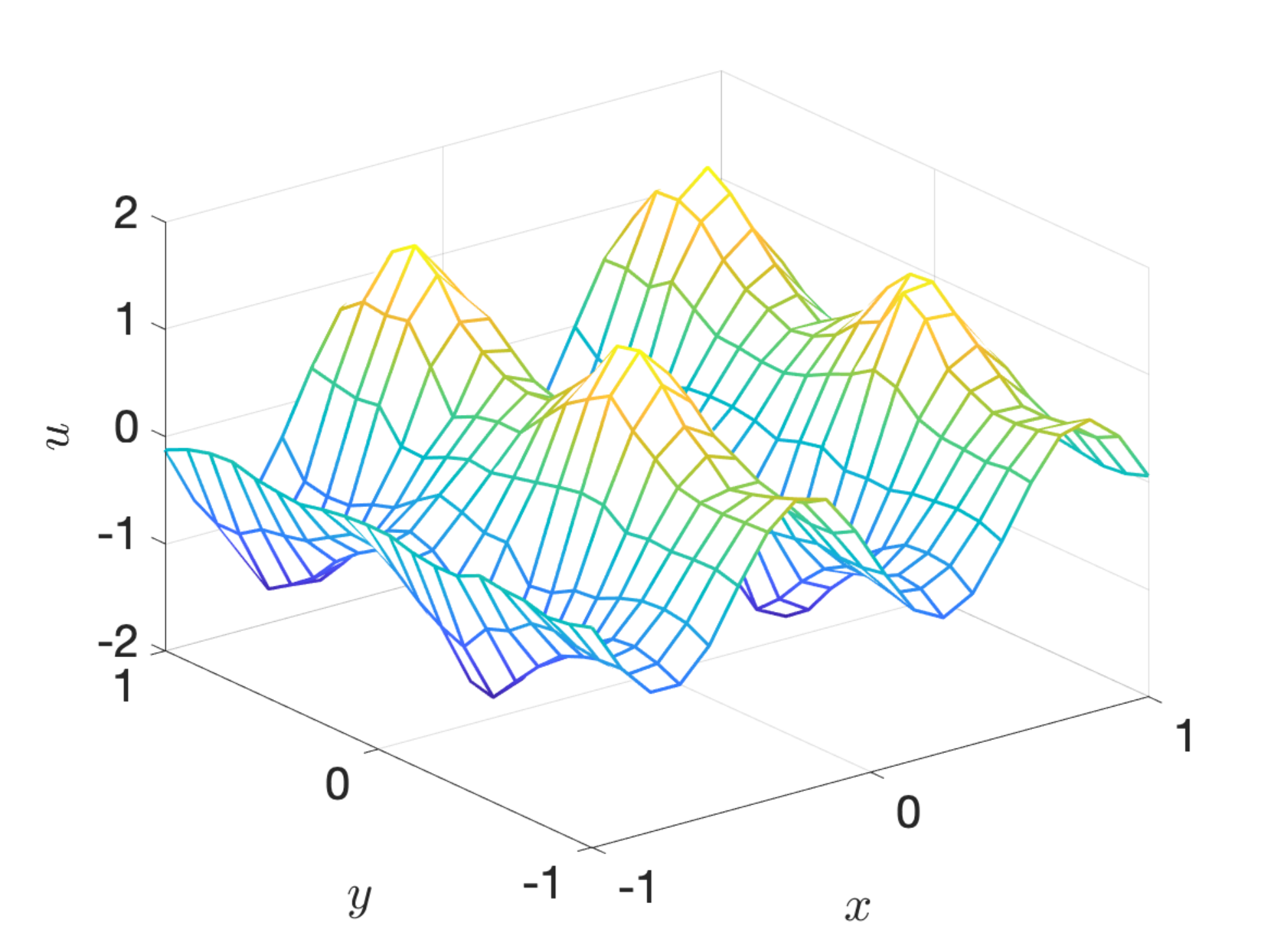}
    		\caption{Weak RBF, $t=400$}
    		\label{fig:linear_2d_cubic_random_weak_t=400} 
  	\end{subfigure}%
  	\caption{
  		Numerical results for the two-dimensional linear IVP with periodic BCs. 
		The cubic kernel with $N = 400$ random (uniformly distributed) nodes was used.   	
	}
  	\label{fig:linear_2d_cubic_random}
\end{figure} 

\begin{figure}[!htb]
  	\centering
 	\begin{subfigure}[b]{0.33\textwidth}
    		\includegraphics[width=\textwidth]{%
      		plots/linear_2d_ref}
    		\caption{Reference, $t=20$}
    		\label{fig:linear_2d_quintic_equid_ref_t=20} 
  	\end{subfigure}%
  	~ 
  	\begin{subfigure}[b]{0.33\textwidth}
    		\includegraphics[width=\textwidth]{%
      		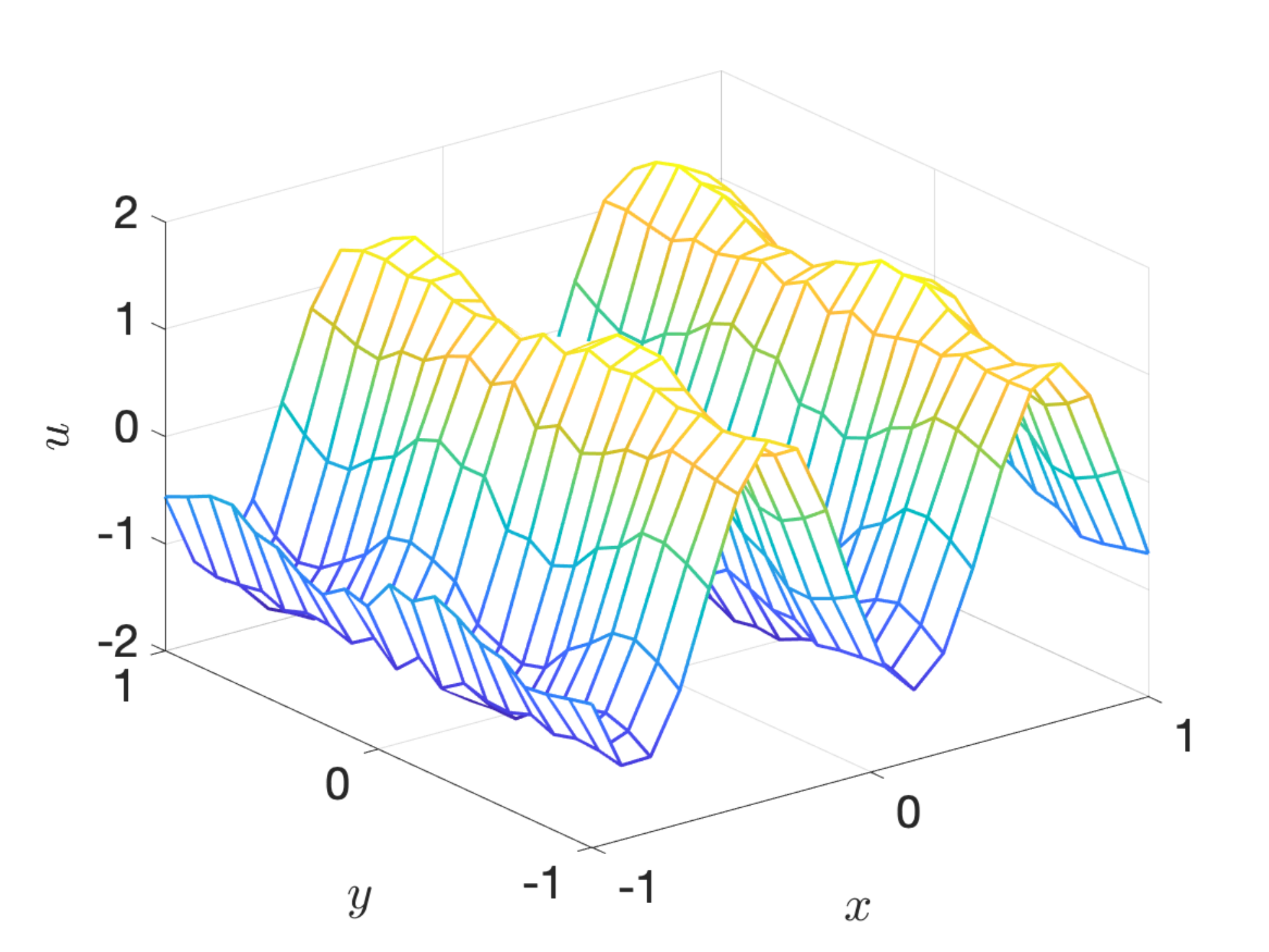}
    		\caption{{Standard} RBF, $t=20$}
    		\label{fig:linear_2d_quintic_equid_usual_t=20} 
  	\end{subfigure}%
	~ 
  	\begin{subfigure}[b]{0.33\textwidth}
    		\includegraphics[width=\textwidth]{%
      		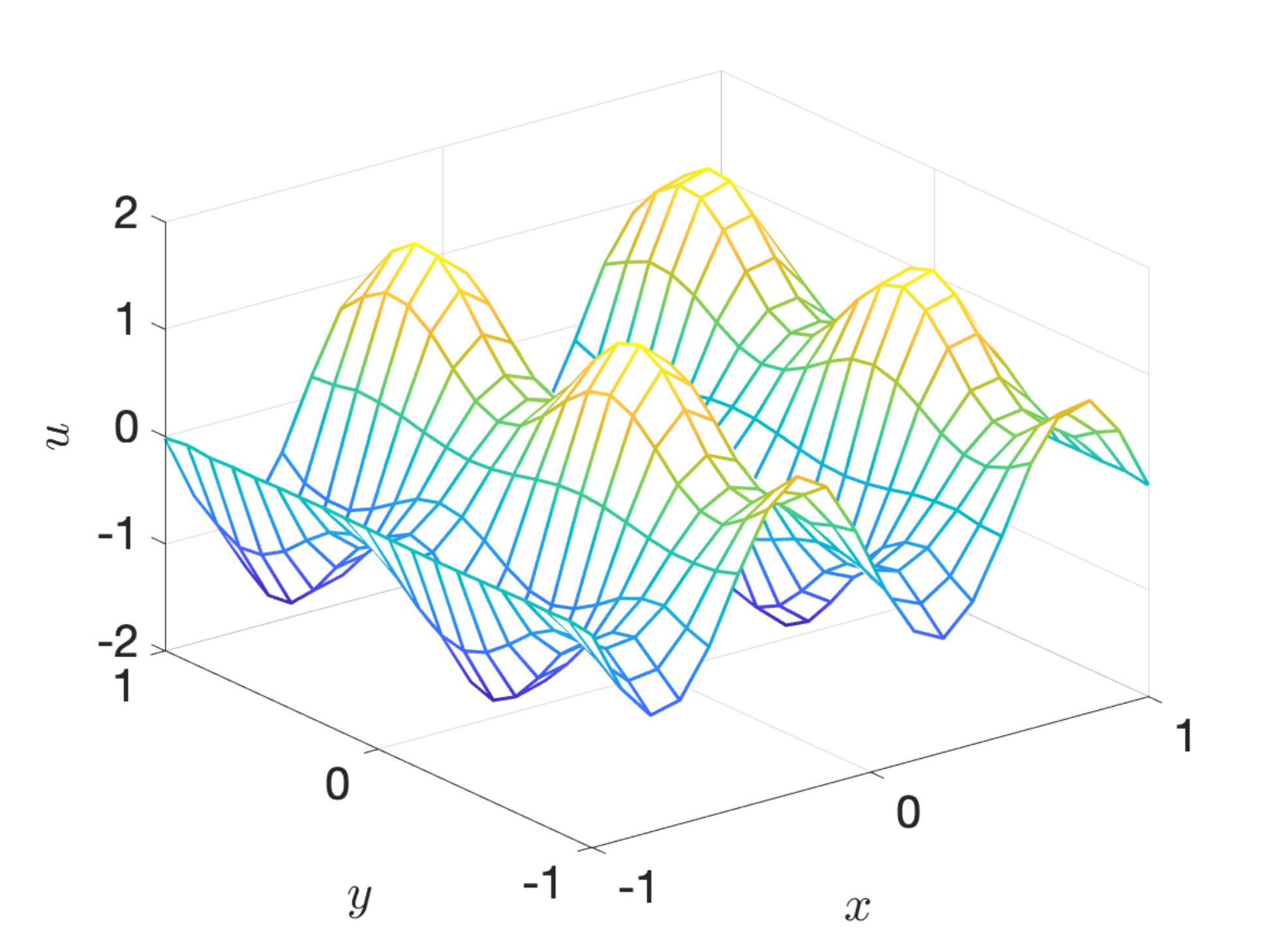}
    		\caption{Weak RBF, $t=20$}
    		\label{fig:linear_2d_quintic_equid_weak_t=20} 
  	\end{subfigure}%
	\\ 
  	\begin{subfigure}[b]{0.33\textwidth}
    		\includegraphics[width=\textwidth]{%
      		plots/linear_2d_ref}
    		\caption{Reference, $t=1600$}
    		\label{fig:linear_2d_quintic_equid_ref_t=1600} 
  	\end{subfigure}%
  	~ 
  	\begin{subfigure}[b]{0.33\textwidth}
    		\includegraphics[width=\textwidth]{%
      		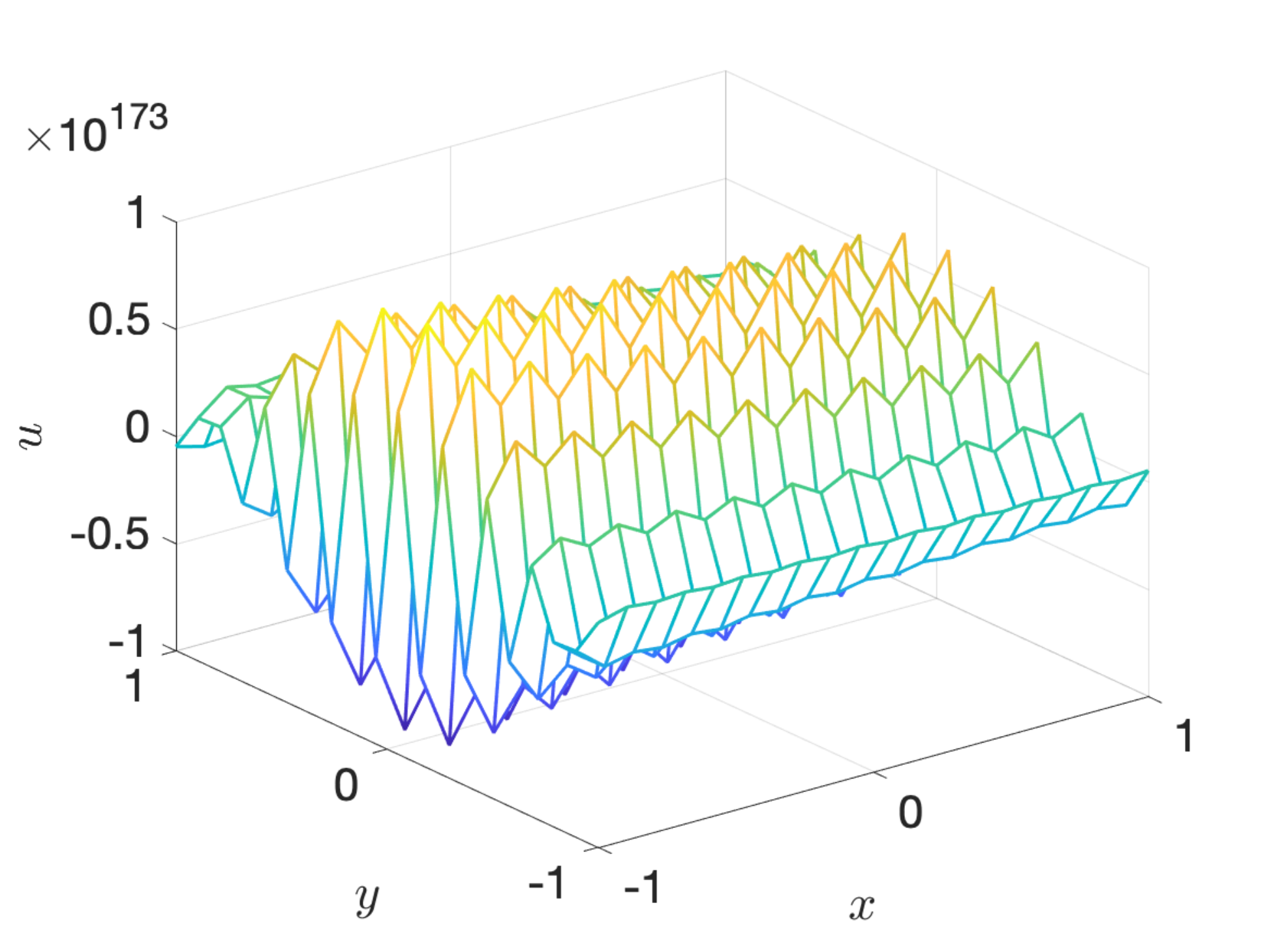}
    		\caption{{Standard} RBF, $t=1600$}
    		\label{fig:linear_2d_quintic_equid_usual_t=1600} 
  	\end{subfigure}%
	~ 
  	\begin{subfigure}[b]{0.33\textwidth}
    		\includegraphics[width=\textwidth]{%
      		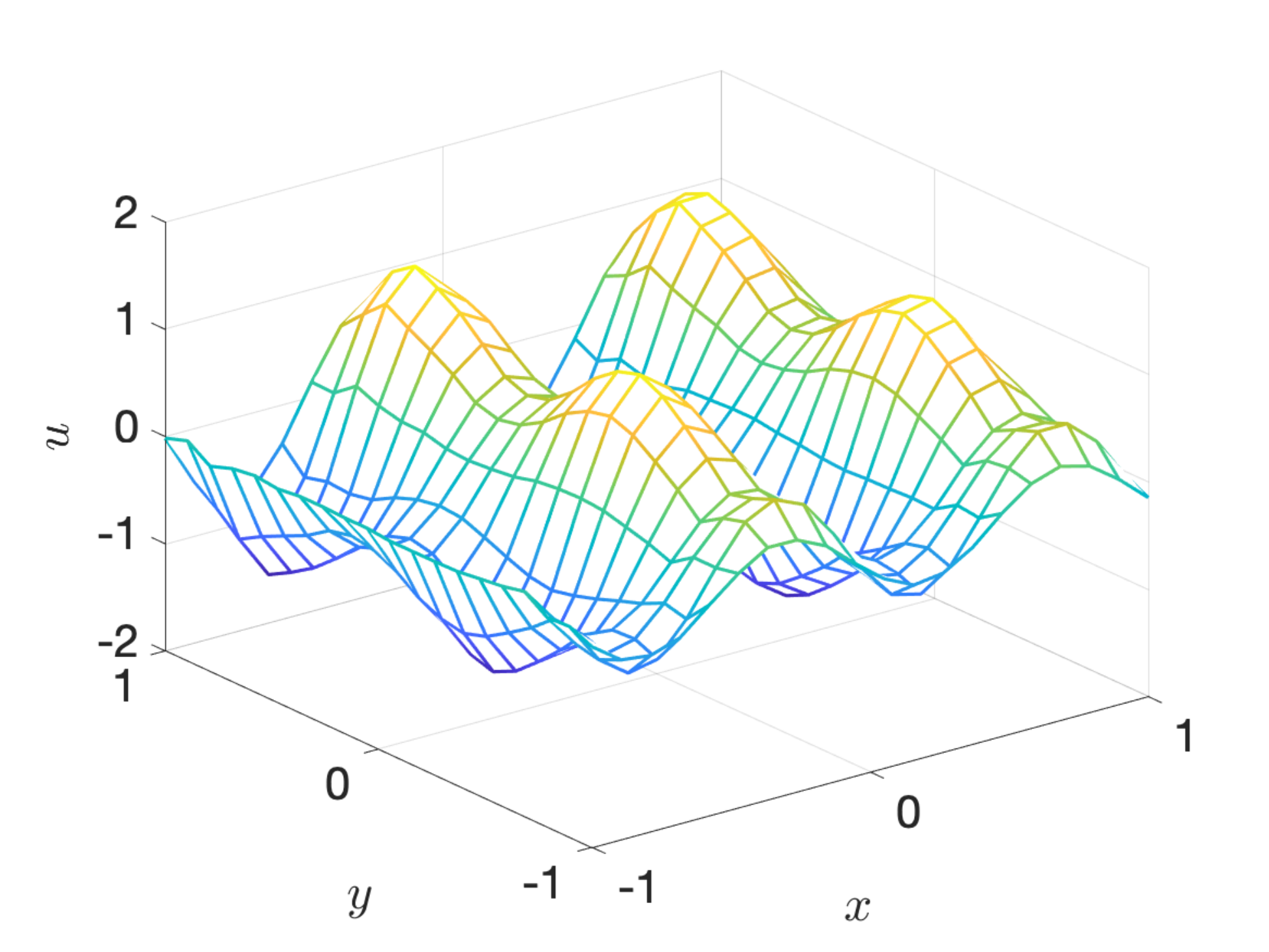}
    		\caption{Weak RBF, $t=1600$}
    		\label{fig:linear_2d_quintic_equid_weak_t=1600} 
  	\end{subfigure}%
  	\caption{
  		Numerical results for the two-dimensional linear IVP with periodic BCs. 
		The quintic kernel with $N = 400$ equidistant nodes was used.   	
	}
  	\label{fig:linear_2d_quintic_equid}
\end{figure} 

\begin{figure}[!htb]
  	\centering
 	\begin{subfigure}[b]{0.33\textwidth}
    		\includegraphics[width=\textwidth]{%
      		plots/linear_2d_ref_random}
    		\caption{Reference, $t=20$}
    		\label{fig:linear_2d_quintic_random_ref_t=20} 
  	\end{subfigure}%
  	~ 
  	\begin{subfigure}[b]{0.33\textwidth}
    		\includegraphics[width=\textwidth]{%
      		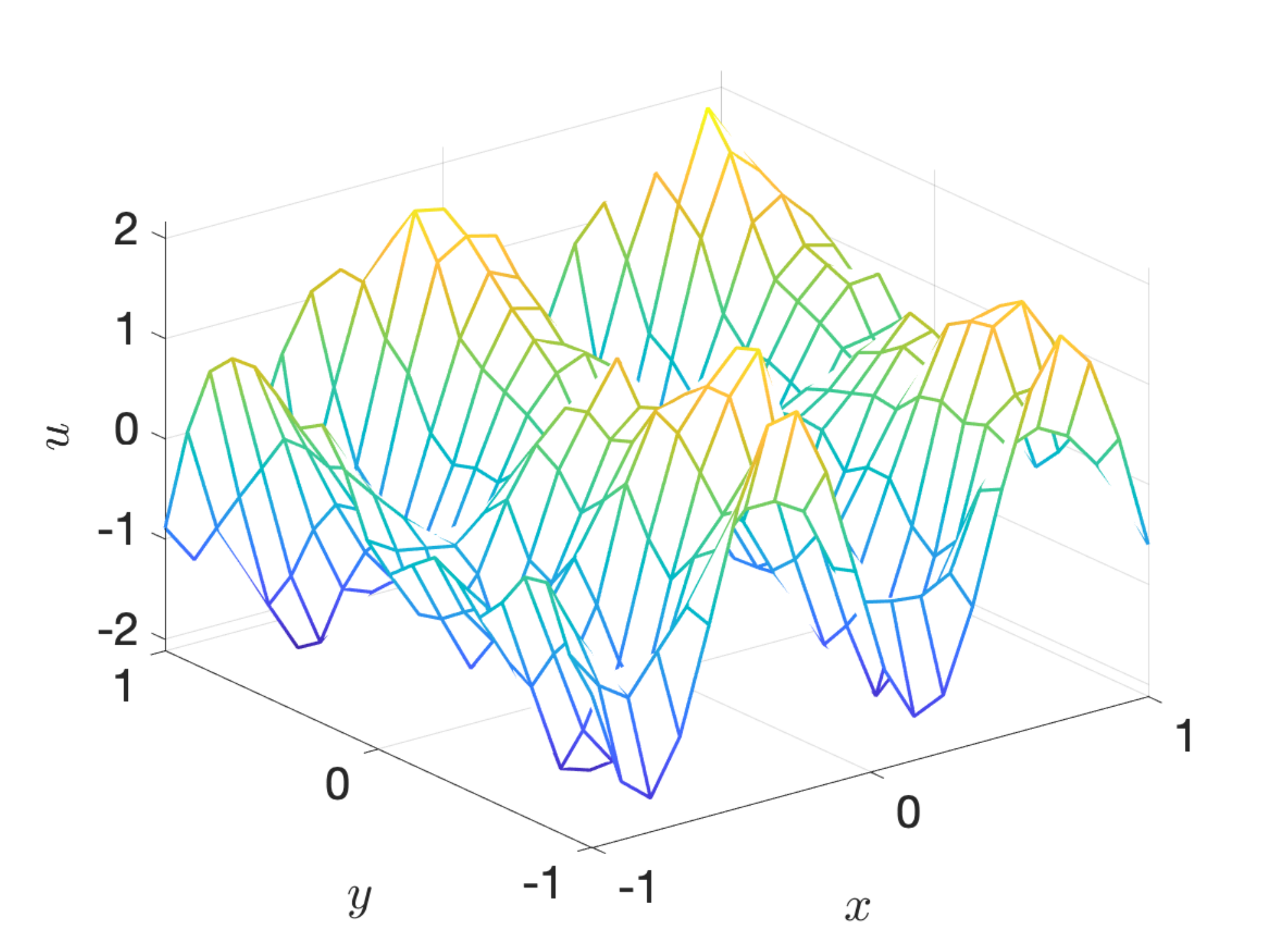}
    		\caption{{Standard} RBF, $t=20$}
    		\label{fig:linear_2d_quintic_random_usual_t=20} 
  	\end{subfigure}%
	~ 
  	\begin{subfigure}[b]{0.33\textwidth}
    		\includegraphics[width=\textwidth]{%
      		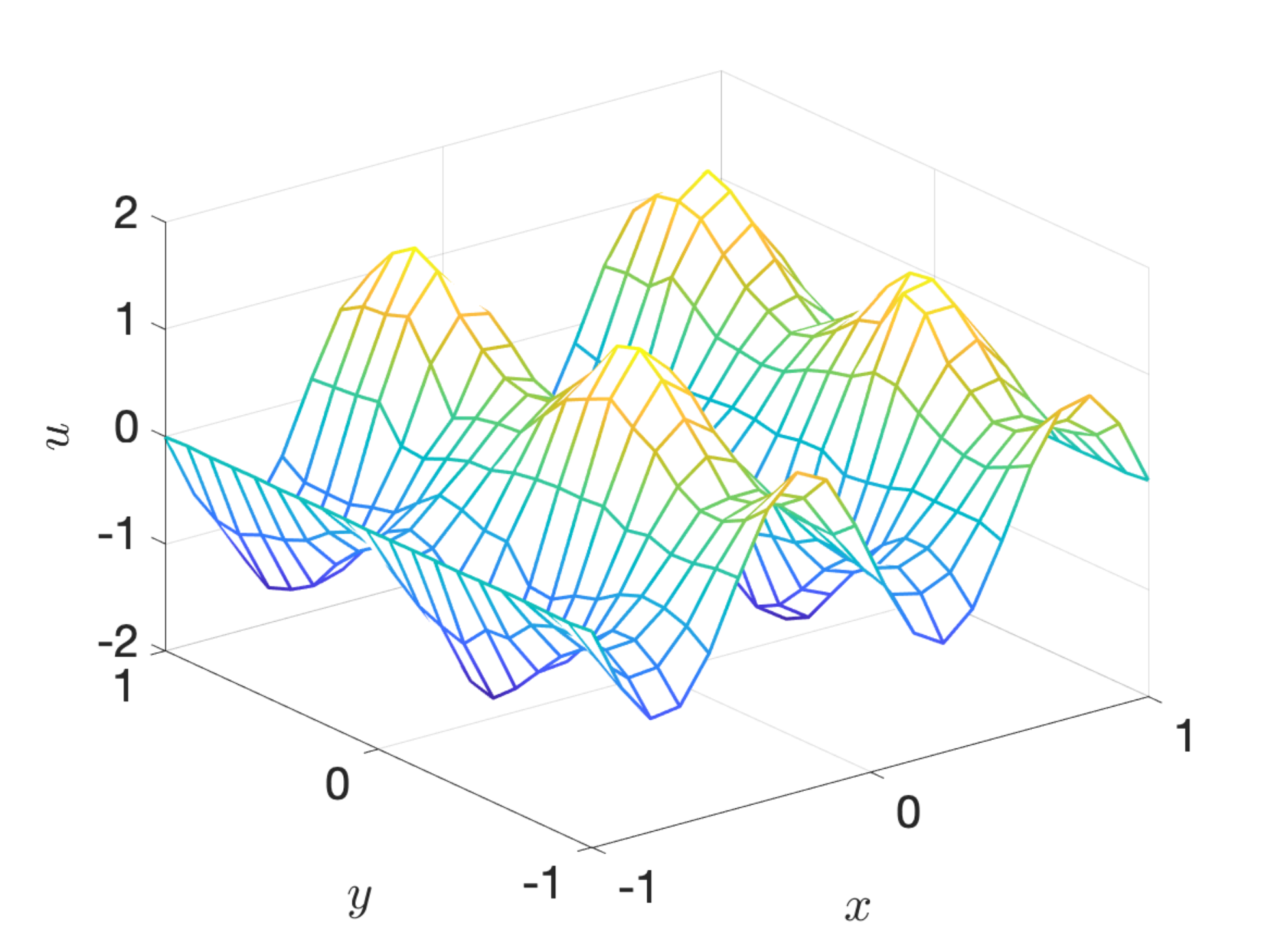}
    		\caption{Weak RBF, $t=20$}
    		\label{fig:linear_2d_quintic_random_weak_t=20} 
  	\end{subfigure}%
	\\ 
  	\begin{subfigure}[b]{0.33\textwidth}
    		\includegraphics[width=\textwidth]{%
      		plots/linear_2d_ref_random}
    		\caption{Reference, $t=1600$}
    		\label{fig:linear_2d_quintic_random_ref_t=1600} 
  	\end{subfigure}%
  	~ 
  	\begin{subfigure}[b]{0.33\textwidth}
    		\includegraphics[width=\textwidth]{%
      		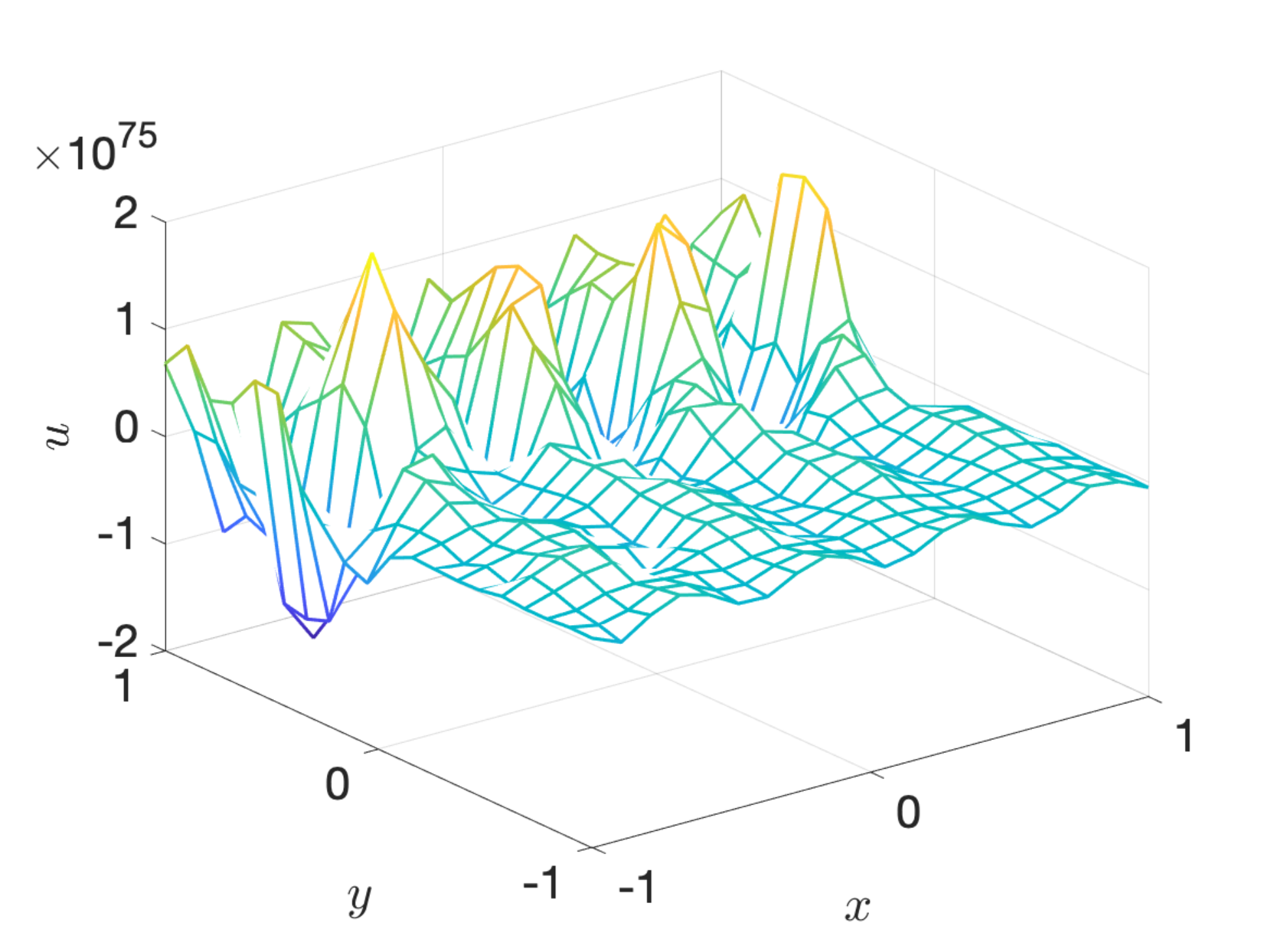}
    		\caption{{Standard} RBF, $t=1600$}
    		\label{fig:linear_2d_quintic_random_usual_t=1600} 
  	\end{subfigure}%
	~ 
  	\begin{subfigure}[b]{0.33\textwidth}
    		\includegraphics[width=\textwidth]{%
      		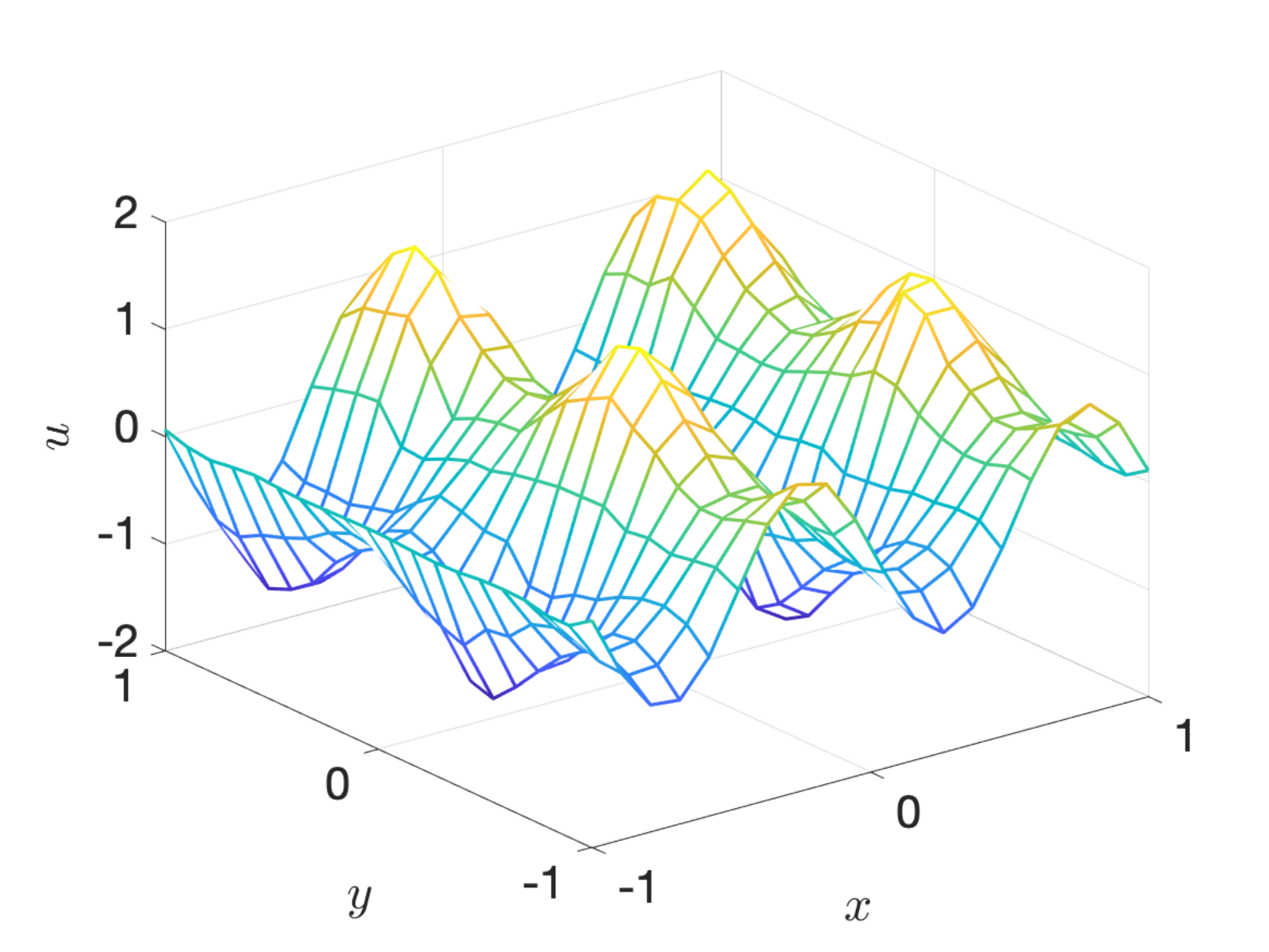}
    		\caption{Weak RBF, $t=1600$}
    		\label{fig:linear_2d_quintic_random_weak_t=1600} 
  	\end{subfigure}%
  	\caption{
  		Numerical results for the two-dimensional linear IVP with periodic BCs. 
		The quintic kernel with $N = 400$ random (uniformly distributed) nodes was used.   	
	}
  	\label{fig:linear_2d_quintic_random}
\end{figure} 

Figures \ref{fig:linear_2d_cubic_equid} and \ref{fig:linear_2d_cubic_random} respectively illustrate the results for the cubic kernel and $N=400$ equidistant and uniformly distributed nodes. 
Figures \ref{fig:linear_2d_quintic_equid} and \ref{fig:linear_2d_quintic_random} present the same result for the quintic kernel. 
In all computations the `exact' integration, performed by MATLAB's \emph{integral2}, was too cost prohibitive.   We therefore replaced it by a tensor product based two-dimensional trapezoidal rule (using $J=1000$ quadrature points in one dimension). 
Based on the results in \S \ref{subsub:num_int}, we believe that a significantly smaller number of quadrature points would have been sufficient.
We used $P = 1$ for the weak RBF method.

The standard RBF method blew up after comparatively small times in all test cases.  By contrast,
the weak  RBF method produced highly accurate results even for long time simulations.  This was true for both equidistant and nonequidistant points. 
After long simulation times, the weak RBF method is seen to decrease in accuracy, which appears to be unrelated to instability. 
Rather it seems that dissipation introduced by the numerical (full-upwind) fluxes is blurring the solution over long times. 
The weak RBF method remained stable for computations up to at least $t=1600$, at which we point we concluded our experiment. 
Future investigation will include using an energy conserving flux, such as a central flux, to determine if this will alleviate the long term dissipation.

\section{Concluding Remarks}
\label{sec:summary} 

In this work we investigated the conservation and energy stability properties of RBF methods. 
In the process we demonstrated that traditional RBF methods based on the strong form of hyperbolic 
CLs, including strong enforcement of BCs, violate these properties and might therefore produce physically unreasonable solutions. 
As an alternative we proposed a weak enforcement of BCs by building RBF schemes based on the weak form of the hyperbolic CL. 
We proved that the resulting methods are conservative assuming that (at least) constants are included in the RBF space. 
Furthermore, these methods were also shown to be energy stable when appropriate 
numerical (E-) fluxes are included in the discretization. 
In case of the weak RBF collocation method this was shown for linear advection when appropriate numerical (E-) fluxes are included in the discretization. 
Thus, the weak RBF methods are able to provide numerical solutions with physically reasonable mass 
and energy profiles. 
A drawback of this approach might be potentially ill-conditioned mass matrices, which arise from the weak form of the CL, \cite[Chapter 7.2.7]{glaubitz2020shock}. 
This may be overcome by choosing sufficiently large shape parameters.  
For more sophisticated applications requiring other kernels, it might be better to use orthonormal basis functions instead.

Future work will focus on the application of the proposed weak RBF method to nonlinear problems and, in particular, on the adaptation of different methods 
\cite{tadmor1990shock,krivodonova2004shock,hesthaven2008filtering,klockner2011viscous,ranocha2018stability,glaubitz2018application,glaubitz2019smooth,glaubitz2019high}
from DG and related methods to further stabilize the weak RBF method in the presence of (shock) discontinuities. 
Moreover, in a forthcoming work \cite{glaubitz2021towards}, the weak enforcement of BCs was also investigated in the context of RBF methods for linear advection problems based on their strong form. 
Finally, in addition to the energy stability analysis provided here, it would be useful to perform a (linear) eigenvalue stability analysis.

\section*{Acknowledgements}
The authors would like to thank Simon-Christian Klein for helpful advice.

\bibliographystyle{siamplain}
\bibliography{literature}

\begin{thebibliography}{100}

\bibitem{abgrall2019high}
{\sc R.~Abgrall, P.~Bacigaluppi, and S.~Tokareva}, {\em High-order residual
  distribution scheme for the time-dependent {E}uler equations of fluid
  dynamics}, Computers \& Mathematics with Applications, 78 (2019),
  pp.~274--297.

\bibitem{abgrall2019analysis}
{\sc R.~Abgrall, J.~Nordstr{\"o}m, P.~{\"O}ffner, and S.~Tokareva}, {\em
  Analysis of the {SBP-SAT} stabilization for finite element methods part i:
  Linear problems}, Journal of Scientific Computing, 85 (2020), pp.~1--29.

\bibitem{abgrall2019analysis_2}
{\sc R.~Abgrall, J.~Nordstr{\"o}m, P.~{\"O}ffner, and S.~Tokareva}, {\em
  Analysis of the {SBP-SAT} stabilization for finite element methods part ii:
  Entropy stability}, arXiv:1912.08390,  (2020).
\newblock Accepted in {C}ommunications on {A}pplied {M}athematics and
  {C}omputation.

\bibitem{buhmann2000radial}
{\sc M.~D. Buhmann}, {\em Radial basis functions}, Acta Numerica, 9 (2000),
  pp.~1--38.

\bibitem{buhmann2003radial}
{\sc M.~D. Buhmann}, {\em Radial Basis Functions: Theory and Implementations},
  vol.~12, Cambridge University Press, 2003.

\bibitem{caflisch1998monte}
{\sc R.~E. Caflisch}, {\em Monte {C}arlo and quasi-{M}onte {C}arlo methods},
  Acta Numerica, 1998 (1998), pp.~1--49.

\bibitem{canuto2006spectral}
{\sc C.~Canuto, M.~Y. Hussaini, A.~Quarteroni, and T.~A. Zang}, {\em Spectral
  Methods}, Springer, 2006.

\bibitem{canuto1982error}
{\sc C.~Canuto and A.~Quarteroni}, {\em Error estimates for spectral and
  pseudospectral approximations of hyperbolic equations}, SIAM Journal on
  Numerical Analysis, 19 (1982), pp.~629--642.

\bibitem{chen2017entropy}
{\sc T.~Chen and C.-W. Shu}, {\em Entropy stable high order discontinuous
  {G}alerkin methods with suitable quadrature rules for hyperbolic conservation
  laws}, Journal of Computational Physics, 345 (2017), pp.~427--461.

\bibitem{cheng2014positivity}
{\sc J.~Cheng and C.-W. Shu}, {\em Positivity-preserving {L}agrangian scheme
  for multi-material compressible flow}, Journal of Computational Physics, 257
  (2014), pp.~143--168.

\bibitem{cockburn1990runge}
{\sc B.~Cockburn, S.~Hou, and C.-W. Shu}, {\em The {R}unge--{K}utta local
  projection discontinuous {G}alerkin finite element method for conservation
  laws. iv. {T}he multidimensional case}, Mathematics of Computation, 54
  (1990), pp.~545--581.

\bibitem{cockburn1989tvb2}
{\sc B.~Cockburn, S.-Y. Lin, and C.-W. Shu}, {\em {TVB} {R}unge--{K}utta local
  projection discontinuous {G}alerkin finite element method for conservation
  laws iii: {O}ne-dimensional systems}, Journal of Computational Physics, 84
  (1989), pp.~90--113.

\bibitem{cockburn1989tvb}
{\sc B.~Cockburn and C.-W. Shu}, {\em {TVB} {R}unge--{K}utta local projection
  discontinuous {G}alerkin finite element method for conservation laws. ii.
  {G}eneral framework}, Mathematics of Computation, 52 (1989), pp.~411--435.

\bibitem{cockburn1991runge}
{\sc B.~Cockburn and C.-W. Shu}, {\em The {R}unge--{K}utta local projection
  ${P}^1$-discontinuous-{G}alerkin finite element method for scalar
  conservation laws}, ESAIM: Mathematical Modelling and Numerical Analysis, 25
  (1991), pp.~337--361.

\bibitem{cockburn1998runge}
{\sc B.~Cockburn and C.-W. Shu}, {\em The {R}unge--{K}utta discontinuous
  {G}alerkin method for conservation laws v: {M}ultidimensional systems},
  Journal of Computational Physics, 141 (1998), pp.~199--224.

\bibitem{dafermos2005hyperbolic}
{\sc C.~M. Dafermos}, {\em Hyperbolic Conservation Laws in Continuum Physics},
  vol.~3, Springer, 2005.

\bibitem{davis2007methods}
{\sc P.~J. Davis and P.~Rabinowitz}, {\em Methods of Numerical Integration},
  Courier Corporation, 2007.

\bibitem{dick2013high}
{\sc J.~Dick, F.~Y. Kuo, and I.~H. Sloan}, {\em High-dimensional integration:
  the quasi-{M}onte {C}arlo way}, Acta Numerica, 22 (2013), p.~133.

\bibitem{don2016hybrid}
{\sc W.-S. Don, Z.~Gao, P.~Li, and X.~Wen}, {\em Hybrid compact-{WENO} finite
  difference scheme with conjugate {F}ourier shock detection algorithm for
  hyperbolic conservation laws}, SIAM Journal on Scientific Computing, 38
  (2016), pp.~A691--A711.

\bibitem{duchon1977splines}
{\sc J.~Duchon}, {\em Splines minimizing rotation-invariant semi-norms in
  sobolev spaces}, in Constructive theory of functions of several variables,
  Springer, 1977, pp.~85--100.

\bibitem{engels1980numerical}
{\sc H.~Engels}, {\em Numerical Quadrature and Cubature}, Academic Press, 1980.

\bibitem{fasshauer1996solving}
{\sc G.~E. Fasshauer}, {\em Solving partial differential equations by
  collocation with radial basis functions}, in Proceedings of Chamonix,
  vol.~1997, Vanderbilt University Press Nashville, TN, 1996, pp.~1--8.

\bibitem{fasshauer2007meshfree}
{\sc G.~E. Fasshauer}, {\em Meshfree Approximation Methods with MATLAB},
  vol.~6, World Scientific, 2007.

\bibitem{fernandez2014review}
{\sc D.~C. D.~R. Fern{\'a}ndez, J.~E. Hicken, and D.~W. Zingg}, {\em Review of
  summation-by-parts operators with simultaneous approximation terms for the
  numerical solution of partial differential equations}, Computers \& Fluids,
  95 (2014), pp.~171--196.

\bibitem{flyer2016enhancing}
{\sc N.~Flyer, G.~A. Barnett, and L.~J. Wicker}, {\em Enhancing finite
  differences with radial basis functions: experiments on the
  {N}avier--{S}tokes equations}, Journal of Computational Physics, 316 (2016),
  pp.~39--62.

\bibitem{flyer2016role}
{\sc N.~Flyer, B.~Fornberg, V.~Bayona, and G.~A. Barnett}, {\em On the role of
  polynomials in {RBF-FD} approximations: I. {I}nterpolation and accuracy},
  Journal of Computational Physics, 321 (2016), pp.~21--38.

\bibitem{fornberg2002observations}
{\sc B.~Fornberg, T.~A. Driscoll, G.~Wright, and R.~Charles}, {\em Observations
  on the behavior of radial basis function approximations near boundaries},
  Computers \& Mathematics with Applications, 43 (2002), pp.~473--490.

\bibitem{fornberg2015primer}
{\sc B.~Fornberg and N.~Flyer}, {\em A Primer on Radial Basis Functions With
  Applications to the Geosciences}, SIAM, 2015.

\bibitem{fornberg2007runge}
{\sc B.~Fornberg and J.~Zuev}, {\em The {R}unge phenomenon and spatially
  variable shape parameters in {RBF} interpolation}, Computers \& Mathematics
  with Applications, 54 (2007), pp.~379--398.

\bibitem{funaro1988new}
{\sc D.~Funaro and D.~Gottlieb}, {\em A new method of imposing boundary
  conditions in pseudospectral approximations of hyperbolic equations},
  Mathematics of Computation, 51 (1988), pp.~599--613.

\bibitem{funaro1991convergence}
{\sc D.~Funaro and D.~Gottlieb}, {\em Convergence results for pseudospectral
  approximations of hyperbolic systems by a penalty-type boundary treatment},
  Mathematics of Computation, 57 (1991), pp.~585--596.

\bibitem{gassner2013skew}
{\sc G.~J. Gassner}, {\em A skew-symmetric discontinuous {G}alerkin spectral
  element discretization and its relation to {SBP}-{SAT} finite difference
  methods}, SIAM Journal on Scientific Computing, 35 (2013), pp.~A1233--A1253.

\bibitem{gelb2019numerical}
{\sc A.~Gelb, X.~Hou, and Q.~Li}, {\em Numerical analysis for conservation laws
  using $\ell_1$ minimization}, Journal of Scientific Computing, 81 (2019),
  pp.~1240--1265.

\bibitem{glaubitz2019shock}
{\sc J.~Glaubitz}, {\em Shock capturing by {B}ernstein polynomials for scalar
  conservation laws}, Applied Mathematics and Computation, 363 (2019),
  p.~124593, \url{https://doi.org/10.1016/j.amc.2019.124593}.

\bibitem{glaubitz2020constructing}
{\sc J.~Glaubitz}, {\em Constructing positive interpolatory cubature formulas},
  arXiv preprint arXiv:2009.11981,  (2020).
\newblock Submitted.

\bibitem{glaubitz2020weakRBFcode}
{\sc J.~Glaubitz}, {\em jglaubitz/weak{RBF}}, 2020,
  \url{https://doi.org/10.5281/zenodo.4310328},
  \url{https://doi.org/10.5281/zenodo.4310328}.

\bibitem{glaubitz2020shock}
{\sc J.~Glaubitz}, {\em Shock Capturing and High-Order Methods for Hyperbolic
  Conservation Laws}, Logos Verlag Berlin, 2020,
  \url{https://doi.org/10.30819/5084}.

\bibitem{glaubitz2020stableCF}
{\sc J.~Glaubitz}, {\em Stable high-order cubature formulas for experimental
  data}, arXiv preprint arXiv:2009.03452,  (2020).
\newblock Submitted.

\bibitem{glaubitz2020stable}
{\sc J.~Glaubitz}, {\em Stable high order quadrature rules for scattered data
  and general weight functions}, SIAM Journal on Numerical Analysis, 58 (2020),
  pp.~2144--2164.

\bibitem{glaubitz2019high}
{\sc J.~Glaubitz and A.~Gelb}, {\em High order edge sensors with $\ell^1$
  regularization for enhanced discontinuous {G}alerkin methods}, SIAM Journal
  on Scientific Computing, 41 (2019), pp.~A1304--A1330.

\bibitem{glaubitz2021towards}
{\sc J.~Glaubitz, E.~L. M{\'e}l{\'e}do, and P.~{\"O}ffner}, {\em Towards stable
  radial basis function methods for linear advection problems}, arXiv preprint
  arXiv:2101.09623,  (2021).
\newblock To appear.

\bibitem{glaubitz2019smooth}
{\sc J.~Glaubitz, A.~Nogueira, J.~Almeida, R.~Cant{\~a}o, and C.~Silva}, {\em
  Smooth and compactly supported viscous sub-cell shock capturing for
  discontinuous {G}alerkin methods}, Journal of Scientific Computing, 79
  (2019), pp.~249--272.

\bibitem{glaubitz2020stableDG}
{\sc J.~Glaubitz and P.~{\"O}ffner}, {\em Stable discretisations of high-order
  discontinuous galerkin methods on equidistant and scattered points}, Applied
  Numerical Mathematics, 151 (2020), pp.~98--118,
  \url{https://doi.org/10.1016/j.apnum.2019.12.020}.

\bibitem{glaubitz2016artificial}
{\sc J.~Glaubitz, P.~{\"O}ffner, H.~Ranocha, and T.~Sonar}, {\em Artificial
  viscosity for correction procedure via reconstruction using
  summation-by-parts operators}, in XVI International Conference on Hyperbolic
  Problems: Theory, Numerics, Applications, Springer, 2016, pp.~363--375.

\bibitem{glaubitz2018application}
{\sc J.~Glaubitz, P.~{\"O}ffner, and T.~Sonar}, {\em Application of modal
  filtering to a spectral difference method}, Mathematics of Computation, 87
  (2018), pp.~175--207.

\bibitem{gottlieb2001spectral}
{\sc D.~Gottlieb and J.~S. Hesthaven}, {\em Spectral methods for hyperbolic
  problems}, Journal of Computational and Applied Mathematics, 128 (2001),
  pp.~83--131.

\bibitem{gottlieb1998total}
{\sc S.~Gottlieb and C.-W. Shu}, {\em Total variation diminishing
  {R}unge--{K}utta schemes}, Mathematics of Aomputation, 67 (1998), pp.~73--85.

\bibitem{gottlieb2001strong}
{\sc S.~Gottlieb, C.-W. Shu, and E.~Tadmor}, {\em Strong stability-preserving
  high-order time discretization methods}, SIAM Review, 43 (2001), pp.~89--112.

\bibitem{gustafsson2007high}
{\sc B.~Gustafsson}, {\em High Order Difference Methods for Time Dependent
  PDE}, vol.~38, Springer Science \& Business Media, 2007.

\bibitem{gustafsson1995time}
{\sc B.~Gustafsson, H.-O. Kreiss, and J.~Oliger}, {\em Time Dependent Problems
  and Difference Methods}, vol.~24, John Wiley \& Sons, 1995.

\bibitem{haber1970numerical}
{\sc S.~Haber}, {\em Numerical evaluation of multiple integrals}, SIAM Review,
  12 (1970), pp.~481--526.

\bibitem{hesthaven2008filtering}
{\sc J.~Hesthaven and R.~Kirby}, {\em Filtering in {L}egendre spectral
  methods}, Mathematics of Computation, 77 (2008), pp.~1425--1452.

\bibitem{hesthaven2000spectral}
{\sc J.~S. Hesthaven}, {\em Spectral penalty methods}, Applied Numerical
  Mathematics, 33 (2000), pp.~23--41.

\bibitem{hesthaven2019entropy}
{\sc J.~S. Hesthaven and F.~M{\"o}nkeberg}, {\em Entropy stable essentially
  nonoscillatory methods based on {RBF} reconstruction}, ESAIM: Mathematical
  Modelling and Numerical Analysis, 53 (2019), pp.~925--958.

\bibitem{hesthaven2020two}
{\sc J.~S. Hesthaven and F.~M{\"o}nkeberg}, {\em Two-dimensional {RBF-ENO}
  method on unstructured grids}, Journal of Scientific Computing, 82 (2020),
  pp.~1--24.

\bibitem{hesthaven2007nodal}
{\sc J.~S. Hesthaven and T.~Warburton}, {\em Nodal Discontinuous {G}alerkin
  Methods: Algorithms, Analysis, and Applications}, Springer Science \&
  Business Media, 2007.

\bibitem{hon1998efficient}
{\sc Y.~Hon and X.~Mao}, {\em An efficient numerical scheme for {B}urgers'
  equation}, Applied Mathematics and Computation, 95 (1998), pp.~37--50.

\bibitem{huybrechs2009stable}
{\sc D.~Huybrechs}, {\em Stable high-order quadrature rules with equidistant
  points}, Journal of Computational and Applied Mathematics, 231 (2009),
  pp.~933--947.

\bibitem{huynh2007flux}
{\sc H.~T. Huynh}, {\em A flux reconstruction approach to high-order schemes
  including discontinuous {G}alerkin methods}, AIAA paper, 4079 (2007),
  p.~2007.

\bibitem{iske2003radial}
{\sc A.~Iske}, {\em Radial basis functions: basics, advanced topics and
  meshfree methods for transport problems}, Rend. Sem. Mat. Univ. Pol. Torino,
  61 (2003), pp.~247--285.

\bibitem{iske2011scattered}
{\sc A.~Iske}, {\em Scattered data approximation by positive definite kernel
  functions}, Rend. Sem. Mat. Univ. Pol. Torino, 69 (2011), pp.~217--246.

\bibitem{iske2020ten}
{\sc A.~Iske}, {\em Ten good reasons for using polyharmonic spline
  reconstruction in particle fluid flow simulations}, in Continuum Mechanics,
  Applied Mathematics and Scientific Computing: Godunov's Legacy, Springer,
  2020, pp.~193--199.

\bibitem{iske1996structure}
{\sc A.~Iske and T.~Sonar}, {\em On the structure of function spaces in optimal
  recovery of point functionals for {ENO}-schemes by radial basis functions},
  Numerische Mathematik, 74 (1996), pp.~177--201.

\bibitem{jameson2012non}
{\sc A.~Jameson, P.~E. Vincent, and P.~Castonguay}, {\em On the non-linear
  stability of flux reconstruction schemes}, Journal of Scientific Computing,
  50 (2012), pp.~434--445.

\bibitem{jiang1994cell}
{\sc G.~S. Jiang and C.-W. Shu}, {\em On a cell entropy inequality for
  discontinuous {G}alerkin methods}, Mathematics of Computation, 62 (1994),
  pp.~531--538.

\bibitem{kansa2000circumventing}
{\sc E.~Kansa and Y.~Hon}, {\em Circumventing the ill-conditioning problem with
  multiquadric radial basis functions: {A}pplications to elliptic partial
  differential equations}, Computers and Mathematics with Applications, 39
  (2000), pp.~123--138.

\bibitem{kansa1990multiquadrics}
{\sc E.~J. Kansa}, {\em Multiquadrics—a scattered data approximation scheme
  with applications to computational fluid-dynamics—ii {S}olutions to
  parabolic, hyperbolic and elliptic partial differential equations}, Computers
  \& Mathematics with Applications, 19 (1990), pp.~147--161.

\bibitem{ketcheson2008highly}
{\sc D.~I. Ketcheson}, {\em Highly efficient strong stability-preserving
  {R}unge--{K}utta methods with low-storage implementations}, SIAM Journal on
  Scientific Computing, 30 (2008), pp.~2113--2136.

\bibitem{klockner2011viscous}
{\sc A.~Kl{\"o}ckner, T.~Warburton, and J.~S. Hesthaven}, {\em Viscous shock
  capturing in a time-explicit discontinuous {G}alerkin method}, Mathematical
  Modelling of Natural Phenomena, 6 (2011), pp.~57--83.

\bibitem{kreiss1989initial}
{\sc H.-O. Kreiss and J.~Lorenz}, {\em Initial-Boundary Value Problems and the
  {N}avier--{S}tokes Equations}, vol.~47, SIAM, 1989.

\bibitem{krivodonova2004shock}
{\sc L.~Krivodonova, J.~Xin, J.-F. Remacle, N.~Chevaugeon, and J.~E. Flaherty},
  {\em Shock detection and limiting with discontinuous {G}alerkin methods for
  hyperbolic conservation laws}, Applied Numerical Mathematics, 48 (2004),
  pp.~323--338.

\bibitem{larsson2003numerical}
{\sc E.~Larsson and B.~Fornberg}, {\em A numerical study of some radial basis
  function based solution methods for elliptic pdes}, Computers \& Mathematics
  with Applications, 46 (2003), pp.~891--902.

\bibitem{lax1973hyperbolic}
{\sc P.~D. Lax}, {\em Hyperbolic Systems of Conservation Laws and the
  Mathematical Theory of Shock Waves}, SIAM, 1973.

\bibitem{leveque2002finite}
{\sc R.~J. LeVeque}, {\em Finite Volume Methods for Hyperbolic Problems},
  vol.~31, Cambridge University Press, 2002.

\bibitem{levy1998semidiscrete}
{\sc D.~Levy and E.~Tadmor}, {\em From semidiscrete to fully discrete:
  Stability of {R}unge--{K}utta schemes by the energy method}, SIAM Review, 40
  (1998), pp.~40--73.

\bibitem{madych1992miscellaneous}
{\sc W.~Madych}, {\em Miscellaneous error bounds for multiquadric and related
  interpolators}, Computers \& Mathematics with Applications, 24 (1992),
  pp.~121--138.

\bibitem{martel2016stability}
{\sc J.~M. Martel and R.~B. Platte}, {\em Stability of radial basis function
  methods for convection problems on the circle and sphere}, Journal of
  Scientific Computing, 69 (2016), pp.~487--505.

\bibitem{metropolis1949monte}
{\sc N.~Metropolis and S.~Ulam}, {\em The {M}onte {C}arlo method}, Journal of
  the American Statistical Association, 44 (1949), pp.~335--341.

\bibitem{niederreiter1992random}
{\sc H.~Niederreiter}, {\em Random Number Generation and Quasi-Monte Carlo
  Methods}, SIAM, 1992.

\bibitem{offner2018stability}
{\sc P.~{\"O}ffner, J.~Glaubitz, and H.~Ranocha}, {\em Stability of correction
  procedure via reconstruction with summation-by-parts operators for {B}urgers'
  equation using a polynomial chaos approach}, ESAIM: Mathematical Modelling
  and Numerical Analysis, 52 (2018), pp.~2215--2245.

\bibitem{glaubitz2019analysis}
{\sc P.~{\"O}ffner, J.~Glaubitz, and H.~Ranocha}, {\em Analysis of artificial
  dissipation of explicit and implicit time-integration methods}, International
  Journal of Numerical Analysis and Modeling, 17 (2020), pp.~332--349,
  \url{http://global-sci.org/intro/article_detail/ijnam/16862.html}.

\bibitem{osher1984riemann}
{\sc S.~Osher}, {\em Riemann solvers, the entropy condition, and difference},
  SIAM Journal on Numerical Analysis, 21 (1984), pp.~217--235.

\bibitem{platte2004computing}
{\sc R.~B. Platte and T.~A. Driscoll}, {\em Computing eigenmodes of elliptic
  operators using radial basis functions}, Computers and Mathematics with
  Applications, 48 (2004), pp.~561--576.

\bibitem{platte2005polynomials}
{\sc R.~B. Platte and T.~A. Driscoll}, {\em Polynomials and potential theory
  for {G}aussian radial basis function interpolation}, SIAM Journal on
  Numerical Analysis, 43 (2005), pp.~750--766.

\bibitem{platte2006eigenvalue}
{\sc R.~B. Platte and T.~A. Driscoll}, {\em Eigenvalue stability of radial
  basis function discretizations for time-dependent problems}, Computers \&
  Mathematics with Applications, 51 (2006), pp.~1251--1268.

\bibitem{powell1992theory}
{\sc M.~J. Powell}, {\em The theory of radial basis function approximation in
  1990}, Advances in numerical analysis,  (1992), pp.~105--210.

\bibitem{randall1992numerical}
{\sc J.~L. Randall}, {\em Numerical Methods for Conservation Laws}, Lectures in
  Mathematics ETH Z{\"u}rich, 1992.

\bibitem{ranocha2018stability}
{\sc H.~Ranocha, J.~Glaubitz, P.~{\"O}ffner, and T.~Sonar}, {\em Stability of
  artificial dissipation and modal filtering for flux reconstruction schemes
  using summation-by-parts operators}, Applied Numerical Mathematics, 128
  (2018), pp.~1--23.

\bibitem{ranocha2016summation}
{\sc H.~Ranocha, P.~{\"O}ffner, and T.~Sonar}, {\em Summation-by-parts
  operators for correction procedure via reconstruction}, Journal of
  Computational Physics, 311 (2016), pp.~299--328.

\bibitem{sarra2011numerical}
{\sc S.~A. Sarra, A.~R. Heryudono, and C.~Wang}, {\em A numerical study of a
  technique for shifting eigenvalues of radial basis function differentiation
  matrices}, Tech Report, MU-MTH-TR-2011-1,  (2011).

\bibitem{scarnati2018using}
{\sc T.~Scarnati, A.~Gelb, and R.~B. Platte}, {\em Using $\ell_1$
  regularization to improve numerical partial differential equation solvers},
  Journal of Scientific Computing, 75 (2018), pp.~225--252.

\bibitem{schaback1995creating}
{\sc R.~Schaback}, {\em Creating surfaces from scattered data using radial
  basis functions}, in Mathematical Methods for Curves and Surfaces, University
  Press, 1995, pp.~477--496.

\bibitem{schaback1995error}
{\sc R.~Schaback}, {\em Error estimates and condition numbers for radial basis
  function interpolation}, Advances in Computational Mathematics, 3 (1995),
  pp.~251--264.

\bibitem{schaback1995multivariate}
{\sc R.~Schaback}, {\em Multivariate interpolation and approximation by
  translates of a basis function}, Series In Approximations and Decompositions,
  6 (1995), pp.~491--514.

\bibitem{schaback2005multivariate}
{\sc R.~Schaback}, {\em Multivariate interpolation by polynomials and radial
  basis functions}, Constructive Approximation, 21 (2005), pp.~293--317.

\bibitem{shu1988total}
{\sc C.-W. Shu}, {\em Total-variation-diminishing time discretizations}, SIAM
  Journal on Scientific and Statistical Computing, 9 (1988), pp.~1073--1084.

\bibitem{stroud1971approximate}
{\sc A.~H. Stroud}, {\em Approximate Calculation of Multiple Integrals},
  Prentice-Hall, 1971.

\bibitem{svard2014review}
{\sc M.~Sv{\"a}rd and J.~Nordstr{\"o}m}, {\em Review of summation-by-parts
  schemes for initial--boundary-value problems}, Journal of Computational
  Physics, 268 (2014), pp.~17--38.

\bibitem{tadmor1990shock}
{\sc E.~Tadmor}, {\em Shock capturing by the spectral viscosity method},
  Computer Methods in Applied Mechanics and Engineering, 80 (1990),
  pp.~197--208.

\bibitem{tolstykh2000using}
{\sc A.~I. Tolstykh}, {\em On using {RBF}-based differencing formulas for
  unstructured and mixed structured-unstructured grid calculations}, in
  Proceedings of the 16th IMACS world congress, vol.~228, Lausanne, 2000,
  pp.~4606--4624.

\bibitem{toro2013riemann}
{\sc E.~F. Toro}, {\em Riemann Solvers and Numerical Methods for Fluid
  Dynamics: A Practical Introduction}, Springer Science \& Business Media,
  2013.

\bibitem{trefethen2017cubature}
{\sc L.~N. Trefethen}, {\em Cubature, approximation, and isotropy in the
  hypercube}, SIAM Review, 59 (2017), pp.~469--491.

\bibitem{vincent2011new}
{\sc P.~E. Vincent, P.~Castonguay, and A.~Jameson}, {\em A new class of
  high-order energy stable flux reconstruction schemes}, Journal of Scientific
  Computing, 47 (2011), pp.~50--72.

\bibitem{wendland2004scattered}
{\sc H.~Wendland}, {\em Scattered Data Approximation}, vol.~17, Cambridge
  University Press, 2004.

\bibitem{wilson1970discrete}
{\sc M.~W. Wilson}, {\em Discrete least squares and quadrature formulas},
  Mathematics of Computation, 24 (1970), pp.~271--282.

\bibitem{wilson1970necessary}
{\sc M.~W. Wilson}, {\em Necessary and sufficient conditions for equidistant
  quadrature formula}, SIAM Journal on Numerical Analysis, 7 (1970),
  pp.~134--141.

\end{thebibliography}

\end{document}